\documentclass{article}
\usepackage{graphicx}

\usepackage[letterpaper,top=2cm,bottom=2cm,left=3cm,right=3cm,marginparwidth=1.75cm]{geometry}
\usepackage{amsmath,amssymb,amsthm}
\usepackage[colorlinks=true, allcolors=blue, hypertexnames=false]{hyperref}
\usepackage{comment}
\usepackage{enumitem}
\usepackage{commath}
\usepackage{cleveref}
\usepackage{algorithm}
\usepackage{algpseudocode}
\usepackage{bm}
\usepackage{booktabs}
\usepackage{tikz}
\usetikzlibrary{arrows.meta,positioning,fit,backgrounds}

\usepackage{xcolor}

\newtheorem{thm}{Theorem}[section]
\newtheorem{proposition}{Proposition}[section]
\newtheorem{corollary}{Corollary}[section]
\newtheorem{lemma}{Lemma}[section]
\newtheorem{assumption}{Assumption}[section]
\theoremstyle{definition}
\newtheorem{definition}{Definition}[section]
\newtheorem{example}{Example}[section]
\theoremstyle{remark}
\newtheorem{rmk}{Remark}[section]
\theoremstyle{plain}

\numberwithin{equation}{section}
\numberwithin{figure}{section}
\numberwithin{table}{section}

\crefname{thm}{Theorem}{Theorems}
\Crefname{thm}{Theorem}{Theorems}

\crefname{proposition}{Proposition}{Propositions}
\Crefname{proposition}{Proposition}{Propositions}

\crefname{corollary}{Corollary}{Corollaries}
\Crefname{corollary}{Corollary}{Corollaries}

\crefname{definition}{Definition}{Definitions}
\Crefname{definition}{Definition}{Definitions}

\crefname{example}{Example}{Examples}
\Crefname{example}{Example}{Examples}

\crefname{lemma}{Lemma}{Lemmas}
\Crefname{lemma}{Lemma}{Lemmas}

\crefname{rmk}{Remark}{Remarks}
\Crefname{rmk}{Remark}{Remarks}

\crefname{assumption}{Assumption}{Assumptions}
\Crefname{assumption}{Assumption}{Assumptions}

\newcommand{\Law}{\mathrm{Law}}
\DeclareMathOperator{\supp}{supp}
\def\R{\mathbb{R}}
\def\cP{\mathcal{P}}
\def\cL{\mathcal{L}}

\newcommand{\KL}{\mathrm{KL}}
\newcommand{\bD}{\mathbf{D}}                 
\newcommand{\cX}{\mathcal{X}}                
\newcommand{\Pirho}{\Pi(\rho)}               
\newcommand{\Pizero}{\Pi(\rho;\pi_0)}        
\newcommand{\HY}{H^1(0,T;H)}                 
\newcommand{\cZ}{\mathcal{Z}_+}
\newcommand{\cG}{\mathcal{G}}
\newcommand{\cA}{\mathcal{A}}
\newcommand{\cM}{\mathcal{M}}
\newcommand{\cC}{\mathcal{C}}
\newcommand{\cF}{\mathcal{F}}
\newcommand{\cQ}{\mathcal{Q}}                
\newcommand{\LL}{\mathbb{L}}
\newcommand{\projr}{\operatorname{proj}_\rho} 
\newcommand{\prox}{\operatorname{prox}}      
\newcommand{\Expl}{\operatorname{Gap}_{\mathrm S}}
\newcommand{\GapM}{\operatorname{Gap}_{\mathrm M}}
\newcommand{\pair}[2]{\left\langle #1,#2\right\rangle}
\newcommand{\pp}[2]{\left\langle\!\left\langle #1,#2\right\rangle\!\right\rangle} 
\newcommand{\tn}[1]{\lVert #1\rVert_{H^1}}   
\newcommand{\tns}[1]{\lVert #1\rVert_{(H^1)^*}} 
\newcommand{\nsharp}[1]{\lVert #1\rVert_{\mathrm{TV},2}} 
\newcommand{\nG}[1]{\lVert #1\rVert_{\cG}}
\newcommand{\nH}[1]{\lVert #1\rVert_{H}}

\title{Primal--Dual Kullback--Leibler Mirror-Prox for Dynamic Mean Field Games}
\author{Erhan Bayraktar\thanks{Department of Mathematics, University of Michigan. E-mail addresses: \href{mailto:erhan@umich.edu}{\nolinkurl{erhan@umich.edu}}, \href{mailto:iekren@umich.edu}{\nolinkurl{iekren@umich.edu}}, and \href{mailto:ziqingzh@umich.edu}{\nolinkurl{ziqingzh@umich.edu}}.}\textsuperscript{,\,}\thanks{E. Bayraktar is supported in part by NSF Grant DMS-2602036 and in part by the Susan M. Smith Professorship.}
\and Ibrahim Ekren\footnotemark[2]\textsuperscript{,\,}\thanks{I. Ekren is supported in part by NSF Grant DMS-2406240.}
\and Lu Vy\thanks{Department of Operations Research and Financial Engineering, Princeton University. E-mail address: \href{mailto:lv4809@princeton.edu}{\nolinkurl{lv4809@princeton.edu}}.}
\and Ziqing Zhang\footnotemark[2]}
\date{November 2026}
\hypersetup{
  pdftitle={Primal-Dual Kullback-Leibler Mirror-Prox for Dynamic Mean Field Games},
  pdfauthor={Erhan Bayraktar, Ibrahim Ekren, Lu Vy, Ziqing Zhang},
  pdfsubject={Primal-dual methods for dynamic mean field games},
  pdfkeywords={mean field games, Mirror-Prox, relative entropy, occupation measures}
}

\begin{document}

\maketitle

\begin{abstract}
We compute finite-horizon dynamic mean field equilibria using a primal--dual Kullback--Leibler Mirror-Prox algorithm. Each iteration updates occupation measures by exponential tilting and computes the dual correction by solving a linear PDE, without solving a nonlinear Hamilton--Jacobi--Bellman equation or a best-response control problem. Under assumptions and Lasry--Lions monotonicity, we prove an $O(K^{-1})$ ergodic Minty bound after $K$ iterations and an $O(K^{-1/2})$ approximate-equilibrium bound for the corrected average. We obtain convergence rates for the population flow under strong monotonicity, and for smooth diffusions with uniformly convex control costs under ordinary monotonicity. We also prove convergence of continuous-time Markov chain approximations to diffusion equilibria and, under additional smoothness and convexity assumptions, bound the error in terms of the spatial mesh size and the number of iterations. Galerkin approximations preserve the $O(K^{-1})$ estimates under suitable bounds on the correction errors. Numerical experiments illustrate the method and its convergence.
\end{abstract}

{\setcounter{tocdepth}{1}\scriptsize\tableofcontents}

\section{Introduction}\label{sec:intro}

Mean field games (MFGs), introduced by Lasry and Lions \cite{lasry2007mean} and by Huang, Malham\'e and Caines \cite{huang2006large}, describe the limit of symmetric games with many players. A representative player takes the population flow $\bm\mu=(\mu_t)_{t\in[0,T]}$ as given and chooses a control to minimize an expected cost that depends on this flow. At a mean field equilibrium (MFE), the player has an optimal control under which the state has law $\mu_t$ at each time $t$. The classical PDE formulation couples a backward Hamilton--Jacobi--Bellman (HJB) equation for the player's value function with a forward Fokker--Planck equation for the population flow.

For a fixed population flow, the player's control problem can be written as a linear program over occupation measures \cite{guo2025continuous}. Can Korpelevich's extragradient method for linear programs \cite{korpelevich1976extragradient} be adapted to compute a mean field equilibrium? For $A\in\R^{m\times n}$, $b \in \R^{m}$, and $c \in \R^n$, consider the linear program
\[
    \min_{x\in\R_+^n} c^\top x
    \qquad\text{subject to}\qquad Ax=b.
\]
The method seeks a primal--dual pair $(x^*,u^*)\in\R_+^n\times\R^m$ satisfying the variational inequality
\[
    \langle c-A^\top u^*,\,x-x^*\rangle
    +\langle Ax^*-b,\,u-u^*\rangle \ge 0
    \qquad\forall (x,u)\in\R_+^n\times\R^m.
\]
The solutions of this inequality are the saddle points of the Lagrangian $c^\top x-u^\top(Ax-b)$ on $\R_+^n\times\R^m$. Simultaneous gradient descent--ascent on this Lagrangian can fail to converge \cite{korpelevich1976extragradient}. Korpelevich's method adds a predictor step and uses the gradient at the predicted point for the actual update.

Returning to the game, fix a population flow $\bm\mu$. The variable in the player's linear program is a pair $\pi=(\nu,(\pi_t)_{t\in[0,T]})$,
where $\nu$ is the terminal state law and $\pi_t$ is the joint state--action law at time $t$. Writing $\Pi$ for the set of such pairs of probability laws, the problem is
\[
    \min_{\pi\in\Pi}\langle\cC_{\bm\mu},\pi\rangle
    \qquad\text{subject to}\qquad \cA\pi=\rho.
\]
The pair $\cC_{\bm\mu}$ consists of the terminal and running cost functions. We assume that the state dynamics are independent of the population, so the same forward operator $\cA$ applies to every flow $\bm\mu$. The constraint $\cA\pi=\rho$ imposes the weak forward equation from the initial law $\rho$ (see \Cref{sec:setting}). This equation couples the terminal law $\nu$ and the running state--action laws $(\pi_t)_{t\in[0,T]}$ across time. We dualize this constraint so that each primal step can update these laws separately, as in Korpelevich's treatment of block linear programs \cite[Sec.~2, Observation~2]{korpelevich1976extragradient}.

An equilibrium occupation measure $\pi^*$ solves this linear program for its own flow of state marginals, $\bm\mu=\cM\pi^*$. We therefore seek an occupation measure $\pi^*\in\Pi$ and a multiplier $u^*\in H^1(0,T;H)$ satisfying
\[
    \langle\cC_{\cM\pi^*}-\cA^*u^*,\,\pi-\pi^*\rangle
    +\langle\cA\pi^*-\rho,\,u-u^*\rangle \ge 0
    \qquad\forall (\pi,u)\in\Pi\times H^1(0,T;H),
\]
where $H$ is a Hilbert space of spatial test functions; the multiplier plays the role of the player's value function. If $(\pi^*,u^*)$ solves this inequality, then $\pi^*$ is an equilibrium. Conversely, an equilibrium $\pi^*$ admits a multiplier $u^*$ if the dual problem for $\cM\pi^*$ attains the primal value in $H^1(0,T;H)$ (see \Cref{prop:PDVI-MFE}). Let $\cF(\pi,u)=(\cC_{\cM\pi}-\cA^*u,\cA\pi-\rho)$ denote the operator in this inequality. Since the same forward operator $\cA$ applies to every population flow, the terms involving $\cA$ and $\cA^*$ cancel when we compare two pairs $(\pi,u)$ and $(\pi',u')$, giving
\[
    \langle \cF(\pi,u)-\cF(\pi',u'),(\pi-\pi',u-u')\rangle
    =
    \langle \cC_{\cM\pi}-\cC_{\cM\pi'},\pi-\pi'\rangle.
\]
Lasry--Lions monotonicity of the costs is precisely the requirement that the right-hand side be nonnegative for all $\pi,\pi'\in\Pi$. It is therefore equivalent to monotonicity of the primal--dual operator (\Cref{prop:monotone}), the property that extragradient methods rely on \cite[Thm.~2]{korpelevich1976extragradient}.

Unlike in the linear program, $\cC_{\cM\pi}-\cA^*u$ is a pair of functions while $\pi$ is a pair of measures, so the Euclidean step $\pi-\lambda(\cC_{\cM\pi}-\cA^*u)$ has no meaning. We therefore turn to Nemirovski's Mirror-Prox method \cite{nemirovski2004prox}, which replaces the Euclidean step by a Bregman proximal step. As in our work on static mean field games \cite{KLMP}, we use relative entropy to update the measures by exponential tilting. The tilt moves mass toward state--action pairs with lower HJB residuals $\cC_{\cM\pi}-\cA^*u$ (see \eqref{eq:slacks}), which we can evaluate without solving the player's control problem. For the multiplier we use a Sobolev distance, so the dual update is $u^+=u-\tau w$, where $w$ is the Riesz representative of the residual $\cA\pi-\rho$. With our choice of inner product, $w$ is explicit: we integrate the residual in time and apply the Riesz map of $H$ once at each time, which on a periodic grid is one forward and one inverse FFT (\Cref{lem:riesz} and \eqref{eq:riesz-h}). \Cref{sec:algo} gives the full predictor--corrector iteration (\Cref{alg:PDKLMP}).

\subsection{Main results}

Because the forward equation is enforced only through the multiplier, the averaged measure \(\bar\pi_K\) obtained after \(K\) iterations may violate it. Summing the dual updates gives the identity
\[
    \tns{\cA\bar\pi_K-\rho}
    =\frac{\tn{u_K-u_0}}{\tau K}
\]
from \eqref{eq:dual-drift}: the change in the multiplier measures exactly how far $\bar\pi_K$ is from satisfying the forward equation. For Lipschitz, Lasry--Lions monotone costs and small enough step sizes, we obtain a Minty estimate against every feasible $\pi'$ of finite relative entropy with respect to the initialization $\pi_0$: when the population flow is fixed at $\cM\pi'$, the cost of $\bar\pi_K$ is at most that of $\pi'$ plus $O(K^{-1})$. Because the costs are bounded, the same estimate also gives an $O(K^{-1/2})$ bound on the forward-equation residual (see \Cref{cor:primal-minty-residual}). If a finite-entropy equilibrium has a multiplier, the multipliers $u_K$ stay bounded and the identity above improves this to $O(K^{-1})$.

These estimates have two consequences. First, every limit point of $\bar\pi_K$ is an equilibrium under a finite-entropy approximation assumption (\Cref{thm:cluster-points}). Second, we can obtain an approximate equilibrium at a finite iteration by correcting $\bar\pi_K$ to satisfy the forward equation (\Cref{subsec:projection}). Under suitable bounds on the correction and relative entropy, the corrected measure has Minty gap $O(K^{-1})$ and Stampacchia gap $O(K^{-1/2})$. A feasible measure is an $\varepsilon$-MFE if no deviation lowers the player's cost by more than $\varepsilon$ while the population flow is held fixed, that is, if its Stampacchia gap is at most $\varepsilon$. Thus $O(\varepsilon^{-2})$ iterations suffice to obtain an $\varepsilon$-MFE (see \Cref{thm:quantitative-eps-MFE}).

For continuous-time finite-state games \cite{GomesMohrSouza2013,BayraktarCohen2018FiniteState} with finite action sets, the correction and entropy bounds hold when $\pi_0$ is uniformly positive (see \Cref{thm:finite-state}). The correction keeps the feedback of $\bar\pi_K$, its action law given the state, and recomputes the state laws from $\rho$ under that feedback. The finite-state HJB system also supplies a multiplier for every equilibrium.

How close is the averaged population flow $\cM\bar\pi_K$ to an equilibrium flow? Under strong Lasry--Lions monotonicity, all equilibria share one flow, and if a finite-entropy equilibrium has a multiplier, $\cM\bar\pi_K$ converges to that flow at rate $O(K^{-1/2})$ in a Wasserstein distance (see \Cref{thm:strong-LL}). For a class of smooth one-dimensional diffusions with uniformly convex control costs, we also obtain a Wasserstein convergence rate under ordinary Lasry--Lions monotonicity (see \Cref{thm:convex-diffusion}).

\paragraph{Numerical implementation.} We study two ways to implement the method for one-dimensional periodic diffusions. The first replaces the diffusion by a continuous-time Markov chain on a spatial grid, as in Kushner and Dupuis \cite{KushnerDupuis2001} (\Cref{subsec:mesh}). Bayraktar, Budhiraja and Cohen \cite{BayraktarBudhirajaCohen2018Numerical} also use a Markov chain for the reflected diffusions arising from the queueing models of \cite{BayraktarBudhirajaCohen2019RateControl}, but solve the discrete control problem at every population update; our iterations never solve a control problem.

With a discrete Sobolev norm of order $m\ge3$ on the multiplier, the step sizes can be chosen independently of the grid (\Cref{thm:mesh}); with the Euclidean norm they shrink as the grid is refined. If the number of iterations grows faster than the logarithm of the number of grid points, every limit point of the averaged occupation measures is an equilibrium of the diffusion game (\Cref{prop:mesh-limit}). For smooth models with uniformly convex control costs, we also bound the Wasserstein error in terms of the grid spacing and the number of iterations (\Cref{thm:convex-grid,cor:convex-grid-balanced}).

The second implementation (\Cref{subsec:inexact}) uses the diffusion generator directly and approximates the Riesz correction in a function basis. If the approximation errors are summable over the iterations and a finite-entropy equilibrium has a multiplier, the $O(K^{-1})$ Minty estimate and residual bound survive (\Cref{thm:inexact}). In \Cref{sec:experiments} we study a finite-state congestion game and illustrate both implementations for a periodic diffusion; the Fourier--Galerkin runs keep a fixed set of modes, so they enforce the forward equation only on those modes.

\subsection{Related Work}

\emph{Occupation measures and best responses.}
Occupation measures turn a stochastic control problem into a linear program; see Bhatt and Borkar \cite{BhattBorkar1996} and Kurtz and Stockbridge \cite{KurtzStockbridge1998}. Guo, Hu, Zhang and Zhang \cite{guo2025continuous} use the representative player's linear program and its dual over HJB subsolutions to characterize continuous-time mean field equilibria. Bayraktar, Hu and Zhang \cite{BayraktarHuZhang2026StochasticDuality} verify equilibria with measurable dual functions and, under further assumptions, prove dual attainment and characterize all equilibria. Our estimates that use an optimal multiplier require one in $H^1(0,T;H)$. A closely related computational approach is the linear programming fictitious play of Dumitrescu, Leutscher and Tankov \cite{DumitrescuLeutscherTankov2023} for games with optimal stopping or absorption. As in other fictitious-play methods \cite{cardaliaguet2017learning,hadikhanloo2019finite,perrin2020fictitious,yu2024convergence}, each iteration computes a best response to a given population flow. For Dumitrescu, Leutscher and Tankov this means solving an occupation-measure linear program. We use the same linear programming principle, but take an extragradient step in the occupation measure and its multiplier together.

\emph{PDE methods and potential games.}
Numerical approaches to the HJB--Fokker--Planck system include finite differences \cite{achdou2010mean,achdou2013mean,BonnansLiuPfeiffer2023Theta} and policy iteration \cite{cacace2021policy}. Probabilistic methods work with the associated forward--backward stochastic differential equations \cite{angiuli2019cemracs,CarmonaDelarue2018,ChassagneuxCrisanDelarue2019Numerical}. For surveys, see \cite{achdou2020mean,lauriere2021numerical}. These methods take the value function and the density as unknowns and keep the Hamiltonian, hence the nonlinearity of the HJB equation, inside the iteration. When the game has a potential, its variational formulation also admits primal--dual and proximal splitting methods \cite{briceno2018proximal,liu2021computational}. In particular, Bonnans, Lavigne and Pfeiffer \cite{BonnansLavignePfeiffer2023Discrete} use a relative-entropy update in a primal--dual method for discrete-time potential games. In continuous time, Lavigne and Pfeiffer \cite{lavigne2023generalized} compute each conditional-gradient direction by solving an HJB equation for a best response and then propagating the state law. Liu and Pfeiffer \cite{LiuPfeiffer2025Mesh} obtain mesh-independent convergence bounds for the discretized method. In our formulation the dual variable comes from the representative player's forward equation. This allows us to treat games without a potential. The Hamiltonian never appears in our iteration, and the generator enters only through residuals, so no parabolic equation is solved (\Cref{sec:riesz-bvp}).

\emph{Extragradient methods for MFGs.}
Meynard \cite{meynard2026extragradient} applies extragradient to open-loop controls under displacement monotonicity. For separable costs, he also identifies a Lasry--Lions monotone variational inequality on the flows of state--action laws induced by admissible controls \cite[Sec.~2.1.3]{meynard2026extragradient}. Our choice to dualize the forward equation allows the terminal and running laws to be updated separately. It also means that their averages can violate that equation, so an approximate solution of the variational inequality is not yet an approximate equilibrium. We compare only against feasible occupation measures, using restricted gap estimates in the spirit of \cite{Nesterov2007}, correct the average, and then bound a player's gain from deviating as in \cite{bigi2023approximate}.

\emph{Learning dynamics and policy evaluation.}
Learning dynamics for monotone games include online mirror descent and dual averaging \cite{hadikhanloo2017learning,hadikhanloo2022learning}. For static games, related approaches evolve mixed strategies \cite{hynd2023evolution}, approximate equilibria by Langevin particle systems \cite{tangpi2025particle}, or learn from payoff observations with finite-sample guarantees \cite{yardim2025variational}. In the online mirror descent method of P\'erolat et al.~\cite{perolat2022scaling}, each update evaluates the current policy through its $Q$-function. Isobe, Abe and Ariu \cite{isobe2026last} prove last-iterate convergence for a KL proximal-point method in which each exact step solves a regularized MFG. Our measure update is an explicit tilt. Its exponent is the HJB residual of the current multiplier, which we can evaluate for each action without computing the value of the current policy.

\emph{Reinforcement learning and entropy regularization.}
Reinforcement learning methods for MFGs are developed in \cite{guo2019learning,lauriere2022scalable,magnino2026solving} and surveyed in \cite{lauriere2022learning}. One use of entropy in this literature is to regularize the player's objective \cite{xie2021learning,guo2022entropy,cui2021approximately}. In entropy-regularized mean-field control, Bayraktar, Hernandez, Yan and Zhu \cite{bayraktar2026policygradient} compute policy gradients using instantaneous advantages derived from the current policy's value function. We use relative entropy only to move from one measure to the next, so the equilibria we compute are those of the unregularized game.

\emph{Riesz corrections and diffusion grids.}
Gevorgyan and Gomes \cite{GevorgyanGomes2026} obtain mesh-independent step sizes for a projected extragradient method with a Riesz correction in a first-order price-formation game. Al Abdulaziz, Ashrafyan, Gevorgyan and Gomes \cite{abdulaziz2026bregman} use Bregman projections for regularized stationary MFG systems. In both, the unknowns are a density and a value function, monotonicity comes from the convexity of the Hamiltonian together with a monotone coupling, and the density is kept nonnegative by projection. Our operator is affine apart from the costs, so its monotonicity is exactly Lasry--Lions monotonicity, and the KL tilts keep every primal iterate a probability measure. For diffusion problems, the order of the spatial norm matters. Meng, Liu, Li and Osher \cite{MengLiuLiOsher2024}, who apply the primal--dual hybrid gradient method of Chambolle and Pock \cite{ChambollePock2011} to viscous Hamilton--Jacobi equations, explain why the second spatial derivative makes the step-size bound in their $H^1$ metric depend on the mesh. On our one-dimensional diffusion grids, a spatial Sobolev norm of order $m\ge3$ controls the second-order generator uniformly, so the same step sizes can be used as the grid is refined (\Cref{rmk:gomes-osher}).

\paragraph{Organization.}
The rest of the paper is organized as follows.
\Cref{sec:setting} introduces the occupation measures, the forward equation and the primal--dual variational inequality.
\Cref{sec:algo} presents the KL Mirror-Prox algorithm and derives the dual update.
\Cref{sec:LL} establishes convergence under Lasry--Lions monotonicity.
\Cref{sec:finite} treats finite-state games and Markov chain approximations of diffusions.
\Cref{subsec:inexact} develops the Galerkin implementation, and \Cref{sec:experiments} reports the numerical experiments.
The proofs are collected in Appendices~\ref{sec:proofs-setting} and~\ref{sec:proofs-LL}.
Figure~\ref{fig:map} in Appendix~\ref{app:logical-structure} shows which assumptions enter each result.

\section{Setting and primal--dual formulation}\label{sec:setting}

\paragraph{The game.}
A continuum of symmetric players move in a state space $E$ and choose
actions in $A$ over the horizon $[0,T]$. A player's state $X$ starts
from the law $\cL(X_0)=\rho$ and is a weak controlled process whose generator
$\LL^{a}_t$ depends on the action $a$ in force at time $t$: jump rates on
a finite set, a controlled diffusion on the torus, or a bounded jump
kernel, see \Cref{ex:finite,ex:diffusion,ex:jump} below. A player uses a
relaxed Markov control $\alpha=(\alpha_t(da\mid y))$, a probability on $A$
for each $(t,y)$, under which the generator is
$\LL^{\alpha}_t\phi(y)=\int_A\LL^a_t\phi(y)\,\alpha_t(da\mid y)$. Facing a
population flow $\bm\mu=(\mu_t)_{t\in[0,T]}$ in $\cP(E)$, the player pays
\begin{equation}\label{eq:game-cost}
    J(\alpha;\bm\mu):=\mathbb E\Bigl[\int_0^T\!\!\int_A f(t,X_t,a,\mu_t)\,\alpha_t(da\mid X_t)\,dt+g(X_T,\mu_T)\Bigr],
\end{equation}
where $f$ is a running cost and $g$ a terminal cost, both depending on the
population through the current marginal. A mean field equilibrium is a pair $(\alpha^*,\bm\mu^*)$ such
that $\alpha^*$ minimizes $J(\cdot\,;\bm\mu^*)$ and
$\mu^*_t=\Law(X^{\alpha^*}_t)$ for every $t$. Following
\cite{guo2025continuous}, we encode a control by its occupation measure
\begin{equation}\label{eq:occupation}
    \pi=(\nu,(\pi_t)_{t\in[0,T]}),\qquad \nu:=\Law(X^\alpha_T),\qquad \pi_t(dy,da):=\Law(X^\alpha_t)(dy)\,\alpha_t(da\mid y),
\end{equation}
for which the cost \eqref{eq:game-cost} is linear,
$J(\alpha;\bm\mu)=\int g(\cdot,\mu_T)\,d\nu+\int_0^T\!\!\int f(t,\cdot,\cdot,\mu_t)\,d\pi_t\,dt$,
and the dynamics become the linear constraint that the marginals
$\mu_t=\pi_t(\cdot\times A)$ solve the Kolmogorov forward equation. \Cref{eq:feasible-set} below makes it precise. 

\subsection{Standing assumptions}

Throughout, $T>0$ is fixed. For a compact metric space $S$, $\cP(S)$
denotes the Borel probability measures on $S$ and $\mathcal M_\pm(S)$ the
finite signed Borel measures on $S$. We write
$\lVert \xi\rVert_1:=|\xi|(S)$ for $\xi\in\mathcal M_\pm(S)$. Thus
$\lVert \mu-\mu'\rVert_1$ is twice the total variation distance, as in
\cite{KLMP}. The $1$-Wasserstein distance is $W_1$. Recall that
$W_1(\mu,\mu')\le \tfrac{D_S}{2}\lVert\mu-\mu'\rVert_1$, where $D_S$ is the
diameter of $S$.

\begin{assumption}
\label{ass:standing}
\leavevmode
\begin{enumerate}
\item[(S1)] The state space $E$ and the action space $A$ are compact metric
spaces. Let $D_E:=\operatorname{diam}(E)$. The initial law is
$\rho\in\cP(E)$.
\item[(S2)] $H$ is a separable Hilbert space of real functions on $E$ with
$H\subset C(E)$ and
\[
    \lVert \phi\rVert_\infty\le c_H\nH{\phi},\qquad \phi\in H.
\]
\item[(S3)] (Mean-field-independent dynamics.) For each $(t,a)\in[0,T]\times A$ there is a linear operator $\LL_t^a:H\to C(E)$ with
\[
    \lVert \LL^a_t \phi\rVert_\infty\le \ell_H\nH{\phi},
    \qquad \phi\in H,\ (t,a)\in[0,T]\times A.
\]
For each $\phi\in H$, the map $(t,y,a)\mapsto \LL_t^a\phi(y)$ is Borel, and
$(y,a)\mapsto\LL^a_t\phi(y)$ is continuous for each $t$.
\item[(S4)] (Costs.) The maps $f:[0,T]\times E\times A\times\cP(E)\to\R$ and $g:E\times\cP(E)\to\R$ are Borel and bounded. For each $(t,\mu)$ the map $(y,a)\mapsto f(t,y,a,\mu)$ is continuous, and for each $\mu$ the map $y\mapsto g(y,\mu)$ is continuous.
\item[(S5)] The feasible set $\Pirho$ defined in \eqref{eq:feasible-set} below
is nonempty.
\end{enumerate}
\end{assumption}

The operator $\LL^a_t$ plays the role of the generator
$\LL^{\bm\mu,a}$ in \cite{guo2025continuous}. By (S3) it does not depend
on the population. This is the structural assumption under which the
primal--dual operator is monotone, see \Cref{prop:monotone}. The space $H$
fixes the spatial regularity of the dual variable. Only the two constants
$c_H$ and $\ell_H$ enter the analysis.

\begin{example}[Finite state space, continuous time]
\label{ex:finite}
Let $E=\{1,\ldots,N\}$ carry the discrete metric, so $D_E=1$ and $W_1$ is
the total variation distance. Let $A$ be finite, or compact. Take jump
rates $Q_t(y,y',a)\ge0$ for $y'\ne y$, Borel in $t$ and continuous in $a$,
with
\[
    \bar q:=\sup_{t,y,a}\sum_{y'\neq y}Q_t(y,y',a)<\infty .
\]
Let $\LL^a_t\phi(y)=\sum_{y'\ne y}Q_t(y,y',a)(\phi(y')-\phi(y))$ and let
$H=\R^E$ with the Euclidean inner product. Then $c_H=1$ and
$\ell_H=2\bar q$. This is the class of continuous-time finite-state MFGs of
\cite{GomesMohrSouza2013}, with mean-field-independent rates.
For finite-state games, Bayraktar and Cohen \cite{BayraktarCohen2018FiniteState}
use the master equation to approximate the players' value functions in the finite-player game.
\end{example}

\begin{example}[Diffusions on the torus]
\label{ex:diffusion}
Let $E=\mathbb T^d$ and
$\LL^a_t\phi=b(t,\cdot,a)\cdot\nabla \phi+\frac12\operatorname{tr}(\sigma\sigma^\top(t,\cdot,a)\nabla^2\phi)$,
with $b,\sigma$ bounded, Borel in $t$ and continuous in $(y,a)$. Let
$H=H^m(\mathbb T^d)$ with $m>d/2+2$. Sobolev embedding gives
$\lVert \phi\rVert_{C^2}\le C_m\nH{\phi}$, where
$\lVert \phi\rVert_{C^2}:=\max\{\lVert \phi\rVert_\infty,\lVert|\nabla \phi|\rVert_\infty,\lVert|\nabla^2\phi|_F\rVert_\infty\}$
(Euclidean norm of the gradient, Frobenius norm of the Hessian). Hence (S2)--(S3) hold with
$c_H=C_m$ and
$\ell_H=C_m(\lVert |b|\rVert_\infty+\tfrac12\lVert|\sigma\sigma^\top|_F\rVert_\infty)$.
This is the setting of \cite{guo2025continuous}, under the extra
assumption that $b$ and $\sigma$ do not depend on the population.
\end{example}

\begin{example}[Bounded jump kernels]
\label{ex:jump}
Let $E$ be compact metric and
$Q_t(y,a,E)\le\bar q$, and define
\[
    \LL^a_t\phi(y)=\int_E\bigl(\phi(y')-\phi(y)\bigr)\,Q_t(y,a,dy').
\]
Assume that
$t\mapsto Q_t(y,a,B)$ is Borel for every Borel $B$, and
$(y,a)\mapsto Q_t(y,a,\cdot)$ is weakly continuous. Let $H$ be any
separable reproducing kernel Hilbert space on $E$ with a continuous kernel
$\mathsf k$. Then $c_H=\sup_y \mathsf k(y,y)^{1/2}$ and $\ell_H=2\bar q c_H$.
These bounds do not imply that $H$ determines the full forward equation. For bounded jump generators, uniform density of $H$ in $C(E)$ is a
sufficient test-space condition: uniform approximation of $\phi$ also
approximates $\LL_t^a\phi$, since
$\|\LL_t^a(\phi-\phi_n)\|_\infty\le2\bar q\|\phi-\phi_n\|_\infty$.
\end{example}

\begin{example}[Finite-dimensional test spaces]
\label{ex:finite-dim}
Let $\phi_1,\ldots,\phi_M\in C(E)$ be such that $\LL^a_t\phi_j$ is defined, Borel in
$(t,y,a)$ and continuous in $(y,a)$ with $\sup_{t,a}\lVert\LL^a_t\phi_j\rVert_\infty<\infty$,
and let $H:=\operatorname{span}\{\phi_1,\ldots,\phi_M\}$ with any inner product,
say the one making $(\phi_j)$ orthonormal. Then (S2)--(S3) hold with
$c_H=\sup_{|c|=1}\lVert\sum_jc_j\phi_j\rVert_\infty$ and
$\ell_H=\sup_{t,a}\sup_{|c|=1}\lVert\sum_jc_j\LL^a_t\phi_j\rVert_\infty$. This
is the choice made in computations, see \Cref{sec:algo}. For this $H$ the
algorithm and the results below apply to the game in which the forward
equation is imposed against the test functions
$t\mapsto\sum_jc_j(t)\phi_j$ with $c\in H^1(0,T;\R^M)$.
\end{example}

\subsection{Primal space, cost space and dual space}

\paragraph{Primal variables.}
Let $\cX$ be the space of pairs $\pi=(\nu,(\pi_t)_{t\in[0,T]})$ with
$\nu\in\mathcal M_\pm(E)$ and $(\pi_t)_{t\in[0,T]}$ a measurable family in
$\mathcal M_\pm(E\times A)$ (that is, $t\mapsto \pi_t(B)$ is Borel for every
Borel $B$), identified up to $dt$-null sets and such that $\int_0^T\lVert \pi_t\rVert_1^2\,dt<\infty$. We use the same letter for the pair and for its running family, and $\nu$ for the terminal component. We combine total variation in an $L^2$-in-time norm, including the terminal measure:
\[
    \nsharp{\pi}^2:=\lVert \nu\rVert_1^2+\int_0^T\lVert \pi_t\rVert_1^2\,dt .
\]
The set of probability pairs is
\[
    \Pi:=\bigl\{\pi=(\nu,(\pi_t)_{t\in[0,T]})\in\cX:\ \nu\in\cP(E),\ \pi_t\in\cP(E\times A)\text{ for a.e. }t\bigr\}.
\]
It is convex and $\operatorname{diam}_{\mathrm{TV},2}(\Pi)\le 2\sqrt{1+T}$.
For $\pi\in\Pi$, the population flow $\cM \pi=(\mu_t)_{t\in[0,T]}$ is
\begin{equation}\label{eq:mean-field-map}
    \mu_t:=\pi_t(\cdot\times A)\quad(t\in[0,T)),\qquad \mu_T:=\nu .
\end{equation}
The map $\cM$ is affine, and $\cM$ is defined for every $\pi\in\Pi$.

\paragraph{Cost space.}
Let $\cG$ be the space of pairs $\gamma=(G,\Gamma)$ with $G\in C(E)$ and
$\Gamma:[0,T]\times E\times A\to\R$ Borel. We require that
$\Gamma_t:=\Gamma(t,\cdot,\cdot)\in C(E\times A)$ for a.e.\ $t$ and that
\[
    \nG{\gamma}^2:=\lVert G\rVert_\infty^2+\int_0^T\lVert\Gamma_t\rVert_\infty^2\,dt<\infty .
\]
The map $t\mapsto\lVert\Gamma_t\rVert_\infty$ is measurable, since the
supremum may be taken over a countable dense subset of $E\times A$. The
pairing between $\cG$ and $\cX$ is
\[
    \pair{\gamma}{\pi}:=\int_E G\,d\nu+\int_0^T\!\!\int_{E\times A}\Gamma_t\,d\pi_t\,dt .
\]
By Cauchy--Schwarz, $|\pair{\gamma}{\pi}|\le\nG{\gamma}\nsharp{\pi}$.
For a measurable flow $\bm\mu=(\mu_t)_{t\in[0,T]}$ in $\cP(E)$, the
\emph{cost vector} is
\begin{equation}\label{eq:cost-vector}
    \cC_{\bm\mu}:=\bigl(g(\cdot,\mu_T),\ (t,y,a)\mapsto f(t,y,a,\mu_t)\bigr)\in\cG,
    \qquad
    \nG{\cC_{\bm\mu}}\le C_0:=\bigl(\lVert g\rVert_\infty^2+T\lVert f\rVert_\infty^2\bigr)^{1/2}.
\end{equation}

\paragraph{Dual variables.}
Let $H^1(0,T;H)$ be the space of $u\in L^2(0,T;H)$ with weak
derivative $\partial_tu\in L^2(0,T;H)$. In norm subscripts we abbreviate
$H^1(0,T;H)$ by $H^1$. We use the continuous representative
$u\in C([0,T];H)$ and the inner product
\[
    \pp{u}{v}:=\pair{u(T)}{v(T)}_H+\int_0^T\pair{\partial_tu(t)}{\partial_tv(t)}_H\,dt ,
    \qquad \tn{u}:=\pp{u}{u}^{1/2}.
\]
This inner product induces a Hilbert norm equivalent to the standard
Sobolev norm on $H^1(0,T;H)$. We write $u(t,y):=u(t)(y)$. The dual norm of
$\ell\in H^1(0,T;H)^*$ is $\tns{\ell}:=\sup\{\pair{\ell}{v}:\tn{v}\le1\}$.
The Riesz isomorphisms are $R:H^1(0,T;H)^*\to H^1(0,T;H)$ and
$R_H:H^*\to H$. For each $\ell\in H^*$, $R_H\ell$ is the unique element of
$H$ satisfying $(R_H\ell,\phi)_H=\pair{\phi}{\ell}$ for every $\phi\in H$;
$R$ is defined in the same way using $\pp{\cdot}{\cdot}$. Every
$\xi\in\mathcal M_\pm(E)$ defines an element of $H^*$ by
$\phi\mapsto\int \phi\,d\xi$, with $\lVert \xi\rVert_{H^*}\le c_H\lVert \xi\rVert_1$.

\paragraph{Forward operator.}
For $v\in H^1(0,T;H)$ define
\begin{equation}\label{eq:A-adjoint}
    \cA^*v:=\Bigl(v(T,\cdot),\ -\bigl(\partial_tv(t,y)+\LL^a_tv(t,\cdot)(y)\bigr)\Bigr).
\end{equation}
Define $\cA:\cX\to H^1(0,T;H)^*$ by
\[
    \pair{\cA \pi}{v}:=\pair{\cA^*v}{\pi}
    =\int_Ev(T,y)\,\nu(dy)-\int_0^T\!\!\int_{E\times A}\bigl(\partial_tv+\LL^a_tv\bigr)(t,y)\,\pi_t(dy,da)\,dt .
\]
We identify the initial law $\rho$ with the element
$v\mapsto\int_Ev(0,y)\,\rho(dy)$ of $H^1(0,T;H)^*$ and write
\[
    \pair{\rho}{v}:=\int_Ev(0,y)\,\rho(dy),\qquad v\in H^1(0,T;H),
\]
so that $\cA\pi-\rho\in H^1(0,T;H)^*$ is the residual of the forward equation
started from $\rho$.

\begin{lemma}
\label{lem:elementary}
Under \Cref{ass:standing}, every $u\in\HY$ satisfies $\cA^*u\in\cG$ and
\[
    \sup_{t\in[0,T]}\nH{u(t)}\le\sqrt{1+T}\,\tn{u},\qquad
    \nG{\cA^*u}\le L_A\tn{u},\qquad
    L_A:=\bigl(2c_H^2+2\ell_H^2T(1+T)\bigr)^{1/2}.
\]
Consequently $\tns{\cA \pi}\le L_A\nsharp{\pi}$ for $\pi\in\cX$, and
$|\pair{\rho}{u}|\le c_H\sqrt{1+T}\,\tn{u}$.
\end{lemma}

The feasible set consists of the tested solutions of the forward equation
started from $\rho$,
\begin{equation}\label{eq:feasible-set}
    \Pi(\rho):=\Bigl\{\pi\in\Pi:\ \pair{\cA \pi}{v}=\int_Ev(0,y)\,\rho(dy)\ \ \forall v\in H^1(0,T;H)\Bigr\}
    =\{\pi\in\Pi:\ \cA\pi=\rho\ \text{in }H^1(0,T;H)^*\}.
\end{equation}
Written out, $\pi=(\nu,(\pi_t)_t)\in\Pi(\rho)$ if and only if, for every
$v\in H^1(0,T;H)$,
\[
    \int_Ev(T,y)\,\nu(dy)-\int_Ev(0,y)\,\rho(dy)
    =\int_0^T\!\!\int_{E\times A}\bigl(\partial_tv+\LL^a_tv\bigr)(t,y)\,\pi_t(dy,da)\,dt .
\]
This is the forward equation tested in $H^1(0,T;H)$, corresponding to the (Primal Feasibility) condition of
\cite[Theorem~4.5]{guo2025continuous}.
The set $\Pi(\rho)$ is convex. In the diffusion model of
\Cref{ex:diffusion}, $H^m(\mathbb T^d)$ contains the smooth periodic
spatial test functions, so this constraint determines the full weak
Fokker--Planck equation. Disintegrating
$\pi_t(dy,da)=\mu_t(dy)\alpha_t(da\mid y)$ and applying the superposition
argument of \cite[Props.~{7.2--7.3}]{guo2025continuous} gives a weak
relaxed-control state--policy pair with state marginals $\mu_t$ and
terminal law $\nu$. Thus each feasible pair is realized as an occupation
measure in this model. The state flow is understood through its narrowly
continuous representative, with $\mu_0=\rho$ and $\mu_T=\nu$.
For general generators and test spaces, the relation to controlled
processes is specified in \Cref{rmk:MFE-process}.

For feasible $\pi$ and a frozen flow $\bm\mu$, the number
$\pair{\cC_{\bm\mu}}{\pi}$ is the cost
\[
    \int_E g(y,\mu_T)\,\nu(dy)+\int_0^T\!\!\int_{E\times A} f(t,y,a,\mu_t)\,\pi_t(dy,da)\,dt
\]
of the representative player, as in \cite{guo2025continuous}.

\subsection{Mean field equilibria}

Here $\cM \pi$, defined in \eqref{eq:mean-field-map}, is the population flow
generated by $\pi$, and $\cC_{\cM \pi}$, defined in \eqref{eq:cost-vector}, is its
cost vector.

\begin{definition}
\label{def:MFE}
A point $\pi^*\in\Pirho$ is a \emph{mean field equilibrium} (MFE) if
\[
    \pair{\cC_{\cM \pi^*}}{\pi^*}\le\pair{\cC_{\cM \pi^*}}{\pi}\qquad\forall \pi\in\Pirho .
\]
For $\pi\in\Pirho$, the \emph{Stampacchia gap}, also called exploitability, is
\[
    \Expl(\pi)
    :=\sup_{\pi'\in\Pirho}
       \pair{\cC_{\cM\pi}}{\pi-\pi'}.
\]
The \emph{Minty gap} is
\begin{equation}\label{eq:Minty-gap}
 \GapM(\pi)
 :=\sup_{\pi'\in\Pirho}
       \pair{\cC_{\cM\pi'}}{\pi-\pi'}.
\end{equation}
Both gaps are nonnegative, since both suprema include $\pi'=\pi$, and $\pi$ is an MFE if and only if $\Expl(\pi)=0$. For $\varepsilon\ge0$, we call $\pi$ an \emph{$\varepsilon$-MFE} if $\Expl(\pi)\le\varepsilon$.
\end{definition}

For an occupation measure $\pi$ of a relaxed control $\alpha$,
\eqref{eq:occupation}--\eqref{eq:cost-vector} give
$\cM \pi=(\Law(X_t^\alpha))_{t\in[0,T]}$ and
$\pair{\cC_{\bm\mu}}{\pi}=J(\alpha;\bm\mu)$. Therefore, whenever $\Pirho$
consists of occupation measures, let $\alpha^*$ be the control represented by
$\pi^*$ and set $\bm\mu^*=\cM \pi^*$. The defining inequality is exactly the
best-response condition $J(\alpha^*;\bm\mu^*)\le J(\alpha;\bm\mu^*)$ for every
admissible $\alpha$, while $\mu_t^*=\Law(X_t^{\alpha^*})$ is the consistency
condition. These are the usual MFG equilibrium conditions.

\begin{rmk}
\label{rmk:MFE-process}
For a controlled process satisfying Dynkin's formula on $H$, its occupation measure \eqref{eq:occupation} belongs to $\Pirho$. The converse holds for the finite-state model of \Cref{ex:finite}, whose test space is
all of $\R^E$. For diffusions and bounded jump kernels, it requires a
determining test class and the applicable representation hypotheses
\cite{BhattBorkar1996,KurtzStockbridge1998,guo2025continuous}. In
\Cref{ex:diffusion}, the Sobolev space contains the smooth spatial tests;
in \Cref{ex:jump}, uniform density in $C(E)$ suffices to recover all
continuous spatial tests. For a general test space $H$, \Cref{def:MFE} describes the problem with
the forward equation imposed against $H$. When the preceding representation
hypotheses hold, these feasible pairs are occupation measures of controlled
processes. This distinction also applies to the finite-dimensional test
spaces in \Cref{ex:finite-dim}; see the discussion after
\Cref{lem:riesz-explicit}. Under this identification, the Stampacchia gap is the gain from the best admissible deviation against the flow $\cM\pi$.
\end{rmk}

\subsection{The primal--dual variational inequality}

Let $\cZ:=\Pi\times H^1(0,T;H)$. For $z=(\pi,u)\in\cZ$ define
\begin{equation}\label{eq:PD-operator}
    \cF(z):=\bigl(\cC_{\cM \pi}-\cA^*u,\ \cA \pi-\rho\bigr)\in\cG\times\HY^*,
\end{equation}
with pairing
$\pair{\cF(z)}{z'-z''}:=\pair{\cC_{\cM \pi}-\cA^*u}{\pi'-\pi''}+\pair{\cA \pi-\rho}{u'-u''}$.
Explicitly, with $\bm\mu=\cM \pi$, the first component of $\cF(z)$ consists
of the terminal and running \emph{HJB residuals}
\begin{equation}\label{eq:slacks}
    G^z(y):=g(y,\mu_T)-u(T,y),\qquad
    \Gamma^z(t,y,a):=f(t,y,a,\mu_t)+\partial_t u(t,y)+\LL^a_t u(t,\cdot)(y).
\end{equation}

For an HJB subsolution these residuals are nonnegative. A zero residual
means that the corresponding HJB inequality holds with equality.
At a general iterate the residuals can have either sign.

\begin{definition}
\label{def:PDVI}
A point $z^*=(\pi^*,u^*)\in\cZ$ solves the \emph{primal--dual VI} if
$\pair{\cF(z^*)}{z-z^*}\ge0$ for all $z\in\cZ$.
\end{definition}

\begin{proposition}
\label{prop:PDVI-MFE}
Suppose \Cref{ass:standing} holds.
\begin{enumerate}
\item[(i)] If $(\pi^*,u^*)$ solves the PD-VI, then $\pi^*$ is an MFE.
\item[(ii)] Let $\pi^*$ be an MFE. Suppose there is an \emph{optimal dual variable}
$u^*\in\HY$, that is, a function such that
$G^{(\pi^*,u^*)}\ge0$ on $E$, $\Gamma^{(\pi^*,u^*)}(t,\cdot,\cdot)\ge0$
on $E\times A$ for a.e.\ $t$, and
$\pair{\cC_{\cM \pi^*}}{\pi^*}=\int_E u^*(0,y)\,\rho(dy)=\pair{\rho}{u^*}$ (zero duality gap). Then
$(\pi^*,u^*)$ solves the PD-VI.
\item[(iii)] A PD-VI multiplier $u^*$ need not itself be an optimal
dual variable in the sense of (ii). If additionally $1\in H$ and
$\LL_t^a1=0$, it can be normalized to one. In fact, writing
$G^*=G^{(\pi^*,u^*)}$ and $\Gamma^*=\Gamma^{(\pi^*,u^*)}$, set
\begin{align*}
    \underline G&:=\min_{y\in E}G^*(y),\\
    \underline\Gamma(t)&:=\min_{(y,a)\in E\times A}\Gamma^*(t,y,a)
      \quad\text{for a.e. }t,\\
    c(t)&:=\underline G+\int_t^T \underline\Gamma(s)\,ds .
\end{align*}
Then $\widehat u^*(t,y):=u^*(t,y)+c(t)$ belongs to $\HY$ and is an
optimal dual variable with nonnegative HJB residuals and the same value
as the primal problem.
\end{enumerate}
\end{proposition}

The distinction in (iii) comes from the probability normalizations in
$\Pi$: spatially constant, time-dependent shifts of $u$ leave the
PD-VI unchanged when $1\in H$ and $\LL_t^a1=0$, but do change the HJB
residuals. The normalization in part (iii) gives an HJB subsolution
that attains the dual value.

\section{Primal--dual KL Mirror-Prox}\label{sec:algo}

\paragraph{Reference measures and relative entropy.}
Choose an initial probability pair
$\pi_0=(\nu_0,(\pi_{0,t})_t)\in\Pi$. This pair is the reference for the
relative entropy used in the updates and determines their support. A convenient product choice is
\begin{equation}\label{eq:initial-reference}
    \nu_0=r_T,\qquad
    \pi_{0,t}(dy,da)=r_t(dy)\,\zeta(da),
\end{equation}
where $(r_t)_t$ is a measurable family of state laws and
$\zeta\in\cP(A)$ is a full-support action reference. For a finite action
set one may take $\zeta(\{a\})=1/|A|$; on a compact action space one may
choose any full-support probability law. The pair $\pi_0$ need not belong to $\Pi(\rho)$: the prescribed initial law $\rho$ is
imposed by the forward constraint.

For $\pi=(\nu,(\pi_t)_t)$ and $\pi'=(\nu',(\pi'_t)_t)$ in $\Pi$, the dynamic KL divergence is
\[
    \bD(\pi'\|\pi):=\KL(\nu'\|\nu)+\int_0^T\KL(\pi'_t\|\pi_t)\,dt\in[0,\infty],
\]
where
\[
    \KL(\eta\|\xi):=
    \begin{cases}
      \displaystyle\int\log\!\left(\frac{d\eta}{d\xi}\right)d\eta,
          &\eta\ll\xi,\\
      +\infty,&\text{otherwise}.
    \end{cases}
\]
For the product reference \eqref{eq:initial-reference} and the
disintegration $\pi_t(dy,da)=\mu_t(dy)\alpha_t(da\mid y)$, the entropy
chain rule gives
\[
    \KL(\pi_t\|r_t\otimes\zeta)
    =\KL(\mu_t\|r_t)
      +\int_E\KL\bigl(\alpha_t(\cdot\mid y)\|\zeta\bigr)\,\mu_t(dy).
\]
Thus the joint occupation-measure reference specifies both a state
reference and an action reference. The integrand defining $\bD$ is
measurable in $t$ because $\KL$ is jointly lower semicontinuous.
By Pinsker's
inequality applied to $\nu$ and to each $\pi_t$,
\begin{equation}\label{eq:pinsker-dyn}
    \bD(\pi'\|\pi)\ge\tfrac12\nsharp{\pi'-\pi}^2 .
\end{equation}
Fix step sizes $\lambda,\tau>0$ and define
\[
    \Phi(z',z):=\frac1\lambda\bD(\pi'\|\pi)+\frac1{2\tau}\tn{u'-u}^2,
    \qquad z=(\pi,u),\ z'=(\pi',u').
\]

\begin{definition}
\label{def:prox}
Let $z=(\pi,u)\in\cZ$ with $\pi=(\nu,(\pi_t)_t)$, and let $(\gamma,\ell)\in\cG\times H^1(0,T;H)^*$ with $\gamma=(G,\Gamma)$. Define $\prox_z(\gamma,\ell):=(\pi^+,u^+)$, where $\pi^+=(\nu^+,(\pi^+_t)_t)$ and
\[
    \frac{d\nu^+}{d\nu}=\frac{e^{-\lambda G}}{\int e^{-\lambda G}d\nu},\qquad
    \frac{d\pi^+_t}{d\pi_t}=\frac{e^{-\lambda\Gamma_t}}{\int e^{-\lambda\Gamma_t}d\pi_t}
    \quad(\text{a.e. }t),\qquad
    u^+=u-\tau R\ell,
\]
where $R$ is the Riesz isomorphism of $\HY$. In the algorithm below, $G$ and $\Gamma_t$ are the HJB residuals defined in \eqref{eq:slacks}, evaluated at the current iterate for the predictor and at the predicted pair for the corrector.
\end{definition}

\Cref{lem:PD-three-point} shows that $\prox_z(\gamma,\ell)$ is the unique
minimizer over $\cZ$ of
$z'\mapsto\pair{\gamma}{\pi'}+\pair{\ell}{u'}+\Phi(z',z)$.

\paragraph{The algorithm.}
To solve the PD-VI of \Cref{def:PDVI}, we apply Mirror-Prox to $\cF$ using the proximal map in \Cref{def:prox}. Given $z_0=(\pi_0,u_0)\in\cZ$, the resulting \emph{primal--dual KL Mirror-Prox} (PD-KL-MP) iteration is
\begin{equation}\label{eq:PD-KL-MP}
    \underbrace{z_{k+1/2}:=\prox_{z_k}\bigl(\cF(z_k)\bigr)}_{\text{predictor}},\qquad
    \underbrace{z_{k+1}:=\prox_{z_k}\bigl(\cF(z_{k+1/2})\bigr)}_{\text{corrector}},\qquad k\ge0,
\end{equation}
with ergodic averages
\begin{equation}\label{eq:ergodic-avg}
    \bar z_K=(\bar \pi_K,\bar u_K):=\frac1K\sum_{k=0}^{K-1}z_{k+1/2}.
\end{equation}
The average $\bar \pi_K$ is taken componentwise, in $\nu$ and in $\pi_t$ for each $t$, so $\bar \pi_K\in\Pi$. Moreover $\cM\bar \pi_K$ is the average of the flows, because $\cM$ is affine.

More explicitly, the predictor is obtained from the two regularized minimization problems
\[
\begin{aligned}
 \pi_{k+1/2}&=\operatorname*{argmin}_{\pi\in\Pi}
 \left\{\pair{\cC_{\cM \pi_k}-\cA^*u_k}{\pi}+\frac1\lambda\bD(\pi\|\pi_k)\right\},\\
 u_{k+1/2}&=\operatorname*{argmin}_{u\in\HY}
 \left\{\pair{\cA \pi_k-\rho}{u}+\frac1{2\tau}\tn{u-u_k}^2\right\}.
\end{aligned}
\]
This is the dynamic analogue of the regularized linear minimization in the
static KL Mirror-Prox method \cite{KLMP}. Solving these problems gives the
exponential tilts and the Riesz correction below.
Both predictor and corrector use $\bD(\,\cdot\,\|\pi_k)$. Each update has a strictly positive density relative to that base;
hence all iterates remain mutually absolutely continuous with $\pi_0$
componentwise. The initial reference therefore fixes their support.

\noindent\textbf{Predictor.} Write $\pi_k=(\nu_k,(\pi_{k,t})_t)$ and
$\bm\mu_k=\cM \pi_k$. By \eqref{eq:slacks}, the predictor step reads
\begin{equation}\label{eq:explicit-primal}
\begin{aligned}
    \nu_{k+1/2}(dy)&\propto \exp\bigl(-\lambda[g(y,\mu_{k,T})-u_k(T,y)]\bigr)\,\nu_k(dy),\\
    \pi_{k+1/2,t}(dy,da)&\propto \exp\bigl(-\lambda[f(t,y,a,\mu_{k,t})+\partial_t u_k(t,y)+\LL^a_t u_k(t,\cdot)(y)]\bigr)\,\pi_{k,t}(dy,da),\\
    u_{k+1/2}&=u_k-\tau R(\cA \pi_k-\rho).
\end{aligned}
\end{equation}
\noindent\textbf{Corrector.} Set $\bm\mu_{k+1/2}=\cM \pi_{k+1/2}$ and
evaluate the HJB residuals and the forward-equation residual at $z_{k+1/2}$. The
corrector step reads
\[
\begin{aligned}
    \nu_{k+1}(dy)&\propto \exp\bigl(-\lambda[g(y,\mu_{k+1/2,T})-u_{k+1/2}(T,y)]\bigr)\,\nu_k(dy),\\
    \pi_{k+1,t}(dy,da)&\propto \exp\Bigl(-\lambda\bigl[f(t,y,a,\mu_{k+1/2,t})+\partial_t u_{k+1/2}(t,y)\\
    &\hspace{6em}+\LL^a_t u_{k+1/2}(t,\cdot)(y)\bigr]\Bigr)\,\pi_{k,t}(dy,da),\\
    u_{k+1}&=u_k-\tau R(\cA \pi_{k+1/2}-\rho).
\end{aligned}
\]
The base measures remain $\nu_k$ and $\pi_{k,t}$, and the dual correction is
subtracted from $u_k$, not from $u_{k+1/2}$. Thus both steps start from
$z_k$; only the point where the residuals are evaluated changes. This is
Korpelevich's predictor--corrector scheme
\cite{korpelevich1976extragradient} on $\Pi\times\HY$, with relative
entropy for the measure updates and the squared Hilbert norm for the
dual updates \cite{nemirovski2004prox}. The primal step increases the
weight of state--action pairs with lower HJB residuals. At a normalized
PD-VI solution, the HJB inequalities hold with equality on the occupied
state--action pairs, as in \Cref{prop:PDVI-MFE}(iii).
The dual step is made explicit by the following lemma.

\begin{lemma}
\label{lem:riesz}
Let $\pi=(\nu,(\pi_t)_t)\in\Pi$ and $\bm\mu=\cM \pi$. Define $\LL_t^*\pi_t\in H^*$ by
\[
 \pair{\phi}{\LL_t^*\pi_t}:=\int_{E\times A}\LL^a_t\phi(y)\,\pi_t(dy,da),
\]
so that
$\lVert \LL_t^*\pi_t\rVert_{H^*}\le\ell_H$. Define the \emph{integrated forward residuals}
\[
    d_s:=\mu_s-\rho-\int_0^s \LL_t^*\pi_t\,dt\in H^*\ (\text{a.e. }s\in[0,T]),\qquad
    d_T:=\nu-\rho-\int_0^T \LL_t^*\pi_t\,dt\in H^* .
\]
Then $\pair{\cA \pi-\rho}{v}=\pair{v(T)}{d_T}-\int_0^T\pair{\partial_sv(s)}{d_s}\,ds$
for all $v\in H^1(0,T;H)$. Moreover
\begin{equation}\label{eq:time-exact}
    \bigl(R(\cA \pi-\rho)\bigr)(t)=R_Hd_T+\int_t^TR_Hd_s\,ds,
    \qquad
    \tns{\cA \pi-\rho}^2=\lVert d_T\rVert_{H^*}^2+\int_0^T\lVert d_s\rVert_{H^*}^2\,ds .
\end{equation}
In particular, $\pi\in\Pi(\rho)$ if and only if $d_T=0$ and $d_s=0$ for a.e.\ $s$.
\end{lemma}

The dual update therefore integrates the quantities $d_s$ in time and
applies $R_H$. The squared norm of the forward-equation residual is the
sum of the terminal term $\|d_T\|_{H^*}^2$ and the integral of
$\|d_s\|_{H^*}^2$.

\begin{algorithm}[t]
\caption{PD-KL Mirror-Prox (one iteration)}\label{alg:PDKLMP}
\begin{algorithmic}[1]
\Require $z_k=(\pi_k,u_k)$ with $\pi_k=(\nu_k,(\pi_{k,t})_t)$, step sizes $\lambda,\tau$
\State \textbf{Predictor:} $\bm\mu_k\gets\cM \pi_k$; compute the HJB residuals $G^{z_k},\Gamma^{z_k}$ from \eqref{eq:slacks} and the integrated residuals $d^{(k)}$ of $\pi_k$ (\Cref{lem:riesz})
\State $\nu_{k+1/2}\propto e^{-\lambda G^{z_k}}\nu_k$;
       $\pi_{k+1/2,t}\propto e^{-\lambda\Gamma^{z_k}_t}\pi_{k,t}$
\State $u_{k+1/2}(t)\gets u_k(t)-\tau \bigl[R_Hd^{(k)}_T+\int_t^T R_H  d^{(k)}_s\,ds\bigr]$
\State \textbf{Corrector:} $\bm\mu_{k+1/2}\gets\cM \pi_{k+1/2}$; compute HJB residuals and the integrated residuals at $z_{k+1/2}$
\State $\nu_{k+1}\propto e^{-\lambda G^{z_{k+1/2}}}\nu_k$;
       $\pi_{k+1,t}\propto e^{-\lambda\Gamma^{z_{k+1/2}}_t}\pi_{k,t}$
\State $u_{k+1}(t)\gets u_k(t)-\tau\bigl[R_Hd^{(k+1/2)}_T+\int_t^TR_Hd^{(k+1/2)}_s\,ds\bigr]$
\end{algorithmic}
\end{algorithm}

\subsection{The dual step as a linear boundary value problem}\label{sec:riesz-bvp}

\Cref{lem:riesz} gives the dual correction $w=R(\cA \pi-\rho)$ in closed
form. We now read that formula as the solution of a linear two-point
boundary value problem in time whose spatial part is the Riesz map of $H$.
Two facts follow. The generator $\LL$ enters only through the data, so no
parabolic equation is solved. And the time problem is solved exactly, so
the dual variable needs no discretization in $t$.

\paragraph{Weak and strong form.}
By definition of $R$, $w$ is the unique element of $H^1(0,T;H)$ with
\begin{equation}\label{eq:continuum-riesz-weak}
    (w(T),v(T))_H+\int_0^T(\partial_tw,\partial_tv)_H\,dt
    =\pair{\cA \pi-\rho}{v}
    \qquad\text{for all }v\in H^1(0,T;H) .
\end{equation}
Recall from \Cref{lem:riesz} that $\LL_t^*\pi_t\in H^*$,
$\pair{\phi}{\LL_t^*\pi_t}=\int\LL_t^a\phi\,d\pi_t$, is the adjoint of the relaxed
generator applied to $\pi_t$, and that the integrated forward residual is
$d_s=\mu_s-\rho-\int_0^s\LL_t^*\pi_t\,dt$. Let $R_H^{-1}:H\to H^*$,
$\phi\mapsto(\phi,\cdot)_H$, be the inverse Riesz map. Suppose that
$t\mapsto\mu_t$ is absolutely continuous into $H^*$, with one-sided limits
$\mu_{0^+}$ and $\mu_{T^-}$. Then $\partial_td_t=\partial_t\mu_t-\LL_t^*\pi_t$,
$d_{0^+}=\mu_{0^+}-\rho$ and $d_T-d_{T^-}=\nu-\mu_{T^-}$, and
\Cref{lem:riesz} ($\partial_tw=-R_Hd_t$, $w(T)=R_Hd_T$) says that $w$ solves
\begin{equation}\label{eq:riesz-strong}
\begin{aligned}
    -R_H^{-1}\partial_{tt}w(t)&=\partial_t\mu_t-\LL_t^*\pi_t
    \qquad\text{for a.e. }t\in(0,T),\\
    R_H^{-1}\partial_tw(0)&=\rho-\mu_{0^+},\qquad
    R_H^{-1}\bigl[w(T)+\partial_tw(T)\bigr]=\nu-\mu_{T^-}.
\end{aligned}
\end{equation}
This is a second-order linear equation in time with a Neumann condition at
$t=0$ carrying the initial-law mismatch, a Robin condition at $t=T$
carrying the terminal-law mismatch, and the Fokker--Planck residual as
source. Conversely, \eqref{eq:riesz-strong} is the Euler--Lagrange equation
of \eqref{eq:continuum-riesz-weak}, obtained by integrating the time
integral by parts. For a general $\pi\in\Pi$ the flow $t\mapsto\mu_t$ is
merely in $L^2$, the strong form is formal, and
\eqref{eq:continuum-riesz-weak} together with \Cref{lem:riesz} is its
meaning.

The next lemma lists $R_H$ for the spaces of \Cref{ex:finite,ex:diffusion,ex:jump,ex:finite-dim}.

\begin{lemma}
\label{lem:riesz-explicit}
Let $\ell\in H^*$.
\begin{enumerate}
\item[(a)] If $E$ is finite and $H=\R^E$ with the Euclidean inner product, write $\xi_y=\ell(\mathbf 1_{\{y\}})$. Then $R_H\ell=\xi$ as a vector and $\lVert \ell\rVert_{H^*}=|\xi|_2$.
\item[(b)] If $E=\mathbb T^d$ and $H=H^m(\mathbb T^d)$ with $(\phi,\phi')_{H^m}=\sum_{i=0}^m\binom mi(\nabla^i\phi,\nabla^i\phi')_{L^2}$, then $R_H^{-1}=(I-\Delta)^m$, so $R_H\ell=(I-\Delta)^{-m}\ell$ has Fourier coefficients $\hat\ell(n)/(1+|n|^2)^m$ and $\lVert \ell\rVert_{H^*}=\lVert \ell\rVert_{H^{-m}}$.
\item[(c)] If $H$ is a reproducing kernel Hilbert space with continuous kernel $\mathsf k$, put $\mathsf k_x:=\mathsf k(\cdot,x)$. Then $(R_H\ell)(x)=\ell(\mathsf k_x)$ for every $x\in E$. If $\ell$ is integration against a finite signed measure $\xi$, this gives $R_H\ell=\int_E\mathsf k(\cdot,y)\,\xi(dy)$ and
\[
 \lVert\ell\rVert_{H^*}^2
 =\int_E\!\int_E\mathsf k(y,y')\,\xi(dy)\,\xi(dy').
\]
\item[(d)] If $\phi_1,\ldots,\phi_M$ form a basis of $H$ with Gram matrix $\mathbf H_{jj'}=(\phi_j,\phi_{j'})_H$, then $R_H\ell=\sum_jc_j\phi_j$ with $\mathbf H\mathbf c=\mathbf b(\ell)$, $\mathbf b(\ell)_j:=\ell(\phi_j)$, and $\lVert \ell\rVert_{H^*}^2=\mathbf b(\ell)^\top\mathbf H^{-1}\mathbf b(\ell)$.
\end{enumerate}
\end{lemma}

\begin{proof}
In each case $(R_H\ell,\phi)_H=\pair{\phi}{\ell}$ is checked directly: in (a) and (d)
by testing against the basis, in (b) by Parseval, in (c) by the reproducing
property $(\mathsf k(\cdot,y),\phi)_H=\phi(y)$. The isometry
$\lVert \ell\rVert_{H^*}=\lVert R_H\ell\rVert_H$ gives the norms.
\end{proof}

In case (b), $R_H=(I-\Delta)^{-m}$ maps $H^{-m}$ to $H^m$.
The correction therefore belongs to $H^1(0,T;H)$, and
$\cA^*u_k\in\cG$ by \Cref{lem:elementary}. In case (c), testing
the integrated forward residual from \Cref{lem:riesz} on $\mathsf k_x$ gives
\begin{equation}\label{eq:rkhs-integrated-residual}
 (R_Hd_s)(x)
 =\int_E \mathsf k_x\,d(\mu_s-\rho)
   -\int_0^s\!\int_{E\times A}\LL_t^a\mathsf k_x(y)\,\pi_t(dy,da)\,dt .
\end{equation}
The same formula gives $R_Hd_T$ with $s=T$ and $\nu$ in place of $\mu_s$.
For a diffusion, $\LL_t^a$ differentiates the kernel section in the state variable $y$;
$d_s$ need not be a signed measure.
In case (d), which is \Cref{ex:finite-dim}, the
dual iterate is $u_k(t)=\sum_jc_{k,j}(t)\phi_j$ with
$c_k\in H^1(0,T;\R^M)$, the dual step is
\begin{equation}\label{eq:dual-step-coords}
    \mathbf c_{k+1}(t)=\mathbf c_k(t)-\tau\,\mathbf H^{-1}\Bigl[\mathbf b(d_T)+\int_t^T\mathbf b(d_s)\,ds\Bigr],
    \qquad
    \partial_t\mathbf c_{k+1}(t)=\partial_t\mathbf c_k(t)+\tau\,\mathbf H^{-1}\mathbf b(d_t),
\end{equation}
with $\mathbf b(d_s)_j=\int \phi_j\,d\mu_s-\int \phi_j\,d\rho-\int_0^s\!\int\LL^a_t\phi_j\,d\pi_t\,dt$,
and one factorization of the $M\times M$ matrix $\mathbf H$ serves every
iteration. Feasibility in this case means that the moments
$\int \phi_j\,d\mu_s$ satisfy the forward equation,
\[
    \int_E\phi_j\,d\mu_s=\int_E\phi_j\,d\rho+\int_0^s\!\!\int_{E\times A}\LL^a_t\phi_j\,d\pi_t\,dt
    \qquad(j=1,\ldots,M,\ \text{a.e. }s),
\]
together with the same identity at $s=T$ with $\nu$ in place of $\mu_s$.
When $H$ is large enough this is the standard Fokker-Planck equation. 

\begin{rmk}
\label{rmk:no-discretization}
Successive corrector updates give explicit densities relative to $\pi_0$:
\[
    \nu_k\propto \exp\Bigl(-\lambda\sum_{j<k}G^{z_{j+1/2}}\Bigr)\,\nu_0,\qquad
    \pi_{k,t}\propto e^{-\lambda\Sigma_k(t,\cdot,\cdot)}\,\pi_{0,t},\qquad
    \Sigma_k:=\sum_{j<k}\Gamma^{z_{j+1/2}},
\]
with the HJB residuals of \eqref{eq:slacks}, so $\Sigma_k$ is the accumulated running HJB residual. The dual iterates are given by \eqref{eq:time-exact} and \Cref{lem:riesz-explicit}. The only quantities that have to be evaluated numerically are integrals against the iterates, namely the normalizing constants of the tilts and the moments $\int \phi\,d\mu_{k,s}$ and $\int\LL^a_t\phi\,d\pi_{k,t}$. If $\pi_0$ is the occupation measure of a reference process $X^0$ with relaxed feedback $\alpha^0$, every such integral is an expectation under the reference law reweighted by the explicit tilt,
\[
    \int\phi\,d\pi_{k,t}
    =\frac{\mathbb E\bigl[\sum_a\alpha^0(a\mid t,X^0_t)\,\phi(X^0_t,a)\,e^{-\lambda \Sigma_k(t,X^0_t,a)}\bigr]}
          {\mathbb E\bigl[\sum_a\alpha^0(a\mid t,X^0_t)\,e^{-\lambda \Sigma_k(t,X^0_t,a)}\bigr]},
\]
and can be computed by quadrature when the reference law is explicit, as in
\Cref{sec:experiments}, or by Monte Carlo on one fixed sample of reference
paths with the time variable sampled uniformly. Feasible measures, which
are needed for the correction $\projr$ of \Cref{subsec:projection} and for
the Stampacchia gap, are laws of controlled processes. When those laws are
absolutely continuous with respect to a common reference law---for example,
for bounded drift controls with a fixed uniformly nondegenerate diffusion---
they may be evaluated by Girsanov reweighting of the reference paths.
\end{rmk}

\section{Convergence under Lasry--Lions monotonicity}\label{sec:LL}

\begin{assumption}
\label{ass:LL}
For all $\pi,\pi'\in\Pi$,
$\pair{\cC_{\cM \pi}-\cC_{\cM \pi'}}{\pi-\pi'}\ge0$.
\end{assumption}

\begin{rmk}
\label{rmk:LL}
Monotonicity is imposed on all of $\Pi$ because the iterates need not
satisfy the forward equation. It holds in the classical separable case
$f(t,y,a,\mu)=f_0(t,y,\mu)+f_1(t,y,a)$ whenever $f_0(t,\cdot,\cdot)$ and
$g$ are Lasry--Lions monotone on $\cP(E)$ \cite{lasry2007mean}, that is,
$\int(f_0(t,y,\mu)-f_0(t,y,\mu'))(\mu-\mu')(dy)\ge0$, and similarly for
$g$. Indeed,
\[
    \pair{\cC_{\cM \pi}-\cC_{\cM \pi'}}{\pi-\pi'}
    =\pair{g(\cdot,\nu)-g(\cdot,\nu')}{\nu-\nu'}
    +\int_0^T\pair{f_0(t,\cdot,\mu_t)-f_0(t,\cdot,\mu'_t)}{\mu_t-\mu'_t}\,dt\ge0 .
\]
\end{rmk}

\begin{assumption}
\label{ass:Lip}
There exist constants $L_f,L_g\ge0$, with $L_C:=\max\{L_f,L_g\}>0$, such that for every $t\in[0,T]$ and every $\mu,\mu'\in\cP(E)$,
\[
\begin{aligned}
    \sup_{y,a}|f(t,y,a,\mu)-f(t,y,a,\mu')|
    &\le L_f\lVert\mu-\mu'\rVert_1,\\
    \sup_y|g(y,\mu)-g(y,\mu')|
    &\le L_g\lVert\mu-\mu'\rVert_1.
\end{aligned}
\]
\end{assumption}

We assume $L_C>0$ because $L_C=0$ removes all mean-field interaction.
Costs that are Lipschitz in $W_1$ with constants $L_f',L_g'$ satisfy
\Cref{ass:Lip} with $L_f=\frac{D_E}2L_f'$ and $L_g=\frac{D_E}2L_g'$, since
$W_1\le\frac{D_E}{2}\lVert\cdot\rVert_1$. The total variation form is the
one used in the proofs, and it does not involve the diameter of $E$. For
measurable flows $\bm\mu,\bm\mu'$ in $\cP(E)$ we define the $L^2$-in-time
total variation norm
\begin{equation}\label{eq:flow-distance}
\nsharp{\bm\mu-\bm\mu'}^2:=\lVert\mu_T-\mu'_T\rVert_1^2+\int_0^T\lVert\mu_t-\mu'_t\rVert_1^2\,dt ,
\end{equation}
the same construction as the norm of $\cX$. It is the distance in which
\Cref{ass:Lip} makes $\bm\mu\mapsto\cC_{\bm\mu}$ Lipschitz (\Cref{prop:monotone}).

\begin{proposition}
\label{prop:monotone}
Under \Cref{ass:standing}, for all $z=(\pi,u)$ and $z'=(\pi',u')$ in $\cZ$,
\[
    \pair{\cF(z)-\cF(z')}{z-z'}=\pair{\cC_{\cM \pi}-\cC_{\cM \pi'}}{\pi-\pi'} .
\]
Hence $\cF$ is monotone under \Cref{ass:LL}. Under \Cref{ass:Lip}, for
all measurable flows $\bm\mu,\bm\mu'$ and all $\pi,\pi'\in\Pi$,
\[
    \nG{\cC_{\bm\mu}-\cC_{\bm\mu'}}\le L_C\,\nsharp{\bm\mu-\bm\mu'},
    \qquad
    \nsharp{\cM \pi-\cM \pi'}\le\nsharp{\pi-\pi'},
\]
and consequently $\nG{\cC_{\cM \pi}-\cC_{\cM \pi'}}\le L_C\nsharp{\pi-\pi'}$.
\end{proposition}

\begin{proposition}
\label{prop:eps-relations}
Let $\pi\in\Pirho$.
\begin{enumerate}
\item[(i)] Under \Cref{ass:LL}, $\GapM(\pi)\le\Expl(\pi)$.
\item[(ii)] Under \Cref{ass:Lip}, suppose $\GapM(\pi)\le\varepsilon$ with
$0\le\varepsilon\le4L_C(1+T)$.\newline
Then $\Expl(\pi)\le4\sqrt{(1+T)L_C\,\varepsilon}$.
In particular, if $\GapM(\pi)=0$, then $\pi$ is an MFE.
\end{enumerate}
\end{proposition}

\subsection{Three-point identity and one-step estimate}\label{subsec:PD-one-step}

\begin{lemma}
\label{lem:PD-three-point}
Let $z=(\pi,u)\in\cZ$, $(\gamma,\ell)\in\cG\times\HY^*$ and
$z^+=\prox_z(\gamma,\ell)$. Then $z^+\in\cZ$ and $\bD(\pi^+\|\pi)<\infty$. For every
$z'=(\pi',u')\in\cZ$ with $\bD(\pi'\|\pi)<\infty$ we have
$\bD(\pi'\|\pi^+)<\infty$ and
\[
    \pair{\gamma}{\pi^+-\pi'}+\pair{\ell}{u^+-u'}
    =\Phi(z',z)-\Phi(z',z^+)-\Phi(z^+,z).
\]
In particular, $z^+$ is the unique minimizer over $\cZ$ of
$z'\mapsto\pair{\gamma}{\pi'}+\pair{\ell}{u'}+\Phi(z',z)$.
\end{lemma}

\begin{lemma}
\label{lem:PD-one-step}
Suppose Assumptions~\ref{ass:standing} and~\ref{ass:Lip} hold, and let
$(z_k,z_{k+1/2})_{k\ge0}$ be generated by \eqref{eq:PD-KL-MP}. Let
$z'=(\pi',u')\in\cZ$ with $\bD(\pi'\|\pi_0)<\infty$. Then
$\bD(\pi'\|\pi_k)<\infty$ and $\bD(\pi'\|\pi_{k+1/2})<\infty$ for every $k$, and
\[
    \pair{\cF(z_{k+1/2})}{z_{k+1/2}-z'}
    \le\Phi(z',z_k)-\Phi(z',z_{k+1})
    -(1-\theta)\bigl[\Phi(z_{k+1/2},z_k)+\Phi(z_{k+1},z_{k+1/2})\bigr],
\]
where $\theta:=\lambda L_C+\sqrt{\lambda\tau}\,L_A$.
\end{lemma}

The two terms in $\theta$ control the change in the mean-field costs
and the coupling through the forward equation. The factor
$\lambda L_C$ comes from the mean-field dependence of the costs, exactly as
$\lambda L_mD_\Omega$ in \cite{KLMP}. The factor $\sqrt{\lambda\tau}L_A$
comes from the bilinear coupling through $\cA$, as in Mirror-Prox for
bilinear saddle-point problems \cite{nemirovski2004prox}. The step-size condition is
\begin{equation}\label{eq:step-size}
    \theta:=\lambda L_C+\sqrt{\lambda\tau}\,L_A\le1 .
\end{equation}

\subsection{Minty estimate for the averages and the comparison set}\label{subsec:PD-ergodic}

\begin{thm}
\label{thm:PD-ergodic-minty}
Suppose Assumptions~\ref{ass:standing}, \ref{ass:LL} and~\ref{ass:Lip} hold,
and let \eqref{eq:step-size} hold. Then, for every $z'=(\pi',u')\in\cZ$
with $\bD(\pi'\|\pi_0)<\infty$ and every $K\ge1$,
\[
    \pair{\cF(z')}{\bar z_K-z'}
    \le\frac{1}{K}\Bigl(\frac{\bD(\pi'\|\pi_0)}{\lambda}+\frac{\tn{u'-u_0}^2}{2\tau}\Bigr).
\]
\end{thm}

{\begin{rmk}
\label{rmk:comparison-set}
For $z'=(\pi',u')$, the definition of
$\cF$ gives
\begin{equation}\label{eq:gap-expansion}
    \pair{\cF(z')}{\bar z_K-z'}
    =\pair{\cC_{\cM \pi'}}{\bar \pi_K-\pi'}-\pair{u'}{\cA\bar\pi_K-\rho}+\pair{\cA \pi'-\rho}{\bar u_K}.
\end{equation}
Two observations follow.
\begin{enumerate}
\item[(a)] \emph{Primal comparisons.} The gaps of \Cref{def:MFE} use only
feasible comparisons $\pi'\in\Pirho$. For these the last term of
\eqref{eq:gap-expansion} vanishes, and \eqref{eq:gap-expansion} splits
into the \emph{primal Minty term} $\pair{\cC_{\cM \pi'}}{\bar \pi_K-\pi'}$
and the \emph{feasibility term} $-\pair{u'}{\cA\bar\pi_K-\rho}$.
\item[(b)] \emph{Entropy.} As in \cite{KLMP}, the tilts cannot create mass
outside $\supp \pi_{0,t}$, and the bound of \Cref{thm:PD-ergodic-minty}
is finite only for comparisons with $\bD(\pi'\|\pi_0)<\infty$. We
therefore compare against
\[
    \Pizero:=\{\pi\in\Pirho:\bD(\pi\|\pi_0)<\infty\}.
\]
If $E$ and $A$ are finite and $\nu_0$ and $\pi_{0,t}$ give every point
mass bounded below uniformly in $t$, then $\Pizero=\Pirho$.
\end{enumerate}
\end{rmk}}

\begin{corollary}
\label{cor:primal-minty-residual}
Under the assumptions of \Cref{thm:PD-ergodic-minty}, for every $K\ge1$:
\begin{enumerate}[label=(\alph*)]
\item (Minty estimates.) For every $\pi'\in\Pizero$,
\begin{equation}\label{eq:master}
\begin{aligned}
    \pair{\cC_{\cM\pi'}}{\bar\pi_K-\pi'}
    &\le \frac{\bD(\pi'\|\pi_0)}{\lambda K}
    +\pair{\cA\bar\pi_K-\rho}{u_0}
    -\frac{\tau K}{2}\tns{\cA\bar\pi_K-\rho}^2\\
    &\le \frac{\bD(\pi'\|\pi_0)}{\lambda K}
    +\frac{\tn{u_0}^2}{2\tau K}
    =O(K^{-1}).
\end{aligned}
\end{equation}
\item (Forward-equation residual.) For every
$\pi'\in\Pizero$,
\[
    \tns{\cA\bar\pi_K-\rho}\le\frac{2\tn{u_0}}{\tau K}
    +\left(\frac{2}{\tau K}\Bigl(\frac{\bD(\pi'\|\pi_0)}{\lambda K}+2C_0\sqrt{1+T}\Bigr)\right)^{1/2}
    =O(K^{-1/2}).
\]
\item (Forward-equation residual with an optimal dual function.) If there is an
MFE $\pi^*\in\Pizero$ with an optimal dual variable $u^*\in\HY$ as in
\Cref{prop:PDVI-MFE}(ii), then
\[
    \tns{\cA\bar\pi_K-\rho}\le\frac{2(\tn{u_0}+\tn{u^*})}{\tau K}+\frac{1}{K}\Bigl(\frac{2\bD(\pi^*\|\pi_0)}{\lambda\tau}\Bigr)^{1/2}
    =O(K^{-1}).
\]
\end{enumerate}
\end{corollary}

\begin{rmk}
\label{rmk:dual-drift}
The corrector update gives $u_{k+1}=u_k-\tau R(\cA \pi_{k+1/2}-\rho)$. Summing over $k=0,\ldots,K-1$ and using that $\cA$ is affine on averages yields the exact identity
\begin{equation}\label{eq:dual-drift}
    u_K-u_0=-\tau K\,R (\cA\bar\pi_K-\rho),
    \qquad\text{hence}\qquad
    \tns{\cA\bar\pi_K-\rho}=\frac{\tn{u_K-u_0}}{\tau K}.
\end{equation}
In particular, the comparison dual $u'=u_0-\tau KR (\cA\bar\pi_K-\rho)$ used in the proof of \eqref{eq:master} is the last dual iterate $u_K$. The dual update therefore records the forward-equation residual. Its average is determined by the change in the dual variable divided by $\tau K$. Thus the feasibility rate is decided by the growth of the dual
iterates.
\begin{enumerate}
\item[(i)] Suppose an optimal dual variable $u^*\in\HY$ exists for an MFE $\pi^*\in\Pizero$. Then \Cref{lem:PD-one-step} with $z'=z^*$ and monotonicity give $\Phi(z^*,z_{k+1})\le\Phi(z^*,z_k)$. So $\tn{u_K-u^*}\le(2\tau\Phi(z^*,z_0))^{1/2}$ stays bounded, and
\eqref{eq:dual-drift} gives $\tns{\cA\bar\pi_K-\rho}=O(1/K)$. This is another route to the $O(K^{-1})$ residual bound. The optimal dual function $u^*$ bounds the possible cost reduction from violating the forward equation: a residual $r$ can lower the cost by at most $\tn{u^*}\tns{r}$.
\item[(ii)] Without an optimal dual variable in $\HY$, the cost reduction obtainable from a violation of size $r$ need not be $O(r)$. Fix an MFE $\pi^*$ with $\bm\mu^*=\cM \pi^*$, and let $\mathrm{val}(r):=\inf\{\pair{\cC_{\bm\mu^*}}{\pi}:\pi\in\Pi,\ \cA \pi-\rho=r\}$ be the optimal cost when the forward equation has residual $r$. Then $\mathrm{val}(0)=\pair{\cC_{\bm\mu^*}}{\pi^*}$. An optimal dual function gives a linear lower bound for $\mathrm{val}$ at $0$. Without such a function, $\mathrm{val}(0)-\mathrm{val}(r)$ may be larger than any constant times $\tns r$ near zero. The norms of $u_k$ may then keep growing, and \eqref{eq:master} only
balances the quadratic penalty $\frac{\tau K}{2}\tns{\cA\bar\pi_K-\rho}^2$ against a bounded cost gain $2C_0\sqrt{1+T}$, which gives the general $O(K^{-1/2})$ residual bound. More generally, since $\pair{\cC_{\bm\mu^*}}{\bar \pi_K-\pi^*}\ge \mathrm{val}(\cA\bar\pi_K-\rho)-\mathrm{val}(0)$, a growth bound $\mathrm{val}(0)-\mathrm{val}(r)\lesssim\tns r^{\varsigma}$ with $\varsigma\in[0,1]$, inserted into \eqref{eq:master} with $\pi'=\pi^*$, gives $\tns{\cA\bar\pi_K-\rho}=O(K^{-1/(2-\varsigma)})$.
\end{enumerate}
\end{rmk}

\subsection{Cluster points}\label{subsec:cluster}

We identify $\pi=(\nu,(\pi_t)_t)\in\Pi$ with the pair consisting of
$\nu\in\cP(E)$ and the finite Borel measure $dt\,\pi_t(dy,da)$ on
$[0,T]\times E\times A$, give $\Pi$ the product of the weak topologies,
and write $\pi^n\Rightarrow\pi$ for convergence in this topology.

\begin{lemma}
\label{lem:topology}
Under (S1), $\Pi$ is compact and metrizable for $\Rightarrow$. Let $\pi^n\Rightarrow \pi$
in $\Pi$, and let $\Gamma:[0,T]\times E\times A\to\R$ be Borel with
$\Gamma_t\in C(E\times A)$ for a.e.\ $t$ and
$\int_0^T\lVert\Gamma_t\rVert_\infty dt<\infty$. Then $\int_0^T\!\int\Gamma_t\,d\pi^n_t\,dt\to\int_0^T\!\int\Gamma_t\,d\pi_t\,dt$.
Consequently, $\pi\mapsto\pair{\gamma}{\pi}$ is continuous on $\Pi$ for each
$\gamma\in\cG$, and $\pi\mapsto\pair{\cA \pi-\rho}{u}$ is continuous for each
$u\in\HY$.
\end{lemma}

\begin{assumption}
\label{ass:density}
For every $\pi'\in\Pirho$ there are $\pi'_n\in\Pizero$ such that
$\pi'_n\Rightarrow \pi'$ and $\nsharp{\cM \pi'_n-\cM \pi'}\to0$, with $\nsharp{\cdot}$ as in
\eqref{eq:flow-distance}.
\end{assumption}

\Cref{ass:density} holds trivially, with $\pi'_n=\pi'$, when $E$ and $A$ are
finite and $\pi_0$ has full support with every point mass bounded below uniformly in $t$.
For general compact spaces, it requires that every feasible occupation
measure can be approximated by feasible measures with finite entropy
relative to $\pi_0$, with their population flows converging in
$L^2(dt;\lVert\cdot\rVert_1)$.

\begin{thm}
\label{thm:cluster-points}
Suppose Assumptions~\ref{ass:standing}, \ref{ass:LL} and~\ref{ass:Lip} and
\eqref{eq:step-size} hold, and suppose $\Pizero\neq\emptyset$. Then:
\begin{enumerate}
\item[(i)] $(\bar \pi_K)_{K\ge1}$ has $\Rightarrow$-limit points, and every
limit point $\bar \pi$ satisfies $\bar \pi\in\Pirho$ and
\[
    \pair{\cC_{\cM \pi'}}{\bar \pi-\pi'}\le0\qquad\forall \pi'\in\Pizero.
\]
\item[(ii)] If in addition \Cref{ass:density} holds, every limit point is
an MFE. If moreover the MFE is unique, then $\bar \pi_K\Rightarrow \pi^*$.
\end{enumerate}
\end{thm}

This result does not require an optimal dual function in $\HY$. The averaged dual variables
$\bar u_K$ may diverge, while the primal averages still converge to
equilibria. In this sense the scheme keeps the advantage of the primal--dual characterization discussed in \cite{guo2025continuous}: it
needs no convex Hamiltonian, no unique minimizer, and no selection of a
feedback. As in \cite[Rem.~3.1]{KLMP}, the theorem also gives a constructive existence proof of an MFE.

\subsection{Correcting the forward equation and approximate equilibria}\label{subsec:projection}

The average $\bar \pi_K$ need not belong to $\Pirho$. We therefore correct it
to a point satisfying the forward equation for the chosen test functions.
For the total-variation estimate in \Cref{thm:quantitative-eps-MFE},
we assume that the size of this correction is controlled by the
forward-equation residual.

\begin{assumption}
\label{ass:projection}
There are a map $\projr:\Pi\to\Pirho$ and a constant $C_{\mathrm{proj}}\ge0$ such that
\[
    \nsharp{\projr(\pi)-\pi}\le C_{\mathrm{proj}}\,\tns{\cA \pi-\rho}\qquad\forall \pi\in\Pi .
\]
\end{assumption}

In \Cref{prop:cost-correction}, we give a weaker sufficient condition for the approximate-equilibrium bound.

\begin{assumption}
\label{ass:entropy-bound}
There is $\bar D<\infty$ with $\bD(\pi\|\pi_0)\le\bar D$ for all $\pi\in\Pirho$.
\end{assumption}

\Cref{ass:entropy-bound} is a strong restriction: it fails whenever a feasible occupation measure is singular with respect to the reference, as can happen with deterministic controls in a continuous action space. It is met by finite state and action spaces with every point
mass bounded below uniformly in $t$, see \Cref{thm:finite-state}, where \Cref{ass:projection} is verified as well.

\begin{thm}
\label{thm:quantitative-eps-MFE}
Suppose Assumptions~\ref{ass:standing}, \ref{ass:LL}, \ref{ass:Lip}, \ref{ass:projection} and~\ref{ass:entropy-bound} and \eqref{eq:step-size} hold. For $K\ge1$ set $\hat \pi_K:=\projr(\bar \pi_K)\in\Pirho$. Then:
\begin{enumerate}
\item[(i)] $\tns{\cA\bar\pi_K-\rho}\le\omega_K:=\dfrac{2(\tn{u_0}+C_0C_{\mathrm{proj}})}{\tau K}+\dfrac1K\Bigl(\dfrac{2\bar D}{\lambda\tau}\Bigr)^{1/2}$.
\item[(ii)] $\nsharp{\hat \pi_K-\bar \pi_K}\le C_{\mathrm{proj}}\omega_K$.
\item[(iii)] $\GapM(\hat \pi_K)\le\varepsilon_K:=\dfrac{\bar D}{\lambda K}+\dfrac{\tn{u_0}^2}{2\tau K}+C_0C_{\mathrm{proj}}\omega_K$.
\item[(iv)] If $\varepsilon_K\le 4L_C(1+T)$, then
$\Expl(\hat \pi_K)\le4\sqrt{(1+T)L_C\,\varepsilon_K}$.
\end{enumerate}
In particular, for $u_0=0$,
\[
    \varepsilon_K\le\frac1K\Bigl(\frac{2\bar D}{\lambda}+\frac{5C_0^2C_{\mathrm{proj}}^2}{2\tau}\Bigr),
\]
so $\Expl(\hat\pi_K)=O(K^{-1/2})$. For every $\delta>0$, we have $\Expl(\hat\pi_K)\le\delta$ whenever
\[
    K\ge\Bigl(\frac{2\bar D}{\lambda}+\frac{5C_0^2C_{\mathrm{proj}}^2}{2\tau}\Bigr)
    \max\Bigl\{\frac{1}{4L_C(1+T)},\,\frac{16(1+T)L_C}{\delta^2}\Bigr\}.
\]
\end{thm}

For the approximate-equilibrium conclusion, it is sufficient to control how the correction changes the costs evaluated at feasible comparison populations.

\begin{proposition}\label{prop:cost-correction}
Suppose \Cref{ass:standing,ass:LL,ass:Lip,ass:entropy-bound}
and \eqref{eq:step-size} hold.
Assume there are a map $\projr:\Pi\to\Pirho$, a constant
$C_{\mathrm c}\ge0$ and an exponent $\alpha\in(0,1]$ such that
\begin{equation}\label{eq:cost-correction}
 \sup_{q\in\Pirho}
 \left|\pair{\cC_{\cM q}}{\projr(\pi)-\pi}\right|
 \le C_{\mathrm c}\,\tns{\cA\pi-\rho}^{\alpha},
 \qquad \pi\in\Pi.
\end{equation}

For $K\ge1$, set $\hat\pi_K:=\projr(\bar\pi_K)$ and define
\[
\begin{aligned}
 \omega_K&:=\max\left\{
   \frac1K\sqrt{\frac{6\bar D}{\lambda\tau}},\,
   \frac{6\tn{u_0}}{\tau K},\,
   \left(\frac{6C_{\mathrm c}}{\tau K}\right)^{1/(2-\alpha)}
   \right\},\\
 \varepsilon_K&:=\frac{\bar D}{\lambda K}
       +\frac{\tn{u_0}^2}{2\tau K}
       +C_{\mathrm c}\omega_K^{\alpha}.
\end{aligned}
\]
Then
\[
 \tns{\cA\bar\pi_K-\rho}\le\omega_K,
 \qquad
 \GapM(\hat\pi_K)\le\varepsilon_K.
\]
If $\varepsilon_K\le4L_C(1+T)$, then $\Expl(\hat\pi_K)\le4\sqrt{(1+T)L_C\varepsilon_K}$. In particular,
\[
 \tns{\cA\bar\pi_K-\rho}=O\bigl(K^{-1/(2-\alpha)}\bigr),\quad
 \GapM(\hat\pi_K)
       =O\bigl(K^{-\alpha/(2-\alpha)}\bigr),\quad
 \Expl(\hat\pi_K)
       =O\bigl(K^{-\alpha/[2(2-\alpha)]}\bigr).
\]
\end{proposition}

Assumption~\ref{ass:projection} implies \eqref{eq:cost-correction}
with $\alpha=1$ and $C_{\mathrm c}=C_0C_{\mathrm{proj}}$.
The proof is in \Cref{subsec:proof-cost-correction}. For $\alpha=1$,
the forward-equation residual and Minty gap are $O(K^{-1})$, and
the corrected average is an $O(K^{-1/2})$-MFE.

\begin{rmk}\label{rmk:wasserstein-correction}
A Wasserstein bound provides a sufficient condition for
\eqref{eq:cost-correction}. Equip $E\times A$ with the product distance
$d_E(y,y')+d_A(a,a')$, and for
$\pi=(\nu,(\pi_t)_t)$ and $\pi'=(\nu',(\pi'_t)_t)$ define
\begin{equation}\label{eq:wasserstein-pair-distance}
 d_W(\pi,\pi')^2
 :=W_1(\nu,\nu')^2+
       \int_0^T W_{1,E\times A}(\pi_t,\pi'_t)^2\,dt .
\end{equation}
Suppose there are $\Lambda_g\ge0$ and a nonnegative
$\Lambda_f\in L^2(0,T)$ such that, uniformly in $\mu\in\cP(E)$,
\[
 \operatorname{Lip}_y g(\cdot,\mu)\le\Lambda_g,
 \qquad
 \operatorname{Lip}_{(y,a)} f(t,\cdot,\cdot,\mu)\le\Lambda_f(t)
 \quad\text{for a.e. }t.
\]
Kantorovich--Rubinstein duality and Cauchy--Schwarz give
\[
\begin{aligned}
 \sup_{q\in\Pirho}
 \left|\pair{\cC_{\cM q}}{\pi-\pi'}\right|
 &\le\Lambda_g W_1(\nu,\nu')
      +\int_0^T\Lambda_f(t)
               W_{1,E\times A}(\pi_t,\pi'_t)\,dt\\
 &\le\left(\Lambda_g^2+\int_0^T\Lambda_f(t)^2\,dt\right)^{1/2}
          d_W(\pi,\pi').
\end{aligned}
\]
Consequently, if $\projr:\Pi\to\Pirho$ satisfies
\begin{equation}\label{eq:wasserstein-correction}
 d_W(\projr(\pi),\pi)
 \le C_W\,\tns{\cA\pi-\rho}^{\alpha},
 \qquad \pi\in\Pi,
\end{equation}
for some $C_W\ge0$ and $\alpha\in(0,1]$, then the cost estimate
\eqref{eq:cost-correction} holds with
$C_{\mathrm c}=C_W(\Lambda_g^2+\int_0^T\Lambda_f(t)^2\,dt)^{1/2}$.
The estimates of \Cref{prop:cost-correction} therefore apply.
\end{rmk}

\subsection{Wasserstein convergence under strong Lasry--Lions monotonicity}\label{subsec:strong-LL}

The results so far control the forward-equation residual of
$\bar\pi_K$ and the Minty and Stampacchia gaps of its corrected
average $\hat\pi_K$. They do not say how far the iterates are from an equilibrium. Under a strong form of \Cref{ass:LL}, the one-step estimate of
\Cref{lem:PD-one-step} controls the squared distance between the
iterated population flow and the equilibrium flow, as in
\cite[Section~4]{KLMP} for static games. For measurable flows $\bm\mu,\bm\mu'$ in $\cP(E)$ set
\begin{equation}\label{eq:flat-norm}
    \lVert\bm\mu-\bm\mu'\rVert_{W_1,2}^2
    :=\int_0^TW_1(\mu_t,\mu'_t)^2\,dt+W_1(\mu_T,\mu'_T)^2 .
\end{equation}
The integrand is measurable, since $t\mapsto\mu_t$ is a measurable map
into $\cP(E)$ with the weak topology and $W_1$ is continuous there. By
Kantorovich--Rubinstein duality, $W_1(\mu_t,\mu'_t)$ is the norm of
$\mu_t-\mu'_t$ in the dual of the $1$-Lipschitz functions, so
\eqref{eq:flat-norm} is the squared $L^2$-in-time norm of the signed flow
$\bm\mu-\bm\mu'$; in particular $\lVert\cdot\rVert_{W_1,2}$ satisfies the
triangle inequality and $\pi\mapsto\lVert\cM \pi-\bm\mu'\rVert_{W_1,2}^2$ is
convex on $\Pi$, because $\cM$ is affine. Since
$W_1\le\frac{D_E}{2}\lVert\cdot\rVert_1$, \Cref{prop:monotone} gives
\begin{equation}\label{eq:flat-vs-d1}
    \lVert\cM \pi-\cM \pi'\rVert_{W_1,2}\le\tfrac{D_E}{2}\,\nsharp{\cM \pi-\cM \pi'}\le\tfrac{D_E}{2}\nsharp{\pi-\pi'}
    \qquad(\pi,\pi'\in\Pi).
\end{equation}

\begin{assumption}
\label{ass:strong-LL}
There is $\kappa>0$ such that, for all $\pi,\pi'\in\Pi$,
\[
    \pair{\cC_{\cM \pi}-\cC_{\cM \pi'}}{\pi-\pi'}\ge\kappa\,\lVert\cM \pi-\cM \pi'\rVert_{W_1,2}^2 .
\]
\end{assumption}

\Cref{ass:strong-LL} implies \Cref{ass:LL}. In the separable case of
\Cref{rmk:LL} it holds as soon as $f_0(t,\cdot,\cdot)$ and $g$ are
$\kappa$-strongly monotone in $W_1$ on $\cP(E)$, that is,
$\int(f_0(t,y,\mu)-f_0(t,y,\mu'))(\mu-\mu')(dy)\ge\kappa W_1(\mu,\mu')^2$
and likewise for $g$: the identity in \Cref{rmk:LL} then produces exactly
the right side of \eqref{eq:flat-norm}. On a finite state space with the
discrete metric, $W_1=\frac12\lVert\cdot\rVert_1$, and the congestion costs
$f_0(t,y,\mu)=c_f\,\mu(y)$ and $g(y,\mu)=c_g\mu(y)$ with
$c_f,c_g>0$ satisfy it with $\kappa=4\min\{c_f,c_g\}/N$,
since $c_f\lVert\mu-\mu'\rVert_2^2\ge\frac{c_f}{N}\lVert\mu-\mu'\rVert_1^2$.
If only the running cost is strongly monotone, as in
\Cref{subsec:exp-finite}, everything below holds with the terminal term
dropped from \eqref{eq:flat-norm}.

\begin{thm}
\label{thm:strong-LL}
Suppose Assumptions~\ref{ass:standing}, \ref{ass:Lip} and~\ref{ass:strong-LL} and \eqref{eq:step-size} hold, and let $z^*=(\pi^*,u^*)\in\cZ$ solve the PD-VI with $\bD(\pi^*\|\pi_0)<\infty$. When both $E$ and $A$ are finite and the reference masses satisfy the uniform lower bounds of \Cref{thm:finite-state}, this holds for every MFE $\pi^*$ with $u^*=V$, by parts (a) and (b) of that theorem. Write
$\bm\mu^*:=\cM \pi^*$. Then:
\begin{enumerate}[label=(\roman*)]
\item (Descent inequality.) For every $k\ge0$,
\[
    \Phi(z^*,z_{k+1})
    \le\Phi(z^*,z_k)-\kappa\,\lVert\cM \pi_{k+1/2}-\bm\mu^*\rVert_{W_1,2}^2
    -(1-\theta)\bigl[\Phi(z_{k+1/2},z_k)+\Phi(z_{k+1},z_{k+1/2})\bigr].
\]
In particular $(\Phi(z^*,z_k))_{k\ge0}$ is nonincreasing and $\tn{u_k-u^*}\le(2\tau\Phi(z^*,z_0))^{1/2}$ for all $k$.
\item (Rates.) $\sum_{k\ge0}\kappa\lVert\cM \pi_{k+1/2}-\bm\mu^*\rVert_{W_1,2}^2\le\Phi(z^*,z_0)$,
and hence, for every $K\ge1$,
\[
    \lVert\cM\bar \pi_K-\bm\mu^*\rVert_{W_1,2}^2\le\frac{\Phi(z^*,z_0)}{\kappa K},
    \qquad
    \min_{0\le k<K}\lVert\cM \pi_{k+1/2}-\bm\mu^*\rVert_{W_1,2}^2\le\frac{\Phi(z^*,z_0)}{\kappa K}.
\]
\item (Convergence of the whole sequence.) $\lVert\cM \pi_{k+1/2}-\bm\mu^*\rVert_{W_1,2}\to0$ as $k\to\infty$. If
$\theta<1$, then also $\nsharp{\pi_{k+1/2}-\pi_k}\to0$ and
$\lVert\cM \pi_k-\bm\mu^*\rVert_{W_1,2}\to0$.
\item (Uniqueness of the equilibrium flow.) Under Assumptions~\ref{ass:standing} and~\ref{ass:strong-LL} alone, any two MFE
$\pi^*,\pi^{**}\in\Pirho$ have the same flow, $\cM \pi^*=\cM \pi^{**}$, that is,
$\mu^*_t=\mu^{**}_t$ for a.e.\ $t$ and $\mu^*_T=\mu^{**}_T$. In
particular, the limit flow $\bm\mu^*$ in (ii)--(iii) does not depend on the
chosen PD-VI solution.
\end{enumerate}
\end{thm}

The proof is in \Cref{sec:proofs-LL}. The constant $\Phi(z^*,z_0)$ is the
initial distance that also governs \Cref{thm:PD-ergodic-minty}; on a
finite state space with $u_0=0$ it is at most
$\bar D/\lambda+\tn{V}^2/(2\tau)$, with $\bar D$ from \Cref{thm:finite-state}(a).

\begin{rmk}
\label{rmk:no-linear}
\Cref{thm:strong-LL} gives an inequality but no geometric rate, and none should be expected from \Cref{ass:strong-LL} alone. By
\Cref{prop:monotone},
$\pair{\cF(z)-\cF(z')}{z-z'}=\pair{\cC_{\cM \pi}-\cC_{\cM \pi'}}{\pi-\pi'}$ does
not involve $u$ at all: the component $\cA\pi-\rho$ of $\cF$ does not depend on $u$.
Thus (i) bounds $\tn{u_k-u^*}^2/(2\tau)$ but does not show that
this quantity decreases by a fixed factor at each iteration. Moreover, the decrease in (i) is
measured in $\lVert\cM \pi_{k+1/2}-\bm\mu^*\rVert_{W_1,2}$, which controls only
the flow, that is, the state marginals and the terminal law, and neither
the state--action measures $\pi_{k+1/2,t}$ themselves (hence not the
policy) nor $u_{k+1/2}$.
\end{rmk}

\subsection{Convergence rates for diffusions with convex control costs}\label{subsec:convex-diffusion}

Convexity of the control cost also gives convergence in Wasserstein distance
under ordinary Lasry--Lions monotonicity. We adapt the comparison of
control costs and stability estimates used by Bertucci and Cecchin
\cite[Section~3.2]{BertucciCecchin2024}. Their argument compares
solutions of the MFG systems; here \eqref{eq:master} supplies the
comparison estimate for the Mirror-Prox average, which may violate
the forward equation. A backward linear equation controls this residual in
the metric \eqref{eq:flat-norm}. The following assumptions specify the
regularity and the continuous action space used in this argument.

\begin{assumption}\label{ass:convex-diffusion}
\leavevmode
\begin{enumerate}[label=(C\arabic*)]
\item The state space is the unit torus $E=\mathbb T=\R/\mathbb Z$,
$A=[a_-,a_+]\subset\R$ with $a_-<a_+$, and
\[
    \LL^a\phi=a\,\partial_y\phi+\tfrac12\sigma^2\,\partial_{yy}\phi,
    \qquad \sigma>0.
\]
We take $H=H^m(\mathbb T)$ for an integer $m\ge3$, with its usual
Sobolev norm. Lebesgue measure $dy$ on $\mathbb T$ has mass one.
\item The running cost is separable as in \Cref{rmk:LL},
$f(t,y,a,\mu)=f_0(t,y,\mu)+f_1(t,y,a)$.
The function $f_1$ is twice continuously differentiable in $a$, and
there are constants $0<c_{f_1}\le C_{f_1}<\infty$ such that
\[
    c_{f_1}\le\partial_{aa}f_1(t,y,a)\le C_{f_1}
    \qquad (t,y,a)\in[0,T]\times\mathbb T\times A.
\]
The costs satisfy \Cref{ass:standing,ass:LL,ass:Lip}.
\item There are a narrowly continuous probability flow
$\bm\mu^*=(\mu_t^*)_{0\le t\le T}$, a function
$u^*\in C^{1,2}([0,T]\times\mathbb T)\cap\HY$, and a feedback
$a^*\in C([0,T];W^{m+1,\infty}(\mathbb T))$ such that
\begin{equation}\label{eq:convex-HJB}
\begin{aligned}
 \partial_tu^*+\tfrac12\sigma^2\partial_{yy}u^*+f_0(t,y,\mu_t^*)
   +\min_{a\in A}\{f_1(t,y,a)+a\,\partial_yu^*\}&=0,\\
 u^*(T,y)&=g(y,\mu_T^*),
\end{aligned}
\end{equation}
$a^*(t,y)$ is the minimizer in \eqref{eq:convex-HJB}, and
\[
 \pi^*:=\bigl(\mu_T^*,(\mu_t^*(dy)\delta_{a^*(t,y)}(da))_t\bigr)
 \in\Pirho .
\]
For some $\delta_0\in(0,1]$,
$[a^*(t,y)-\delta_0,a^*(t,y)+\delta_0]\subset A$ for every $(t,y)$.
The measures $\mu_t^*$ have densities $p^*(t,\cdot)$ with
$0\le p^*(t,y)\le M$ for all $t\in[0,T]$ and a.e.\ $y$, for a
constant $M<\infty$.
\item The algorithmic reference has densities
$\nu_0(dy)=w_T(y)\,dy$ and
$\pi_{0,t}(dy,da)=w_t(y,a)\,dy\,da$, where
$w_T\ge c_0$ and $w_t\ge c_0$ a.e., uniformly in $t$, for some $c_0>0$.
\end{enumerate}
\end{assumption}

The HJB equation and feasibility in (C3) make $\pi^*$ an MFE.

\begin{thm}\label{thm:convex-diffusion}
Suppose \Cref{ass:convex-diffusion} holds. Run the exact
continuous-time algorithm with fixed $\lambda,\tau>0$ satisfying
\eqref{eq:step-size}, and any fixed $u_0\in\HY$.
For its predictor average $\bar\pi_K=(\bar\nu_K,(\bar\pi_{K,t})_t)$,
define the squared control error
\[ E_K:=\int_0^T\!\!\int_{\mathbb T\times A}|a-a^*(t,y)|^2\,\bar\pi_{K,t}(dy,da)\,dt. \]
There is a constant $C$, independent of $K$, such that for every $K\ge1$,
\begin{equation}\label{eq:convex-energy-rate}
 E_K\le C\frac{\log(eK)}K,\qquad
 \tns{\cA\bar\pi_K-\rho}\le C\frac{\sqrt{\log(eK)}}K,
\end{equation}
and
\begin{equation}\label{eq:convex-flow-rate}
 \lVert\cM\bar\pi_K-\bm\mu^*\rVert_{W_1,2}
 \le C\left[
       \sqrt{\frac{\log(eK)}K}
       +\left(\frac{\sqrt{\log(eK)}}K\right)^{1/(m+1)}
       \right].
\end{equation}
The same right-hand side bounds
\begin{equation}\label{eq:convex-occupation-distance}
 \left[
 W_1(\bar\nu_K,\mu_T^*)^2+
 \int_0^T W_{1,\mathbb T\times A}(\bar\pi_{K,t},\pi_t^*)^2\,dt
 \right]^{1/2},
\end{equation}
where the product distance is
$d_{\mathbb T}(y,y')+|a-a'|$.
The constants may depend on the fixed model data, the bounds in
\Cref{ass:convex-diffusion}, $\lambda,\tau$, and $\tn{u^*-u_0}$.
\end{thm}

The proof, given in \Cref{subsec:proof-convex-diffusion}, uses feasible
comparisons obtained by spreading each optimal action over a symmetric
interval. Their conditional mean action remains $a^*$, so their
population flow is exactly $\bm\mu^*$. Their entropy grows logarithmically
as the interval shrinks. For $m=3$, \eqref{eq:convex-flow-rate} is
$O(K^{-1/4}(\log(eK))^{1/8})$.

\section{Finite-state games and grid approximations}\label{sec:finite}

We return to \Cref{ex:finite}: $E=\{1,\ldots,N\}$, $A$ compact, jump rates
$Q_t(y,y',a)$ with $\bar q<\infty$, $H=\R^E$ Euclidean, $c_H=1$,
$\ell_H=2\bar q$ and $L_A=(2+8\bar q^2T(1+T))^{1/2}$. For a relaxed
feedback $\alpha=(\alpha_t(\cdot\mid y))$ write
$Q^{\alpha}_t(y,y'):=\int_AQ_t(y,y',a)\,\alpha_t(da\mid y)$ for $y'\ne y$,
with diagonal $-\sum_{y'\ne y}Q^\alpha_t(y,y')$, and $(\mu Q)(y'):=\sum_y\mu(y)Q(y,y')$.

\begin{thm}
\label{thm:finite-state}
Let \Cref{ass:standing} hold with $E$ and $A$ finite, and let $\pi_0=(\nu_0,(\pi_{0,t})_t)$
have full support uniformly in $t$: 
$\underline\nu:=\min_y\nu_0(y)>0$ and $\underline\pi:=\operatorname{ess\,inf}_t\min_{y,a}\pi_{0,t}(y,a)>0$.
Then:
\begin{enumerate}[label=(\alph*)]
\item (Entropy bound.) $\Pizero=\Pirho\neq\emptyset$ and
\Cref{ass:entropy-bound} holds with $\bar D=\log(1/\underline\nu)+T\log(1/\underline\pi)$.
\item (Dual optimizers.) Let $\pi^*\in\Pirho$ be an MFE and $\bm\mu^*=\cM \pi^*$. The HJB system
\begin{equation}\label{eq:hjb-finite}
    -\dot V(t,y)=\min_{a\in A}\Bigl[f(t,y,a,\mu^*_t)+\sum_{y'\ne y}Q_t(y,y',a)\bigl(V(t,y')-V(t,y)\bigr)\Bigr],
    \qquad V(T,\cdot)=g(\cdot,\mu^*_T),
\end{equation}
has a unique absolutely continuous solution, $V\in W^{1,\infty}(0,T;\R^E)\subset\HY$,
and $u^*:=V$ is an optimal dual variable for $\pi^*$ in the sense of
\Cref{prop:PDVI-MFE}(ii). Hence $(\pi^*,V)$ solves the PD-VI, and
\Cref{cor:primal-minty-residual}(c) applies to every MFE.
\item (Correction to an admissible flow.) For $\pi=(\nu,(\pi_t)_t)\in\Pi$ let $\mu_t=\pi_t(\cdot\times A)$,
let $\alpha_t(\cdot\mid y):=\pi_t(\{y\}\times\cdot)/\mu_t(y)$ when
$\mu_t(y)>0$ and an arbitrary fixed law otherwise, let $\tilde\mu$ be the
solution of $\dot{\tilde\mu}_t=\tilde\mu_tQ^{\alpha}_t$, $\tilde\mu_0=\rho$,
and set $\projr(\pi):=(\tilde\mu_T,(\tilde\mu_t\otimes\alpha_t)_t)$. Then
$\projr(\pi)\in\Pirho$ and \Cref{ass:projection} holds with
\[
    C_{\mathrm{proj}}=\sqrt N\,\bigl(1+c^2+(1+2\bar qT)^2\bigr)^{1/2},
    \qquad c:=2\bar q\sqrt T\,(1+2\bar qT).
\]
\end{enumerate}
Consequently, under \Cref{ass:LL,ass:Lip} and the step-size condition
\eqref{eq:step-size}, \Cref{thm:cluster-points,thm:quantitative-eps-MFE}
hold without further hypotheses, and the forward-equation residual of the
averaged iterate satisfies the $O(K^{-1})$ bound of
\Cref{cor:primal-minty-residual}(c) for every MFE $\pi^*$.
\end{thm}

The proof is in \Cref{subsec:proof-finite}. Part (c) keeps the feedback extracted from $\pi$ and solves its forward equation from $\rho$. The bound says that if $\pi$ nearly satisfies the forward equation, this correction
changes $\pi$ only a little. Applied to $\bar \pi_K$, it gives a relaxed Markov policy whose
occupation measure is an approximate equilibrium with the explicit
guarantee of \Cref{thm:quantitative-eps-MFE}.

\subsection{Markov chain approximations of diffusions}\label{subsec:mesh}

Bayraktar, Budhiraja and Cohen \cite{BayraktarBudhirajaCohen2018Numerical}
use Markov-chain approximation to compute equilibria of diffusion MFGs
with reflecting boundaries, arising from the queueing models of
\cite{BayraktarBudhirajaCohen2019RateControl}. Here the chain approximates periodic diffusion
dynamics, and its equilibrium is computed by the primal--dual Mirror-Prox
method rather than by repeatedly solving control problems for fixed population flows.

\paragraph{The diffusion model.}
We consider the one-dimensional diffusion game on
$E=\mathbb T_L:=\R/2L\mathbb Z$. For an $A$-valued control $\alpha$, the
state dynamics are
\[
    dX_t=b(t,X_t,\alpha_t)\,dt+\sigma(t,X_t,\alpha_t)\,dW_t,
    \qquad X_0\sim\rho,
\]
where $W$ is a one-dimensional Brownian motion and the coefficients are
$2L$-periodic in the state. The generator is
$\LL^a_t\phi=b(t,y,a)\phi'+\tfrac12\sigma^2(t,y,a)\phi''$, with $b,\sigma^2$
bounded, Borel in $t$ and continuous in $(y,a)$. The population enters the
costs $f$ and $g$, not the dynamics. We use the relaxed Markov formulation:
under $\alpha_t(da\mid y)$, the generator is
$\LL^\alpha_t\phi(y)=\int_A\LL^a_t\phi(y)\,\alpha_t(da\mid y)$ and the player
minimizes \eqref{eq:game-cost} for a fixed population flow $\bm\mu$.
At equilibrium the control is optimal and $\mu_t=\Law(X^\alpha_t)$, with
the process interpretation subject to \Cref{rmk:MFE-process}.
The approximation below discretizes the state space, not time: it
replaces the diffusion by a continuous-time Markov chain on a finite grid.

For the chain on a grid of mesh $h$, the Euclidean norm on $H=\R^E$
gives $\ell_H=2\bar q$ as in \Cref{ex:finite}. When $\bar q$ is of order
$h^{-2}$, so is $L_A$, and the step-size condition requires $\lambda\tau\lesssim h^4$.
To obtain admissible step sizes independent of $h$, we change the
inner product on $\R^E$. We use the grid version
of the Sobolev inner product in \Cref{lem:riesz-explicit}(b): it measures
finite differences as well as function values and is defined in
\eqref{eq:Hm} below. With this inner product the constants $c_H$, $\ell_H$,
$L_A$, and with them the
step-size condition \eqref{eq:step-size}, are independent of $h$.
For a fixed finite action set and uniform references the entropy bound is
$O(\log(1/h))$. The dual step is one backward integral in time and one FFT
solve per time slice. This is the grid analogue of
\Cref{ex:diffusion}, and it makes precise in which sense the method is
mesh-independent.

\paragraph{Grid and chain.}
Fix
$N\in\mathbb N$, $h:=2L/N$, $E_h:=\{y_i=-L+ih:\ i\in\mathbb Z/N\mathbb Z\}$
and $b^\pm:=\max\{\pm b,0\}$. The chain starts from
$\rho_h:=(\mathrm p_h)_\#\rho$, where $\mathrm p_h$ maps each point to a nearest grid
point, with a fixed choice at ties. The costs $f$ and $g$ are restricted
to grid states and probability laws supported on $E_h$.
The Markov chain approximation of
\cite{KushnerDupuis2001} has rates
\begin{equation}\label{eq:rates}
    Q_t(y_i,y_{i+1},a)=\frac{\sigma^2}{2h^2}+\frac{b^+}{h},\qquad
    Q_t(y_i,y_{i-1},a)=\frac{\sigma^2}{2h^2}+\frac{b^-}{h}
    \qquad(\text{evaluated at }(t,y_i,a)),
\end{equation}
and no other transitions, so that
\begin{equation}\label{eq:qbar}
    \bar q\le\frac{\lVert\sigma^2\rVert_\infty}{h^2}+\frac{\lVert b\rVert_\infty}{h}.
\end{equation}
For $\phi\in\R^{E_h}$ let $D^+\phi_i:=(\phi_{i+1}-\phi_i)/h$, $D^-\phi_i:=(\phi_i-\phi_{i-1})/h$
(periodic indices), $D:=D^+$ and $\Delta_h:=D^+D^-=D^-D^+$. The generator of
the chain,
\begin{equation}\label{eq:gen-h}
    \LL^a_{h,t}\phi(y_i)=\sum_{y'\ne y_i}Q_t(y_i,y',a)\bigl(\phi(y')-\phi(y_i)\bigr)
    =b^+D^+\phi_i-b^-D^-\phi_i+\tfrac12\sigma^2\,\Delta_h\phi_i ,
\end{equation}
is the upwind discretization of $\LL^a_t$. With $E=E_h$ the game is an
instance of \Cref{ex:finite}; the action space is unchanged.

\paragraph{Discrete Sobolev metric.}
Let $(\phi,\phi')_h:=h\sum_i\phi_i\phi'_i$ and $\lVert \phi\rVert_h:=(\phi,\phi)_h^{1/2}$. Since
$(D\phi,\phi')_h=-(\phi,D^-\phi')_h$, the adjoint of $D$ is $-D^-$ and
$(D^j\phi,D^j\phi')_h=((-\Delta_h)^j\phi,\phi')_h$. For $m\ge1$ let
\begin{equation}\label{eq:Hm}
    (\phi,\phi')_H:=\sum_{j=0}^m\binom mj(D^j\phi,D^j\phi')_h=\bigl((I-\Delta_h)^m\phi,\phi'\bigr)_h,
    \qquad H^m_h:=(\R^{E_h},(\cdot,\cdot)_H).
\end{equation}
A measure on $E_h$ is a vector of masses $(\xi_i)$ and
$\pair{\phi}{\xi}=\sum_i\phi_i\xi_i=(\phi,\xi/h)_h$, so
\begin{equation}\label{eq:riesz-h}
    R_H\xi=(I-\Delta_h)^{-m}(\xi/h),\qquad
    \lVert \xi\rVert_{H^*}^2=\bigl(\xi/h,(I-\Delta_h)^{-m}\xi/h\bigr)_h .
\end{equation}
Since $-\Delta_h$ is circulant with eigenvalues $\frac4{h^2}\sin^2(\pi n/N)$,
$n=0,\ldots,N-1$, \eqref{eq:riesz-h} is one forward and one inverse FFT with
the multiplier $(1+\frac4{h^2}\sin^2\frac{\pi n}N)^{-m}$, the grid version of
$(1+|n|^2)^{-m}$ in \Cref{lem:riesz-explicit}(b). By \eqref{eq:time-exact} the
dual step of \Cref{alg:PDKLMP} is
\begin{equation}\label{eq:dual-step-h}
    u_{k+1}(t)=u_k(t)-\tau\,(I-\Delta_h)^{-m}\Bigl[\frac{d_T}{h}+\int_t^T\frac{d_s}{h}\,ds\Bigr],
\end{equation}
a backward cumulative sum of the integrated forward residuals followed by one FFT solve
per time slice, that is $O(N_tN\log N)$ operations per half-step for $N_t$
time nodes.

\begin{thm} \label{thm:mesh}
Let $m\ge3$ and $H=H^m_h$. Then, for every $N$:
\begin{enumerate}[label=(\roman*)]
\item (Standing assumptions.) (S2) and (S3) hold with $c_H=(1+\tfrac1{2L})^{1/2}$ and $\ell_H=c_H\bigl(\lVert b\rVert_\infty+\tfrac12\lVert\sigma^2\rVert_\infty\bigr)$, so $c_H$, $\ell_H$ and $L_A$ (\Cref{lem:elementary}) do not depend on $h$. Consequently \eqref{eq:step-size}, and with it \Cref{lem:PD-one-step,thm:PD-ergodic-minty,cor:primal-minty-residual}, hold with $\lambda,\tau$ independent of the mesh.
\item (Entropy.) If $A$ is finite and the reference masses satisfy the uniform lower bounds of \Cref{thm:finite-state}, then part (a) of that theorem holds. In particular, if $\nu_0(y_i)\ge c/N$ and $\pi_{0,t}(y_i,a)\ge c/(N|A|)$ for all $i,a$ and a.e.\ $t$, with $c>0$ independent of $N$, then $\bar D\le(1+T)\log(N|A|/c)$.
\item (Dual optimizers.) For each MFE $\pi^*$, the solution $V$ of \eqref{eq:hjb-finite} is an optimal dual variable. If $\bD(\pi^*\|\pi_0)<\infty$, \Cref{cor:primal-minty-residual}(c) applies under its stated assumptions. In the grid metric,
\[
    \tn{V}^2=\nH{V(T)}^2+\int_0^T\nH{\dot V(t)}^2\,dt.
\]

\item (Correction to an admissible flow.) The correction $\projr$ of \Cref{thm:finite-state}(c) satisfies \Cref{ass:projection} with
\begin{equation}
\label{eq:CPi-h}
    C_{\mathrm{proj}}=\sqrt{2L}\,\Bigl(1+\frac4{h^2}\Bigr)^{m/2}\bigl(1+(1+4\bar q^2T)(1+2\bar qT)^2\bigr)^{1/2},
\end{equation}
with $\bar q$ from \eqref{eq:qbar}; in particular $C_{\mathrm{proj}}=O(h^{-m-4})$ as $h\to0$.
\end{enumerate}
\end{thm}

The proof is in \Cref{subsec:proof-mesh}. Thus, for finite actions, under the reference bounds with a common $c$ in (ii) and \Cref{ass:LL,ass:Lip} for the chain game, and with $\lambda,\tau$ chosen once for all $h$, the ergodic Minty estimate against every feasible comparison point holds with error $(\bar D/\lambda+\tn{u'-u_0}^2/(2\tau))/K$ and $\bar D=O(\log(1/h))$, cluster points of $(\bar \pi_K)_K$ are equilibria of the chain game, and the corrected average satisfies $\Expl(\hat\pi_K)=O(K^{-1/2})$ in the chain game at each fixed $h$.
For interval actions, the grid limit below uses feasible comparisons with entropy $o(K_N)$, constructed in \Cref{lem:mesh-recovery}, in place of an entropy bound for every feasible measure.

\begin{proposition}
\label{prop:mesh-limit}
Consider the periodic diffusion and the continuous-time chains
\eqref{eq:rates}, with $h=2L/N$ and
$\rho_h=(\mathrm p_h)_\#\rho$. Assume that $A$ is either a finite set or a nondegenerate compact
interval, and let $\zeta$ be the uniform probability on $A$. Suppose
that $b$ is bounded, Borel in time and jointly continuous in the
state and action for each time, and that $\sigma=\sigma(t,y)$ is action-independent, bounded,
Borel in time and uniformly Lipschitz in the state, with
$\inf_{t,y}|\sigma(t,y)|>0$. Assume \Cref{ass:LL,ass:Lip}, and suppose
that the bounded costs $f(t,\cdot,\cdot,\cdot)$ and $g$ are jointly
continuous, for each $t$, when $\cP(E)$ carries its weak topology.

Fix an integer $m\ge3$, use the discrete Sobolev inner products
\eqref{eq:Hm}, and fix $\lambda,\tau>0$ satisfying
\eqref{eq:step-size} with the mesh-independent constants of
\Cref{thm:mesh}(i). For each $N\ge3$, run the exact continuous-time
iteration \eqref{eq:PD-KL-MP}, starting from $u_0^N=0$ and reference
measures satisfying
\[
    \nu_0^N(y_i)\ge\frac cN,
    \qquad
    \pi_{0,t}^N(\{y_i\},da)\ge\frac cN\,\zeta(da)
    \quad\text{for every }i\text{ and a.e. }t,
\]
where $c>0$ is independent of $N$. Write $\bar\pi_{N,K}$ for the
predictor average after $K$ iterations.

If $K_N$ satisfies $\log N/K_N\to0$, every limit point of
$\bar\pi_{N,K_N}$, after embedding grid measures in the torus, is a
mean field equilibrium of the diffusion game. The topology is weak
convergence of the terminal law and of the running occupation measure
$dt\,\bar\pi_{N,K_N,t}(dy,da)$. These pairs are relatively compact;
if the equilibrium occupation pair is unique, the entire sequence
converges to it.

Along every subsequence converging to an equilibrium occupation pair
$\pi^*$, the population flows also satisfy
\begin{equation}\label{eq:mesh-population-convergence}
 \lVert\cM\bar\pi_{N,K_N}-\cM\pi^*\rVert_{W_1,2}\longrightarrow0.
\end{equation}
If the equilibrium population flow is unique, the entire population
sequence converges to it in this metric.
\end{proposition}

The proof is in \Cref{subsec:proof-mesh-limit}. The quantitative
strengthening below uses smoothness and uniform convexity of the
control cost.

\begin{assumption}\label{ass:convex-grid}
Assume (C1)--(C3) of \Cref{ass:convex-diffusion}. In addition:
\begin{enumerate}[label=(G\arabic*)]
\item The equilibrium value function satisfies
$u^*\in H^1(0,T;H^{m+1}(\mathbb T))$.
\item There is $L_W<\infty$ such that, for every $t,y,\mu,\mu'$,
\[
 |f_0(t,y,\mu)-f_0(t,y,\mu')|
       +|g(y,\mu)-g(y,\mu')|
 \le L_W W_1(\mu,\mu').
\]
\item On the grid $E_h$, $h=1/N$, the references satisfy the lower bounds in \Cref{prop:mesh-limit}.
\end{enumerate}
\end{assumption}

\begin{thm}\label{thm:convex-grid}
Suppose \Cref{ass:convex-grid} holds. On $E_h$, $h=1/N$, $N\ge3$,
use the continuous-time rates
\[
 Q_h^a(y_i,y_{i+1})=\frac{\sigma^2}{2h^2}+\frac{a^+}{h},
 \qquad
 Q_h^a(y_i,y_{i-1})=\frac{\sigma^2}{2h^2}+\frac{a^-}{h},
\]
initial law $\rho_h=(\mathrm p_h)_\#\rho$, and the costs restricted
to grid states. Use $H_h^m$ from \eqref{eq:Hm}, start from
$u_0^N=0$, and fix admissible steps $\lambda,\tau>0$ independently
of $N$. Let $\bar\pi_{N,K}$ be the predictor average of the exact
continuous-time iteration after $K$ steps, embedded in
$\mathbb T\times A$. Put
\[
 E_{N,K}:=\int_0^T\!\!\int|a-a^*(t,y)|^2
                         \,\bar\pi_{N,K,t}(dy,da)\,dt,
 \quad
 \eta_{N,K}:=\|\cA_h\bar\pi_{N,K}-\rho_h\|_{Y_h^*},
\]
where $Y_h=H^1(0,T;H_h^m)$ has the inner product specified for the
dual variables. There is $C$ independent of $N$ and $K$ such that
\begin{equation}\label{eq:convex-grid-energy}
 E_{N,K}\le C\left(\frac{\log(eNK)}K+h^{1/3}\right),\qquad
 \eta_{N,K}\le C\left(\frac{\log(eNK)}{K^2}
                         +\frac{h^{1/3}}K\right)^{1/2},
\end{equation}
and
\begin{equation}\label{eq:convex-grid-flow-rate}
 \lVert\cM\bar\pi_{N,K}-\bm\mu^*\rVert_{W_1,2}
 \le C\left[
      \left(\frac{\log(eNK)}K+h^{1/3}\right)^{1/2}
      +\left(\frac{\log(eNK)}{K^2}
                  +\frac{h^{1/3}}K\right)^{1/[2(m+1)]}
      \right].
\end{equation}
The same bound holds for the state--action distance
\eqref{eq:convex-occupation-distance}, with $\bar\pi_K$ replaced
by $\bar\pi_{N,K}$.
\end{thm}

\begin{corollary}\label{cor:convex-grid-balanced}
Under the assumptions of \Cref{thm:convex-grid}, if
$K\ge h^{-m/3}$, both distances in that theorem are at most
$Ch^{1/6}$. In particular, for $m=3$, taking $K\ge N$ gives
an $O(N^{-1/6})$ bound.
\end{corollary}

The proof is in \Cref{subsec:proof-convex-grid}. It adapts the
comparison of control costs used by Bertucci and Cecchin
\cite[Section~3.2]{BertucciCecchin2024}, using a chain occupation measure
that satisfies the forward equation, \eqref{eq:master} and a backward linear equation. The comparison has spatial error
$O(h^{1/3})$; convexity converts this into the $O(h^{1/6})$ term in
the population estimate. The theorem concerns spatial discretization
and proximal iterations, with physical time kept continuous.

\begin{rmk}
\label{rmk:mesh}
The operator bounds and admissible step sizes are uniform in $h$.
Under the reference conditions of \Cref{thm:mesh}(ii), the entropy
term in the ergodic Minty bound grows only like $\log(1/h)$; the
term involving the comparison dual function is uniform only when the corresponding norms
$\tn{\cdot}$ are uniformly bounded. The conversion of the weak
forward-equation residual into a total-variation correction error is not
uniform: \eqref{eq:CPi-h} is an inverse-inequality bound and is not
claimed to be sharp. A global total-variation correction estimate
cannot in general have a mesh-independent constant in this Sobolev
residual norm: increasingly oscillatory terminal discrepancies can
have fixed total variation and a vanishing negative-Sobolev norm.
\Cref{prop:mesh-limit} approximates diffusion occupation measures
by feasible chain occupation measures. The metric
estimates in \Cref{thm:convex-diffusion,thm:convex-grid} use control
convexity and a backward parabolic estimate.

In \Cref{sec:experiments}, the Stampacchia gap of the corrected average
is computed directly, so the correction constant is not used
numerically. For fixed $h$ the limiting problem is a continuous-time
finite-state MFG in the sense of \cite{GomesMohrSouza2013}. Under the
additional continuity and nondegeneracy hypotheses of
\Cref{prop:mesh-limit}, every limit obtained as the mesh tends to zero
with $\log N/K_N\to0$ is a diffusion MFE. This combines Markov-chain approximation
\cite{KushnerDupuis2001} with the Minty estimate for the game.
\end{rmk}

\begin{rmk}
\label{rmk:gomes-osher}
With regard to the recent paper by Gevorgyan and Gomes \cite{GevorgyanGomes2026},
\begin{enumerate}
\item In \cite{GevorgyanGomes2026}, the operator
is formed from derivatives of the Lagrangian, with signs corresponding
to its minimization and maximization variables,
\[
\mathcal J(m,u,\varpi) = \int\!\!\int(H(x,\varpi+u_x)-u_t)\,m\,dx\,dt-\int u(x,0)\bar m\,dx+\int Q\varpi\,dt,
\]
which is affine in the density $m$ and convex in $(u,\varpi)$ because $H$
is convex in $p$; the mean-field coupling enters only through the scalar price $\varpi$, linearly. Monotonicity is that of the convex Hamiltonian, and the nonlinearity $H$ stays in the operator, which is why Lipschitz
bounds need slope bounds and a bounded admissible set. Here the Hamiltonian
is linearized by the occupation-measure relaxation (the player's problem is
a linear program over $(\nu,(\pi_t)_t)$), $\cF$ is affine in $u$ and in $\pi$
except through $\cC_{\cM \pi}$, and monotonicity is exactly Lasry--Lions
monotonicity of a general nonlinear coupling $f(t,y,a,\mu_t)$, $g(y,\mu_T)$
(\Cref{prop:monotone}), with no convexity of the Hamiltonian, no differentiability and no uniqueness of best responses required.
\item \cite{GevorgyanGomes2026} and
\cite{MengLiuLiOsher2024} update the density in $L^2$ by a Euclidean step followed by clipping to a box $[0,M_m]$ or by $(\cdot)_+$; the density is not a probability measure along the iteration, and mass conservation is
recovered only at the solution. Here the primal iterates are probability measures at every step, because the proximal update with relative entropy is an
exponential tilt (\Cref{def:prox}); Pinsker's inequality \eqref{eq:pinsker-dyn} controls total variation
by relative entropy and gives the one-step estimate without an
a priori upper bound on the density.
\item In \cite{GevorgyanGomes2026} the value
function is measured in the space--time $H^1$ norm, and $(L_{tt}+L_{xx})^{-1}$ is applied by solving a linear space--time system, at a cost of $O(N^2N_t)$;
\cite{MengLiuLiOsher2024} apply $(I-\partial_{tt}-\Delta)^{-1}$, the Riesz
map of the same space--time $H^1$ norm, to the continuity residual by a space--time FFT, and observe that in the viscous case their step size scales with the mesh because the $H^1$ metric does not control second derivatives. Here the same linear step is taken using the chosen norm on $H^1(0,T;H^m_h)$: the time part is integrated exactly in \eqref{eq:dual-step-h}, the space part is $(I-\Delta_h)^{-m}$, and $m\ge3$ is what makes the second-order generator bounded on $H$ uniformly in $h$, so no admissible set, no slope bound and no projection are needed. No a priori bound on the primal densities is needed, since the KL tilts keep the iterates in $\Pi$.
\end{enumerate}
\end{rmk}

\section{Galerkin implementation}\label{subsec:inexact}

We retain the diffusion generator and approximate each Riesz correction
in a Galerkin space of finite dimension. The primal measures are represented
by densities relative to a reference measure, and the Mirror-Prox
updates are those of \Cref{alg:PDKLMP}. Theorem~\ref{thm:inexact} quantifies
how the Galerkin energy errors affect convergence. The same analysis
extends to Ritz minimization over a parametrized class of functions.

\paragraph{Representation of the primal measures.}
The primal variable remains a pair of probability measures. Since exponential tilting preserves absolute continuity, the running iterate is represented by its density $d\pi_{k,t}/d\pi_{0,t}\propto e^{-\lambda\Sigma_k(t,\cdot,\cdot)}$ from \Cref{rmk:no-discretization}. The dual function is expanded in a basis, and the primal measures are integrated against the reference measure. Quadrature or samples from that measure give the normalization constants and the moments entering the forward-equation residual. The terminal component is treated in the same way relative to $\nu_0$.

A quadrature rule with nodes $(y_j,a_j)$ replaces $\pi_{0,t}$ by $\sum_jq_{0,t,j}\delta_{(y_j,a_j)}$, with positive masses summing to one, and $\pi_{k,t}$ by $\sum_jq_{k,t,j}\delta_{(y_j,a_j)}$ with $q_{k,t,j}\propto q_{0,t,j}e^{-\lambda\Sigma_k(t,y_j,a_j)}$. A primal half-step with running residual $\Gamma$ updates the masses by
\begin{equation}\label{eq:quadrature-measure-update}
    q^+_{t,j}
      =\frac{q_{k,t,j}\exp[-\lambda\Gamma(t,y_j,a_j)]}
             {\sum_iq_{k,t,i}\exp[-\lambda\Gamma(t,y_i,a_i)]},
    \qquad
    \int\varphi\,d\pi^+_t\approx\sum_jq^+_{t,j}\varphi(y_j,a_j).
\end{equation}
Applying $\LL^{a_j}_t$ to the spatial test functions supplies the generator moments directly. {The relative entropy of the atomic measures, $\sum_jq_{k,t,j}\log(q_{k,t,j}/q_{0,t,j})$, converges to $\KL(\pi_{k,t}\|\pi_{0,t})$ if the quadrature converges for the unnormalized tilt and for that tilt multiplied by its logarithm.}

\paragraph{The auxiliary minimization problem.}
For $r\in H^1(0,T;H)^*$, in particular $r=\cA \pi-\rho$, define
\begin{equation}\label{eq:energy}
    \cQ_r(w):=\tfrac12\tn w^2-\pair{r}{w},\qquad w\in H^1(0,T;H) .
\end{equation}
The unique minimizer is the required correction $Rr$: its
first-order condition is $\pp{w}{v}=\pair rv$ for every
$v\in H^1(0,T;H)$, precisely the defining Riesz equation.
Choosing a finite linear basis and minimizing $\cQ_r$ over its span gives a
Gram-matrix linear system, the \emph{Galerkin} approximation. Minimizing
the same energy over a chosen class of functions is called \emph{Ritz
minimization}; the class may also be nonlinear, for example a neural
network. Expanding the square gives
\begin{equation}\label{eq:energy-identity}
    \cQ_r(w)-\cQ_r(Rr)=\tfrac12\tn{w-Rr}^2,\qquad \cQ_r(Rr)=-\tfrac12\tns r^2 .
\end{equation}
Thus the energy excess measures the squared error in the correction.

Here $H$ is the full test space. The bounds below concern the dual correction; errors in numerically evaluating the primal tilts require separate control.

\paragraph{Effect of errors in the dual correction.}
Let $(\epsilon_k)_{k\ge0}$ and $(\epsilon_{k+1/2})_{k\ge0}$ be nonnegative.
The \emph{inexact} iteration is \Cref{alg:PDKLMP} with the two dual
corrections replaced by elements $\tilde w_k,\tilde w_{k+1/2}\in H^1(0,T;H)$ such that
\begin{equation}\label{eq:inexact}
    \tn{\tilde w_k-R(\cA \pi_k-\rho)}\le\epsilon_k,\qquad
    \tn{\tilde w_{k+1/2}-R(\cA \pi_{k+1/2}-\rho)}\le\epsilon_{k+1/2},
\end{equation}
the primal tilts being unchanged. With
$e_k:=R^{-1}\tilde w_k-(\cA \pi_k-\rho)\in H^1(0,T;H)^*$ and likewise $e_{k+1/2}$,
the Riesz isometry gives $\tns{e_k}\le\epsilon_k$ and
$\tns{e_{k+1/2}}\le\epsilon_{k+1/2}$, and the inexact iteration is
\begin{equation}\label{eq:perturbed}
    z_{k+1/2}=\prox_{z_k}\bigl(\cF(z_k)+(0,e_k)\bigr),\qquad
    z_{k+1}=\prox_{z_k}\bigl(\cF(z_{k+1/2})+(0,e_{k+1/2})\bigr),
\end{equation}
the exact proximal iteration of \Cref{def:prox} for operator values
perturbed in the dual component only. Here $\Phi$,
$\theta$ and $\bar z_K$ are as in \Cref{sec:LL}.

\begin{thm}\label{thm:inexact}
Suppose \Cref{ass:standing,ass:LL,ass:Lip} hold and $\theta<1$.
Let $z^*=(\pi^*,u^*)$ solve the PD-VI with
$\bD(\pi^*\|\pi_0)<\infty$.
Run \eqref{eq:perturbed} with exact primal tilts. Suppose the approximate
dual corrections satisfy the energy bounds
\begin{equation}\label{eq:inexact-energy-tolerances}
\begin{aligned}
 2\cQ_{\cA\pi_k-\rho}(\tilde w_k)
       +\tns{\cA\pi_k-\rho}^2&\le\epsilon_k^2,\\
 2\cQ_{\cA\pi_{k+1/2}-\rho}(\tilde w_{k+1/2})
       +\tns{\cA\pi_{k+1/2}-\rho}^2&\le\epsilon_{k+1/2}^2,
\end{aligned}
\end{equation}
with nonnegative tolerances such that
\[
 \sum_{k\ge0}(\epsilon_k+\epsilon_{k+1/2})<\infty.
\]
There is a constant $C$, independent of $K$ and of the comparison
point, such that for every $K\ge1$,
\begin{equation}\label{eq:inexact-forward-rate}
 \tns{\cA\bar\pi_K-\rho}\le\frac{C}{K},
\end{equation}
and, for every $z'=(\pi',u')\in\cZ$ with
$\bD(\pi'\|\pi_0)<\infty$,
\begin{equation}\label{eq:inexact-main-comparison}
 \pair{\cF(z')}{\bar z_K-z'}
 \le\frac{\Phi(z',z_0)+C(1+\tn{u'-u_0})}{K}.
\end{equation}
For a feasible comparison $\pi'\in\Pizero$, the primal estimate is
\[
 \pair{\cC_{\cM\pi'}}{\bar\pi_K-\pi'}
 \le\frac1K\left[
       \frac{\bD(\pi'\|\pi_0)}{\lambda}
       +\frac{\tn{u_0}^2}{2\tau}
       +C(1+\tn{u_0})\right].
\]
The constant depends on the fixed step sizes, the model bounds,
$\Phi(z^*,z_0)$, $\tn{u^*-u_0}$ and the summable error sequence.
\end{thm}

For example, tolerances
$\epsilon_k,\epsilon_{k+1/2}\le c(k+1)^{-1-\gamma}$ with
$\gamma>0$ satisfy the theorem. The required energy excess at each
half-step is then at most $c^2(k+1)^{-2-2\gamma}/2$.

The proof of \Cref{thm:inexact} uses \eqref{eq:energy-identity} to turn
\eqref{eq:inexact-energy-tolerances} into correction-error bounds, then
applies the stability and telescoping estimates in the next lemma.

\begin{lemma}
\label{lem:inexact}
Let \Cref{ass:standing,ass:Lip} hold, let $\theta<1$, and let $(z_k,z_{k+1/2})$ be generated by \eqref{eq:perturbed}.
\begin{enumerate}
\item[(i)] For every $z'=(\pi',u')\in\cZ$ with $\bD(\pi'\|\pi_0)<\infty$ and
every $k\ge0$,
\begin{equation}\label{eq:one-step-inexact}
\begin{aligned}
    \pair{\cF(z_{k+1/2})}{z_{k+1/2}-z'}
    &\le\Phi(z',z_k)-\Phi(z',z_{k+1})
      -(1-\theta)\Phi(z_{k+1/2},z_k)\\
    &\quad+\frac{\tau(\epsilon_k+\epsilon_{k+1/2})^2}{2(1-\theta)}
      +\epsilon_{k+1/2}\tn{u_{k+1/2}-u'} .
\end{aligned}
\end{equation}
\item[(ii)] If in addition \Cref{ass:LL} holds, then for every such $z'$
and every $K\ge1$, with $B_K(u'):=\max_{0\le k<K}\tn{u_{k+1/2}-u'}$,
\begin{equation}\label{eq:ergodic-inexact}
\begin{aligned}
    \pair{\cF(z')}{\bar z_K-z'}
    &\le\frac1K\Biggl[\Phi(z',z_0)
      +\frac{\tau}{2(1-\theta)}
       \sum_{k<K}(\epsilon_k+\epsilon_{k+1/2})^2\\
    &\hspace{4.5cm}
      +B_K(u')\sum_{k<K}\epsilon_{k+1/2}\Biggr].
\end{aligned}
\end{equation}
\item[(iii)] If in addition \Cref{ass:LL} holds and
$z^*=(\pi^*,u^*)$ solves the PD-VI with
$\bD(\pi^*\|\pi_0)<\infty$ and $\bar\epsilon:=\sum_{k\ge0}(\epsilon_k+\epsilon_{k+1/2})<\infty$,
then, with $M_\tau:=2\tau L_A^2(1+T)$ and
$C_\tau:=\sqrt{2\tau}+\sqrt\tau+\sqrt{\tau/(2(1-\theta))}$,
\begin{align*}
    \sqrt{\Phi(z^*,z_k)+M_\tau}
    &\le\sqrt{\Phi(z^*,z_0)+M_\tau}+C_\tau\bar\epsilon,
      &&k\ge0,\\
    \sup_k\tn{u_{k+1/2}-u^*}
    &\le\sqrt{2\tau}\Bigl(\sqrt{\Phi(z^*,z_0)+M_\tau}
      +C_\tau\bar\epsilon\Bigr)
      +\tau\Bigl(2L_A\sqrt{1+T}+\sup_k\epsilon_k\Bigr).
\end{align*}
\end{enumerate}
\end{lemma}

Under the hypotheses of part (iii), $B_K(u')$ is bounded uniformly in
$K$ for each fixed $u'$. Equation \eqref{eq:ergodic-inexact} then gives
an $O(K^{-1})$ comparison estimate, with a constant depending on
$\bar\epsilon$ and $\sum_k(\epsilon_k+\epsilon_{k+1/2})^2$. One can also take $u'=u_K$
in \eqref{eq:ergodic-inexact}, for which
$B_K(u_K)=\max_{k<K}\tn{u_{k+1/2}-u_K}$ is available along the run.

\begin{lemma}
\label{lem:ritz-error}
Let $r=\cA \pi-\rho$. Let $V_h\subset H^1(0,T;H)$ be a closed subspace and $w_h\in V_h$ the unique element with $\pp{w_h}{v}=\pair rv$ for all $v\in V_h$. Then $w_h=\operatorname{proj}_{V_h}Rr$ is the orthogonal projection, for $\pp{\cdot}{\cdot}$, of the exact correction, $w_h$
minimizes $\cQ_r$ over $V_h$, and
\[
    \tn{w_h-Rr}^2=\tns r^2-\tn{w_h}^2=\min_{v\in V_h}\tn{v-Rr}^2 .
\]
In a basis $\{v_i\}$ of $V_h$ with Gram matrix $\mathbf G$ and load vector $b_i=\pair r{v_i}$, this reads $\tn{w_h-Rr}^2=\tns r^2-\mathbf b^\top\mathbf G^{-1}\mathbf b$.
\end{lemma}

\begin{rmk}
\label{rmk:space-only}
Formula \eqref{eq:time-exact} reduces the Riesz correction to spatial energy minimizations. Indeed, at each time the exact correction $w(t)$ minimizes the spatial energy $\phi\mapsto\frac12\nH\phi^2-\pair{d_T+\int_t^Td_s\,ds}{\phi}$, and the space--time error of a candidate $\tilde w\in H^1(0,T;H)$ decomposes as $\tn{\tilde w-Rr}^2 =\nH{\tilde w(T)-R_Hd_T}^2 +\int_0^T\nH{\partial_t\tilde w(t)+R_Hd_t}^2\,dt$. The spatial solves for $w(T)$ and $-\partial_t w(t)$ have data $d_T$ and $d_t$, respectively. The displayed identity controls the full correction error through its terminal value and time derivative, and can be used to choose the tolerances in
\eqref{eq:inexact}.
\end{rmk}

\begin{rmk}
\label{rmk:deep-ritz}
For $H=H^m(\mathbb T^d)$ with $m>d/2+2$, the spatial Riesz problem is
$(I-\Delta)^m\phi=\ell$, and its energy involves derivatives up to order
$m$. When $w_\vartheta$ is represented by a neural network, these
spatial derivatives are evaluated by automatic differentiation and the
energy is minimized over the network parameters. The data $d_t$ are
evaluated against test functions by Monte Carlo on the reference paths
of \Cref{rmk:no-discretization}.

Identity \eqref{eq:energy-identity} expresses the correction error in the $H^1(0,T;H)$ norm through the energy, and \Cref{lem:inexact} bounds its effect on the iteration. For a reproducing kernel Hilbert space satisfying (S2)--(S3), \eqref{eq:rkhs-integrated-residual} gives a direct alternative: average $\mathsf k_x(y)$ against the state laws and $(\LL_t^a\mathsf k_x)(y)$ against the running state--action laws. Whether one estimates the energy or computes the correction directly from samples, the tolerances in \eqref{eq:inexact} must bound the error in the full $H^1(0,T;H)$ norm.
\end{rmk}

\section{Numerical experiments}\label{sec:experiments}

The experiments implement \eqref{eq:PD-KL-MP} in finite coordinates, with the dual step in the form \eqref{eq:dual-step-coords}: the dual variable, its time derivative and its terminal value are propagated separately, so
the HJB residuals \eqref{eq:slacks} use $\partial_t u_k$ without finite
differences in time. Time integrals use midpoint quadrature. For the diffusion model, \Cref{subsec:exp-diffusion} compares two spatial Riesz maps
on a family of periodic Markov-chain grids, while \Cref{subsec:exp-fourier}
uses a Fourier basis and the diffusion generator directly. Averages use the
predictor iterates $z_{k+1/2}$ as in
\eqref{eq:ergodic-avg}.\footnote{Code and reproduction commands are listed in \url{https://github.com/Lethargy/dynamic-kl-mirror-prox-mfe}.}

\subsection{A finite-state congestion game}\label{subsec:exp-finite}

Let $E=\{1,\ldots,5\}$, $A=\{\text{left},\text{stay},\text{right}\}$ and
$T=1$. Every action has a background rate $0.1$ to each existing
neighbour; left and right add a rate $0.5$ to the left and right
neighbour. Thus $\bar q=0.7$, $\ell_H=1.4$ and $L_A=3.137$. The costs are
\[
    f(t,y,a,\mu)=2\,\mu(y)+0.1\cdot\mathbf 1\{a\ne\text{stay}\},
    \qquad
    g(y,\mu)=\tfrac14(y-3)^2,
    \qquad
    \rho=\tfrac12\delta_1+\tfrac12\delta_5 .
\]
The congestion term is Lasry--Lions monotone (\Cref{rmk:LL}), and
\[
    \sup_y|f(\cdot,\mu)-f(\cdot,\mu')|
    \le2\max_y|\mu(y)-\mu'(y)|
    \le\lVert\mu-\mu'\rVert_1,
\]
so \Cref{ass:Lip} holds with $L_C=1$. The terminal cost pulls the
population to the centre and the congestion cost spreads it. With
$\pi_0$ uniform ($\underline\nu=1/5$, $\underline\pi=1/15$) \Cref{thm:finite-state} gives
$\bar D=4.32$ and $C_{\mathrm{proj}}\approx9.50$; moreover $C_0=2.33$. We take $\lambda=\tau=0.99/(L_C+L_A)=0.239$, so that $\theta=0.99$ in \eqref{eq:step-size}, the midpoint rule on $100$ intervals, and $K\le2\times10^4$. The correction $\projr(\bar \pi_K)$ of \Cref{thm:finite-state}(c) is computed with exact matrix exponentials of
the piecewise-constant feedback, and the Stampacchia gap $\Expl(\hat \pi_K)$ with the value of the HJB system \eqref{eq:hjb-finite} for the frozen flow $\cM\hat \pi_K$, integrated by a fourth-order Runge--Kutta scheme.

\begin{figure}[t]
\centering
\includegraphics[width=\textwidth]{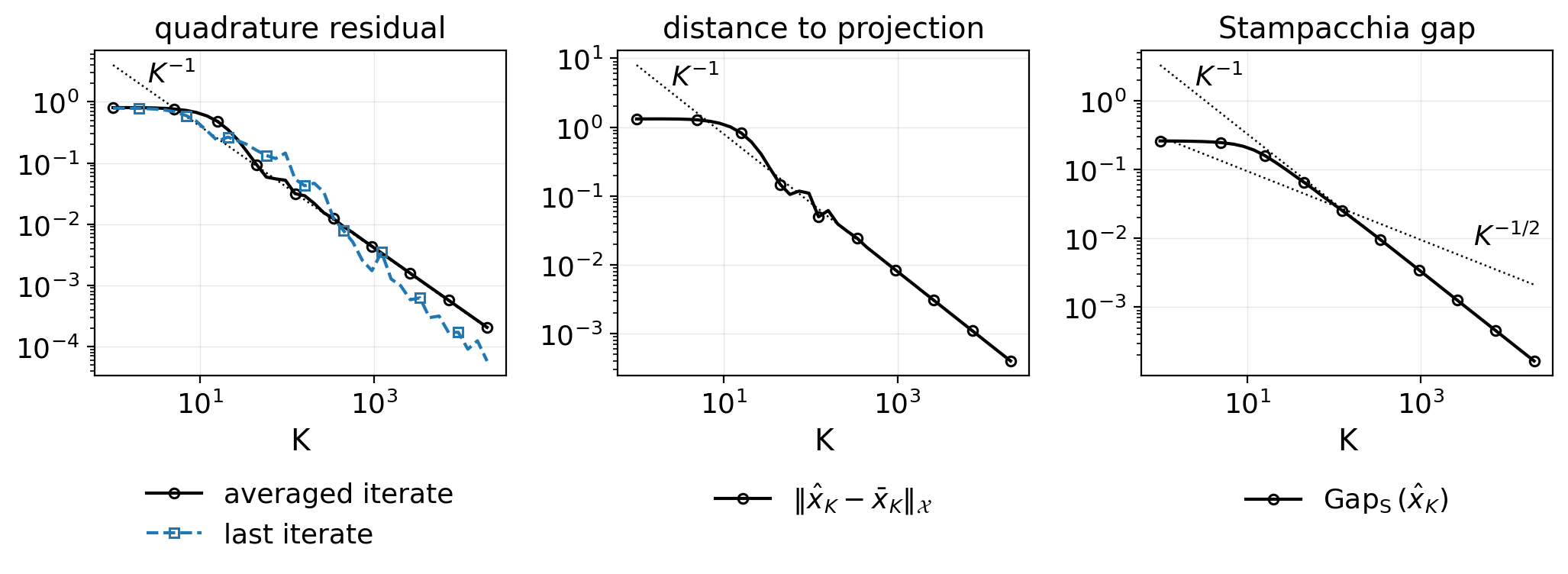}
\caption{Finite-state congestion game. Left: forward-equation residual $\tns{\cA\bar\pi_K-\rho}$ of the averaged iterate and of the last iterate. Middle: size $\nsharp{\hat \pi_K-\bar \pi_K}$ of the correction. Right: Stampacchia gap of the corrected average $\hat \pi_K$. Dotted and dash-dotted lines are the reference slopes $K^{-1}$ and $K^{-1/2}$.}
\label{fig:finite-rates}
\end{figure}

\begin{figure}[t]
\centering
\includegraphics[width=\textwidth]{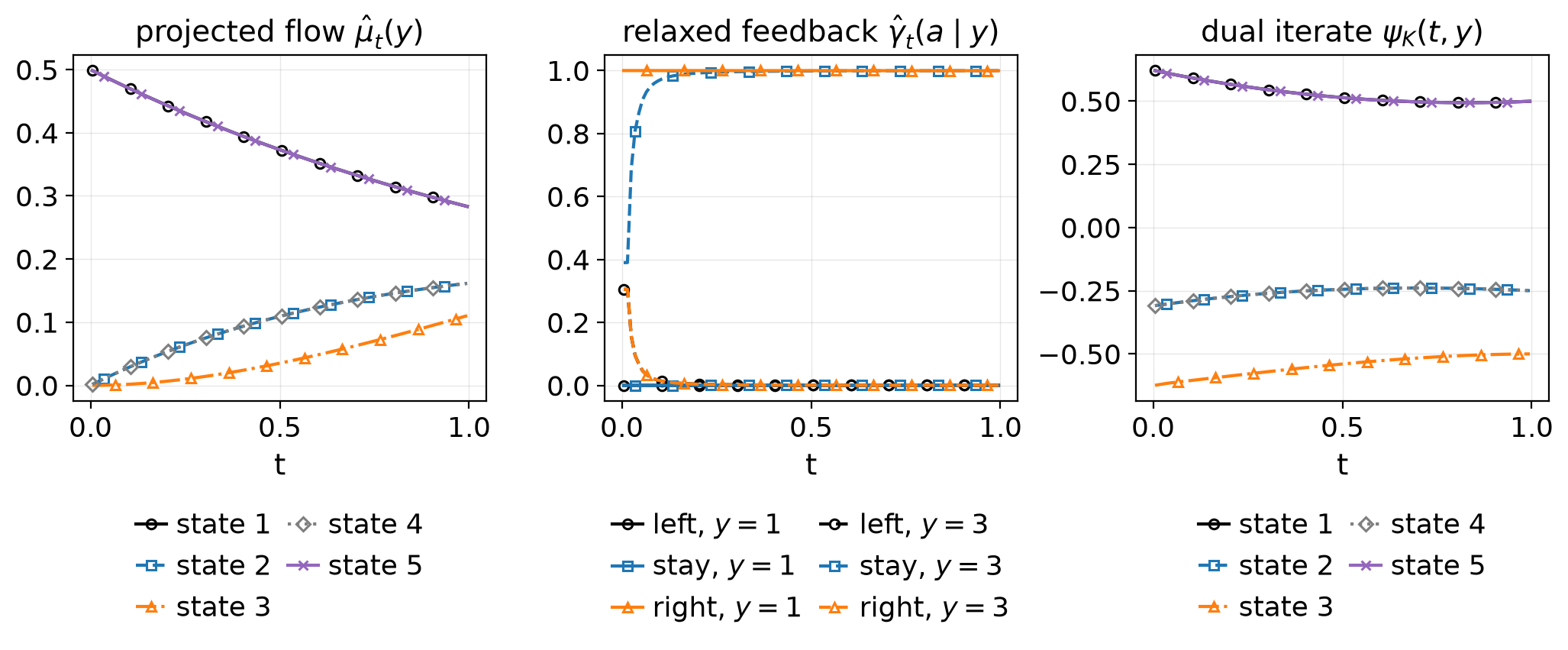}
\caption{Finite-state congestion game at $K=2\times10^4$. Left: the
corrected flow $\hat\mu_t$. Middle: the relaxed feedback extracted from
$\bar \pi_K$ at states $1$ and $3$. Right: the dual iterate $u_K(t,y)$.}
\label{fig:finite-eq}
\end{figure}

\Cref{fig:finite-rates} shows the three quantities controlled by \Cref{thm:quantitative-eps-MFE}. After a transient of about $30$ iterations, all three decay like $K^{-1}$: at $K=2\times10^4$ the residual is $2.1\times10^{-4}$, the size of the correction $4.0\times10^{-4}$ and the Stampacchia gap $1.65\times10^{-4}$, and the fitted constants are $\tns{\cA\bar\pi_K-\rho}\approx4/K$ and $\Expl(\hat \pi_K)\approx3.3/K$. The residual rate is the one predicted by \Cref{cor:primal-minty-residual}(c) and \Cref{thm:finite-state}(b), since an optimal dual variable exists; the norm $\tn{u_K-u_0}$ stabilizes at $0.985$, consistent with the boundedness of \Cref{rmk:dual-drift}(i). The Stampacchia gap decays at the rate $K^{-1}$ of the Minty gap rather than at the worst-case rate $K^{-1/2}$ of \Cref{prop:eps-relations}(ii). The explicit bound of \Cref{thm:quantitative-eps-MFE} for $u_0=0$ is $\varepsilon_K\le(2\bar D/\lambda+5C_0^2C_{\mathrm{proj}}^2/(2\tau))/K\approx5.14\times10^3/K$; the large constant is dominated by the correction term $5C_0^2C_{\mathrm{proj}}^2/(2\tau)$, with $C_{\mathrm{proj}}$ from \Cref{thm:finite-state}(c).
\Cref{fig:finite-eq} shows the computed equilibrium: the population starting at the two ends moves toward the
centre, the feedback at state $1$ is ``right'' throughout, and the feedback at state $3$ switches from a mixed action to ``stay'' within $t\le0.2$, once the congestion at the centre is established. Dual multipliers are
defined only up to spatially constant, time-dependent shifts in the PD-VI; \Cref{prop:PDVI-MFE}(iii) gives the normalization to an optimal HJB subsolution at a solution.

\subsection{A periodic controlled diffusion}\label{subsec:exp-diffusion}

Let $T=1$ and $E=\mathbb T:=\R/(2\pi\mathbb Z)$. Each player starts
from $X_0=0$, chooses $a_t\in\{-1,+1\}$ and follows
\[
    dX_t=a_t\,dt+\sqrt2\,dW_t\pmod{2\pi},
    \qquad \LL^a\phi=a\phi'+\phi''.
\]
The running and terminal costs are smooth periodic functions,
\begin{equation}\label{eq:periodic-example-costs}
    m_\mu(t):=\int_{\mathbb T}\sin(y)\,\mu_t(dy),\qquad
    f(t,y,a,\mu)=\sin(y)m_\mu(t),\qquad
    g(y,\mu)=-\varkappa\sin(y),\quad \varkappa\in\{0,1\}.
\end{equation}
For two occupation measures the Lasry--Lions pairing is
\[
    \pair{\cC_{\cM\pi}-\cC_{\cM\pi'}}{\pi-\pi'}
    =\int_0^1\bigl(m_{\cM\pi}(t)-m_{\cM\pi'}(t)\bigr)^2dt\ge0,
\]
and $|\sin y|\le1$ gives $L_C=1$. This is a compact-state diffusion
of \Cref{ex:diffusion}. For $\varkappa=0$, the balanced feedback
$\alpha(\pm1\mid t,y)=1/2$ generates the periodic heat flow from
$\delta_0$, which is symmetric and has $m_\mu\equiv0$. It is an
exact equilibrium, since every admissible deviation then has zero
cost. For $\varkappa=1$, the terminal reward favors states with positive
sine, while the running coupling penalizes alignment with the population.

\paragraph{Periodic grids and reference.}
We use $y_j=2\pi j/N_\pi$, $j\in\mathbb Z/N_\pi\mathbb Z$, with
$h=2\pi/N_\pi$ and $\rho_h=\delta_0$. All indices are periodic. The
Kushner--Dupuis rates are
\[
    Q_h^a(y_j,y_{j+1})=h^{-2}+a_+/h,\qquad
    Q_h^a(y_j,y_{j-1})=h^{-2}+a_-/h.
\]
Each grid uses all functions on $E_h$ as spatial tests. The terminal and running references are
\[
    \nu_0(y_j)=\frac1{N_\pi},\qquad
    \pi_{0,t}(y_j,a)=\frac1{2N_\pi},
\]
uniformly in time. The references are strictly positive and
\Cref{thm:finite-state,thm:mesh} apply to each continuous-time grid game;
in particular $\bar D=\log N_\pi+\log(2N_\pi)$.

\paragraph{Dual metrics and numerical integration.}
We compare the periodic Sobolev metric $H_h^3$, with Gram matrix
$h(I-\Delta_h)^3$, and the Euclidean metric on $\R^{E_h}$. Here
$\Delta_h$ is the periodic grid Laplacian. The Sobolev Riesz map is
\[
    R_H\xi=\operatorname{IFFT}\left[
      \frac{\operatorname{FFT}(\xi)_n}
      {h\left(1+4h^{-2}\sin^2(\pi n/N_\pi)\right)^3}\right].
\]
We use the three resolutions $(N_\pi,N_t)=(64,25),(128,50),(256,100)$, where $N_t$ is the number of midpoint intervals. We run $K=20{,}000$ iterations in each case. The Sobolev runs use $\lambda=\tau=0.25$ on every mesh; the Euclidean runs use their own admissible steps $\lambda=\tau=0.99/(L_C+L_A)$, with $L_A$ evaluated for each grid. To evaluate the Stampacchia gap, we extract the feedback from the predictor average and solve the forward equation for the periodic chain under that feedback. We solve the HJB system for the resulting population flow using an implicit Euler scheme whose transition operators are adjoint to those in the forward equation. The reported Stampacchia gap is the resulting cost difference.

\begin{figure}[t]
\centering
\includegraphics[width=\textwidth]{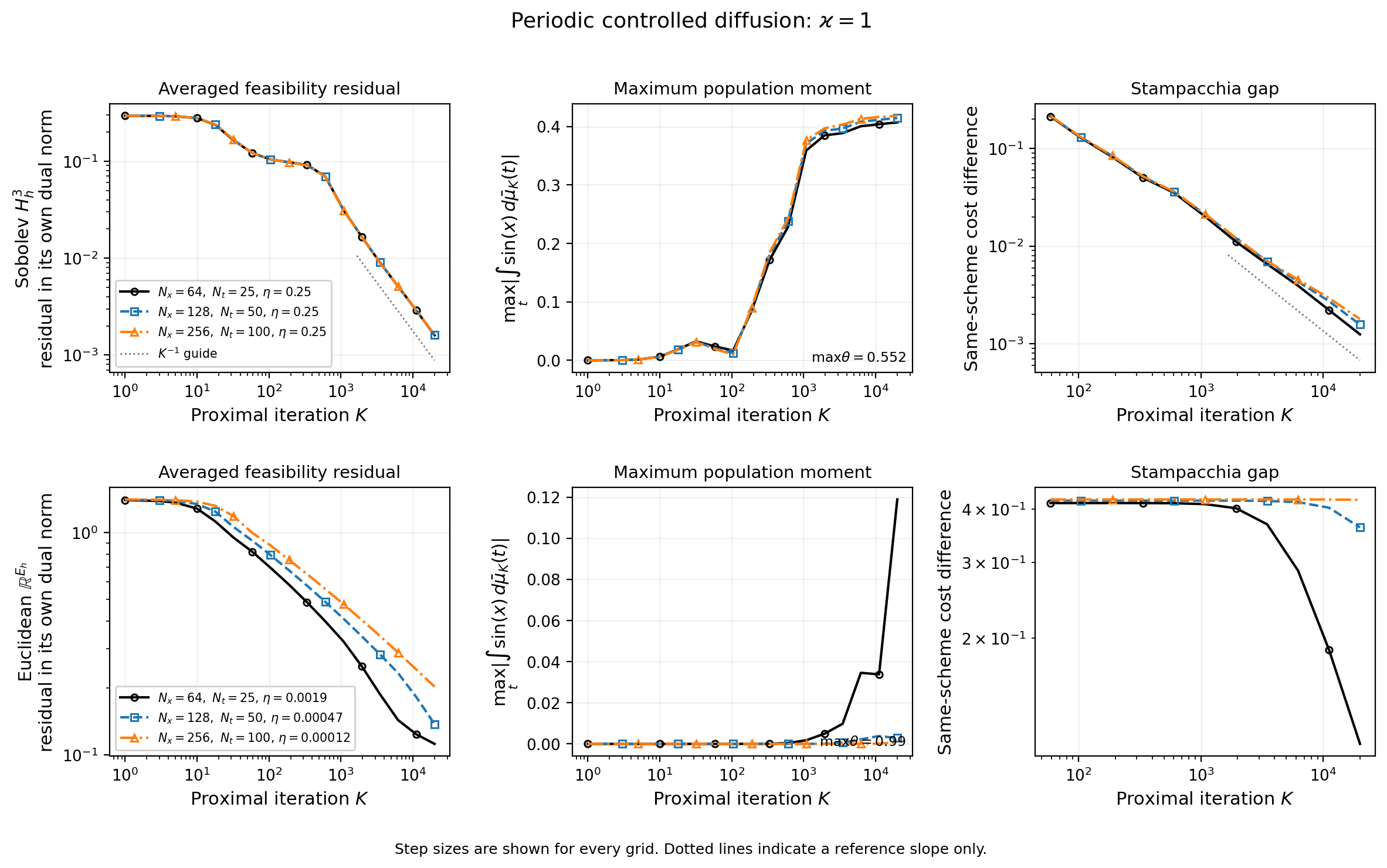}
\caption{Periodic diffusion with $\varkappa=1$. Top: the $H_h^3$
metric with $\lambda=\tau=0.25$ on every mesh. Bottom: the Euclidean
metric with its mesh-dependent admissible step. The columns show the
forward-equation residual, evaluated by quadrature in the corresponding dual norm, the
maximum population moment, and the discrete Stampacchia gap after
re-simulation of the feedback.}
\label{fig:diffusion-grid-g}
\end{figure}

\begin{figure}[t]
\centering
\includegraphics[width=\textwidth]{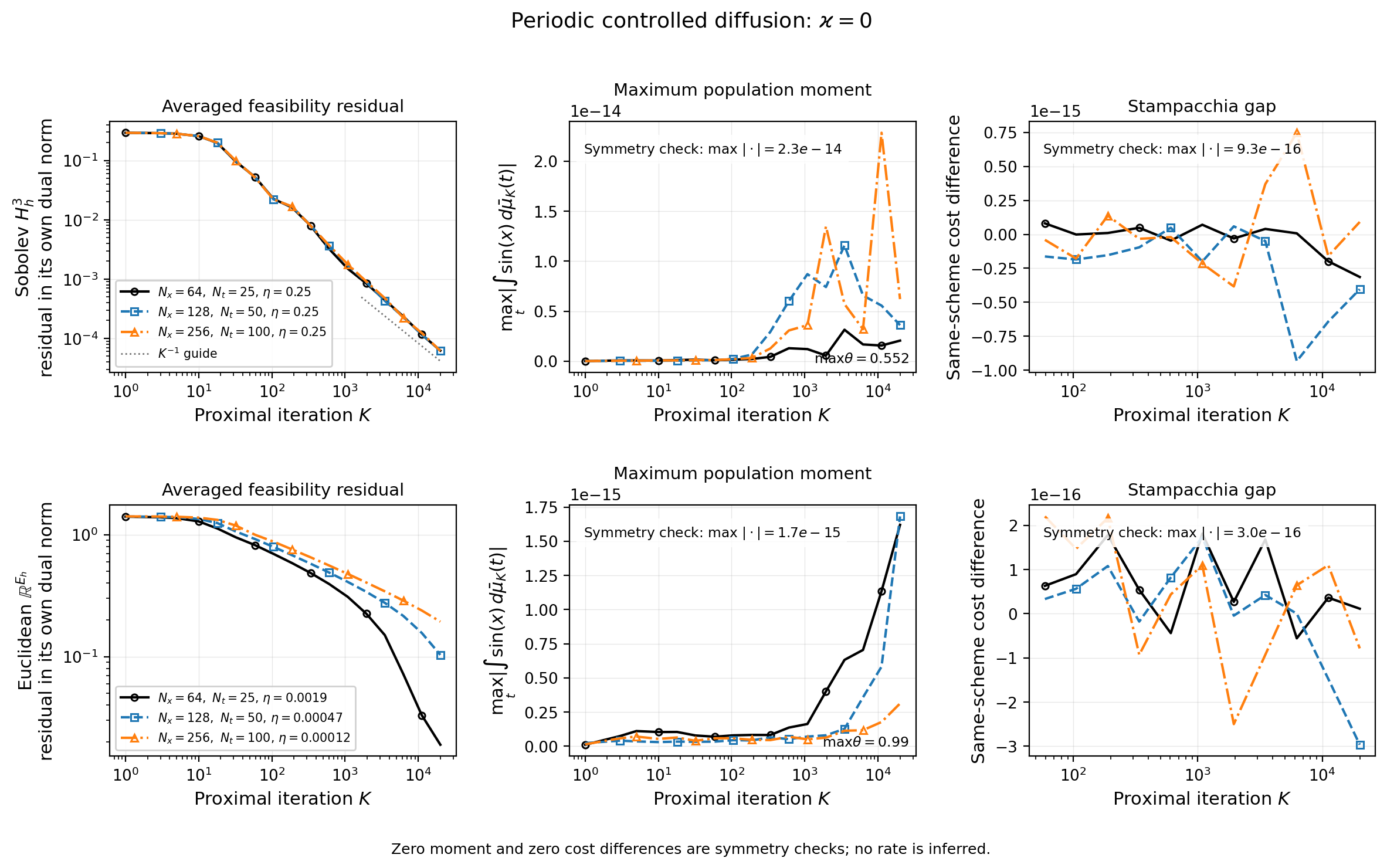}
\caption{Periodic diffusion with $\varkappa=0$ and the same metrics
and steps. The moment and Stampacchia gap panels check the symmetric
zero-cost equilibrium; the residual panel measures the error in the
forward equation.}
\label{fig:diffusion-grid-0}
\end{figure}

\begin{figure}[t]
\centering
\includegraphics[width=\textwidth]{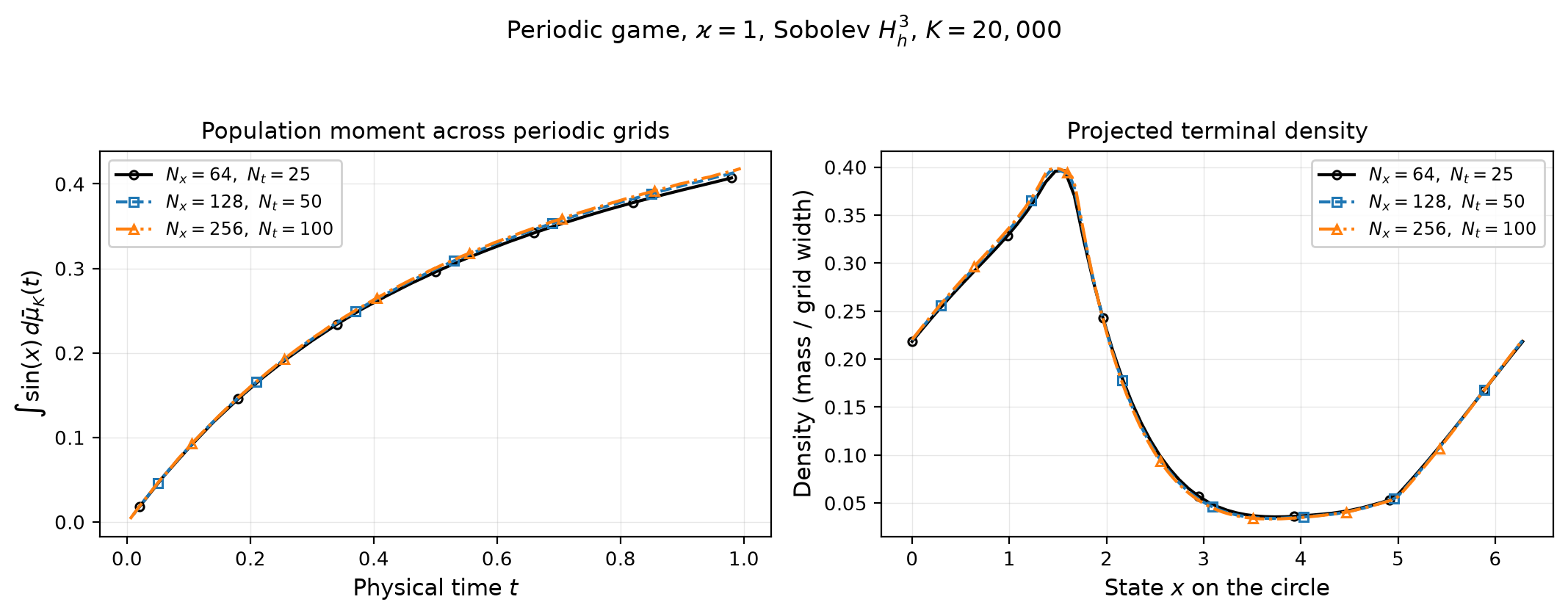}
\caption{Periodic diffusion with $\varkappa=1$ under the $H_h^3$
metric. Population moments and re-simulated terminal densities at the
three resolutions.}
\label{fig:periodic-grid-flow}
\end{figure}

\begin{table}[t]
\centering
\small
\begin{tabular}{cccccc}
\toprule
$(N_\pi,N_t)$ & $\theta_{H_h^3}$ & $\lambda_{\mathrm E}$ &
residual ($H_h^3$) & $\Expl(\hat\pi_K)$ ($H_h^3$) & $\Expl(\hat\pi_K)$ (E)\\
\midrule
$(64,25)$ & $0.5515$ & $1.852\times10^{-3}$ & $1.600\times10^{-3}$ & $1.245\times10^{-3}$ & $0.1133$\\
$(128,50)$ & $0.5480$ & $4.750\times10^{-4}$ & $1.606\times10^{-3}$ & $1.567\times10^{-3}$ & $0.3624$\\
$(256,100)$ & $0.5461$ & $1.202\times10^{-4}$ & $1.616\times10^{-3}$ & $1.784\times10^{-3}$ & $0.4202$\\
\bottomrule
\end{tabular}
\caption{Periodic game with $\varkappa=1$ after $20{,}000$ iterations.
The Sobolev step is $0.25$ on all grids; $\lambda_{\mathrm E}$ is the
admissible Euclidean step. Residuals are measured in each metric's own
dual norm, and Stampacchia gap is evaluated by re-simulation.}
\label{tab:periodic-results}
\end{table}

The exact Fourier operator bounds give $c_H\le0.44943$ and
$\ell_H\le0.51255$ on all three grids.
Consequently the fixed Sobolev step satisfies $\theta<1$ on every grid.
\Cref{tab:periodic-results} shows comparable Sobolev residuals across
resolutions, while the admissible Euclidean step decreases with $h^2$
and the Euclidean runs progress more slowly at the same iteration
count; \Cref{fig:diffusion-grid-g} shows the residual, the
maximum population moment and the re-simulated Stampacchia gap for both
metrics, with the Sobolev runs in the top row and the Euclidean runs in the
bottom row. For $\varkappa=0$ (\Cref{fig:diffusion-grid-0}), the Sobolev residuals at $K=20{,}000$ lie
between $6.15\times10^{-5}$ and $6.22\times10^{-5}$; the sine moments
and the re-simulated Stampacchia gaps are at floating-point roundoff,
as expected from symmetry. \Cref{fig:periodic-grid-flow}
shows, for $\varkappa=1$ and the Sobolev metric, the population moments and
the re-simulated terminal densities at the three resolutions, which nearly coincide.

\subsection{A Fourier--Galerkin alternative}\label{subsec:exp-fourier}

For the same periodic game, we can impose the forward equation against
a finite Fourier basis and evaluate the diffusion generator directly. Take the real Fourier functions
\[
    \phi_0(y)=(2\pi)^{-1/2},\qquad
    \phi_{n,c}(y)=\pi^{-1/2}\cos(ny),\qquad
    \phi_{n,s}(y)=\pi^{-1/2}\sin(ny),\quad 1\le n\le M,
\]
and let $H_M$ be their span with the $H^3(\mathbb T)$ inner product.
Its Gram matrix is diagonal: the constant mode has weight $1$ and the
two modes of frequency $n$ have weight $(1+n^2)^3$. The generator is
evaluated exactly on this basis,
\[
    \LL^a\phi=a\phi'+\phi'',
\]
and the dual correction is computed by dividing each Fourier coefficient
of the integrated forward residual by the corresponding Gram-matrix entry, as in \eqref{eq:dual-step-coords}.
We take $M\in\{4,8,16\}$, $128$ quadrature nodes, $50$ time intervals,
$K=20{,}000$ and $\lambda=\tau=0.25$, with the same uniform reference
as in \Cref{subsec:exp-diffusion}.

\paragraph{Stored measures and updates.}
At each midpoint $t_i$, the code stores the joint state--action masses $q_{k,i,j,a}$ on the nodes $y_j=2\pi j/N$, $N=128$, with $\sum_{j,a}q_{k,i,j,a}=1$, and the terminal masses $\nu_{k,j}$, with $\sum_j\nu_{k,j}=1$. The references are uniform, $q_{0,i,j,a}=1/(2N)$ and $\nu_{0,j}=1/N$. They are quadrature approximations of the uniform circle and action references. The code tilts these mass arrays by \eqref{eq:quadrature-measure-update}, using the Fourier evaluations of $\Gamma$ and the analogous terminal residual. In particular,
\[
    m_{k,i}=\sum_{j,a}q_{k,i,j,a}\sin y_j,\qquad
    \int\LL^a\phi\,d\pi_{k,t_i}
      \approx\sum_{j,a}q_{k,i,j,a}\bigl(a\phi'(y_j)+\phi''(y_j)\bigr).
\]
The state masses are $\mu_{k,i,j}=\sum_aq_{k,i,j,a}$, and the relaxed
feedback is $\alpha_{k,i}(a\mid y_j)=q_{k,i,j,a}/\mu_{k,i,j}$ wherever
the state mass is positive. As in \Cref{ex:finite-dim}, a fixed $M$ imposes the forward equation
only against these Fourier functions.

\begin{figure}[!htbp]
\centering
\includegraphics[width=\textwidth]{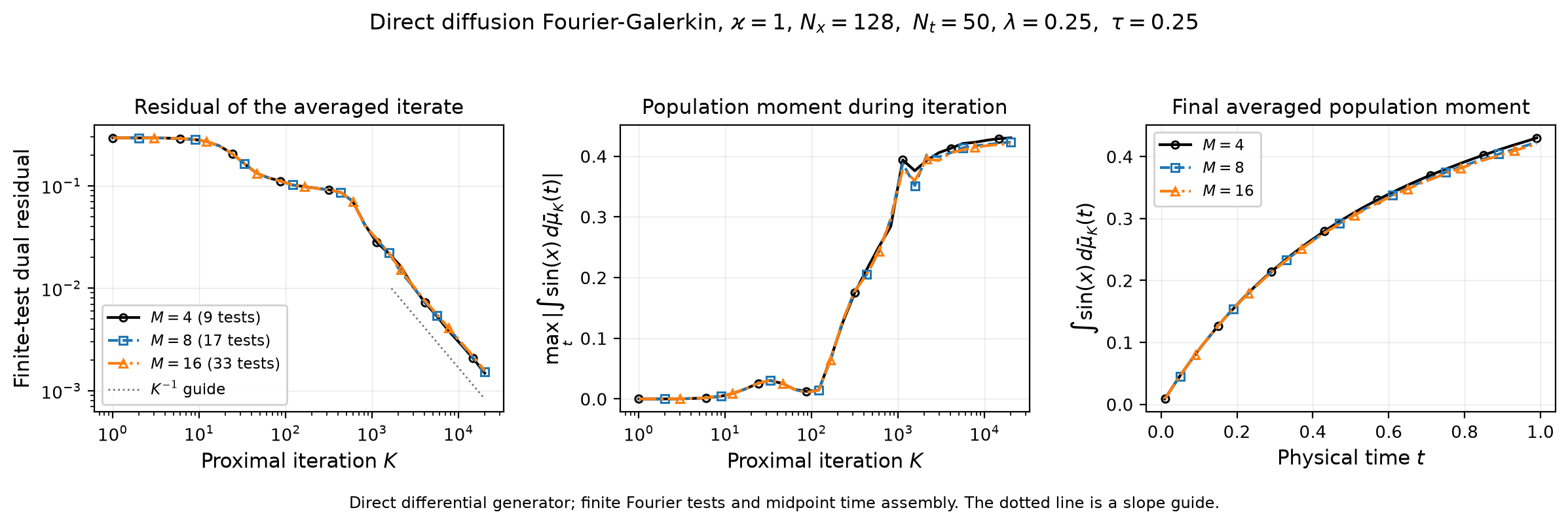}
\caption{Direct Fourier--Galerkin implementation of the periodic game
with $\varkappa=1$ and $M=4,8,16$. The residual is tested in the selected
Fourier space. The population moments use the predictor averages.}
\label{fig:periodic-fourier}
\end{figure}

\Cref{fig:periodic-fourier} shows the residual, tested in the
selected Fourier space, and the population moments for $M=4,8,16$. The averaged residuals at $K=20{,}000$ are
$1.48\times10^{-3}$, $1.53\times10^{-3}$ and $1.57\times10^{-3}$ for
$M=4,8,16$, respectively. The maximum difference between the $M=8$ and
$M=16$ population moments is $2.69\times10^{-3}$. Their distributions
remain sensitive to the test space: testing the $M=8$ result against
$H_{16}$ gives residual $4.04\times10^{-2}$. Thus the small residual in
a fixed Fourier space controls only its selected moments.

\clearpage
\appendix

\section{Proofs for Sections~\ref{sec:setting}, \ref{sec:algo} and~\ref{subsec:inexact}}\label{sec:proofs-setting}

\subsection{Proof of Lemma~\ref{lem:elementary}}
\begin{proof}
Let $u\in\HY$ and $t\in[0,T]$. By Cauchy--Schwarz,
\[
    \nH{u(t)}\le\nH{u(T)}+\int_t^T\nH{\partial_s u}\,ds
    \le\nH{u(T)}+\sqrt T\Bigl(\int_0^T\nH{\partial_s u}^2ds\Bigr)^{1/2}
    \le\sqrt{1+T}\,\tn{u}.
\]

Next we check that $\cA^*u\in\cG$. First, $u(T)\in H\subset C(E)$.
Second, take a strongly measurable representative of
$\partial_t u\in L^2(0,T;H)$. Point evaluations are continuous on $H$, so
$(t,y)\mapsto\partial_t u(t,y)$ is measurable in $t$ and continuous in $y$
for a.e.\ $t$, hence jointly Borel. Approximate $u\in C([0,T];H)$
uniformly by piecewise-constant $H$-valued functions. Since
$\lVert\LL^a_t\phi\rVert_\infty\le\ell_H\nH{\phi}$, the map
$(t,y,a)\mapsto\LL^a_t u(t,\cdot)(y)$ is a uniform limit of Borel maps.
It is therefore Borel, and it is continuous in $(y,a)$ by (S3). For a.e.\ $t$,
\[
    \sup_{y,a}\bigl|\partial_t u(t,y)+\LL^a_t u(t,\cdot)(y)\bigr|
    \le c_H\nH{\partial_t u(t)}+\ell_H\sqrt{1+T}\,\tn{u}.
\]
Hence
\[
    \nG{\cA^*u}^2\le c_H^2\nH{u(T)}^2+2c_H^2\int_0^T\nH{\partial_t u}^2dt+2\ell_H^2T(1+T)\tn{u}^2
    \le L_A^2\tn{u}^2 .
\]
The bound on $\cA$ follows from $|\pair{\cA \pi}{u}|=|\pair{\cA^*u}{\pi}|\le\nG{\cA^*u}\nsharp{\pi}$.
Finally, $|\pair{\rho}{u}|\le\lVert u(0)\rVert_\infty\le c_H\sqrt{1+T}\,\tn{u}$.
\end{proof}

\subsection{Proof of Proposition~\ref{prop:PDVI-MFE}}
\begin{proof}
(i) The second component of the PD-VI reads
$\pair{\cA \pi^*-\rho}{u-u^*}\ge0$ for all $u\in\HY$. Since $\HY$ is a
vector space, this forces $\cA \pi^*=\rho$, so $\pi^*\in\Pirho$. Now take
$z=(\pi,u^*)$ with $\pi\in\Pirho$. Then
$\pair{\cA^*u^*}{\pi-\pi^*}=\pair{\cA \pi-\cA \pi^*}{u^*}=0$, so the first
component gives $\pair{\cC_{\cM \pi^*}}{\pi-\pi^*}\ge0$.

(ii) Write $s^*:=\cC_{\cM \pi^*}-\cA^*u^*=(G^{(\pi^*,u^*)},\Gamma^{(\pi^*,u^*)})$.
Since $\cA \pi^*=\rho$,
\[
    \pair{s^*}{\pi^*}=\pair{\cC_{\cM \pi^*}}{\pi^*}-\pair{\cA \pi^*}{u^*}
    =\pair{\cC_{\cM \pi^*}}{\pi^*}-\pair{\rho}{u^*}=0 .
\]
Hence, for every $(\pi,u)\in\cZ$,
\[
    \pair{\cF(z^*)}{(\pi,u)-z^*}=\pair{s^*}{\pi}-\pair{s^*}{\pi^*}+\pair{\cA \pi^*-\rho}{u-u^*}=\pair{s^*}{\pi}\ge0,
\]
because $s^*\ge0$ and $\nu$, $\pi_t$ are nonnegative measures.

(iii) By \Cref{lem:elementary}, $\Gamma^*$ is a Carath\'eodory integrand
with $t\mapsto\|\Gamma^*(t,\cdot,\cdot)\|_\infty$ in $L^2(0,T)$.
Its minimum $\underline\Gamma$ is measurable (take the infimum over a countable dense
subset of $E\times A$) and lies in $L^2(0,T)$. Thus $c\in H^1(0,T)$,
$c(T)=\underline G$ and $c'=-\underline\Gamma$, so $c1\in\HY$. The HJB residuals after this change are
$G^*-\underline G$ and $\Gamma^*-\underline\Gamma$, hence nonnegative.

The first component of the PD-VI states that $\pi^*$ minimizes
$\pair{s^*}{\pi}$ over $\Pi$, where $s^*:=(G^*,\Gamma^*)$. Since the
terminal law and each running measure have mass one,
\[
    \inf_{\pi\in\Pi}\pair{s^*}{\pi}=\underline G+\int_0^T \underline\Gamma(t)\,dt=c(0).
\]
Indeed, the lower bound is immediate, and measurable approximate
minimizers chosen from a countable dense subset attain the bound in the
limit. By (i), $\cA \pi^*=\rho$, and $\pair{\rho}{c1}=c(0)$. Therefore
\[
    \pair{\cC_{\cM \pi^*}}{\pi^*}-\pair{\rho}{\widehat u^*}
    =\pair{s^*}{\pi^*}-c(0)=0.
\]
Finally, for any $b\in H^1(0,T)$, $\cA^*(b1)=(b(T),-b')$ has zero
pairing with $\pi-\pi^*$ for $\pi,\pi^*\in\Pi$. This proves the stated
invariance of the PD-VI. For example, with zero costs and any feasible
$\pi^*$, the multiplier $u^*\equiv1$ solves the PD-VI but has terminal
residual $-1$, so is not dual feasible.
\end{proof}

\subsection{Proof of Lemma~\ref{lem:riesz}}
\begin{proof}
The map $t\mapsto \LL_t^*\pi_t\in H^*$ is weakly measurable and bounded by
$\ell_H$. Since $H^*$ is separable, it is Bochner integrable, and so is
$s\mapsto d_s$. Let $v\in\HY$. Since $\partial_t v$ does not
depend on $a$,
$\int_0^T\!\int\partial_t v\,d\pi_t\,dt=\int_0^T\pair{\partial_t v(t)}{\mu_t}dt$.
By definition of $\LL_t^*\pi_t$,
$\int_0^T\!\int\LL^a_t v(t)\,d\pi_t\,dt=\int_0^T\pair{v(t)}{\LL_t^*\pi_t}dt$.
Insert $v(t)=v(T)-\int_t^T\partial_s v\,ds$ and use Fubini;
the integrand is bounded by $\ell_H\nH{\partial_s v}\in L^1([0,T]^2)$.
This gives
\[
\begin{aligned}
    \int_0^T\pair{v(t)}{\LL_t^*\pi_t}dt
    &=\Bigl\langle v(T),\int_0^T\LL_t^*\pi_t\,dt\Bigr\rangle-\int_0^T\Bigl\langle\partial_s v(s),\int_0^s\LL_t^*\pi_t\,dt\Bigr\rangle ds,\\
    \pair{v(0)}{\rho}&=\pair{v(T)}{\rho}-\int_0^T\pair{\partial_s v(s)}{\rho}ds .
\end{aligned}
\]
Collecting terms,
\begin{align*}
    \pair{\cA \pi-\rho}{v}
    &=\Bigl\langle v(T),\nu-\rho-\int_0^T\LL_t^*\pi_t\,dt\Bigr\rangle\\
    &\quad+\int_0^T\Bigl\langle\partial_s v(s),
      \rho-\mu_s+\int_0^s\LL_t^*\pi_t\,dt\Bigr\rangle ds\\
    &=\pair{v(T)}{d_T}
      -\int_0^T\pair{\partial_s v(s)}{d_s}\,ds .
\end{align*}
Now define $w(t):=R_Hd_T+\int_t^TR_Hd_s\,ds$. Then $w\in\HY$,
$w(T)=R_Hd_T$ and $\partial_tw=-R_Hd_t$. Hence
$\pp{w}{v}=\pair{v(T)}{d_T}-\int_0^T\pair{\partial_s v}{d_s}ds=\pair{\cA \pi-\rho}{v}$
for all $v\in\HY$, that is, $w=R(\cA \pi-\rho)$.

Moreover, since
\[
    w(T)=R_H d_T,
    \qquad
    \partial_tw(t)=-R_H d_t
    \quad\text{for a.e. }t\in[0,T],
\]
using the definition of $\tn{\cdot}$ and the fact that the
Riesz map $R_H:H^*\to H$ is an isometry,
\[
\begin{aligned}
    \tns{\cA \pi-\rho}^2
    &=\tn{w}^2 \\
    &=\lVert w(T)\rVert_H^2
      +\int_0^T\lVert \partial_tw(t)\rVert_H^2\,dt \\
    &=\lVert R_Hd_T\rVert_H^2
      +\int_0^T\lVert R_Hd_t\rVert_H^2\,dt \\
    &=\lVert d_T\rVert_{H^*}^2
      +\int_0^T\lVert d_t\rVert_{H^*}^2\,dt.
\end{aligned}
\]

The last claim then follows because $\pi\in\Pirho$ if and only if this norm vanishes.
\end{proof}

\subsection{Proof of Lemma~\ref{lem:inexact}}
\begin{proof}
(i) Apply \Cref{lem:PD-three-point} to the corrector with comparison $z'$
and to the predictor with comparison $z_{k+1}$, as in the proof of
\Cref{lem:PD-one-step}, but with the perturbed values
$\hat\cF(z_k):=\cF(z_k)+(0,e_k)$ and $\hat\cF(z_{k+1/2}):=\cF(z_{k+1/2})+(0,e_{k+1/2})$,
for which \eqref{eq:perturbed} is the exact proximal step. Adding gives
\[
    \pair{\hat\cF(z_{k+1/2})}{z_{k+1/2}-z'}
    =\Phi(z',z_k)-\Phi(z',z_{k+1})-\Phi(z_{k+1},z_{k+1/2})-\Phi(z_{k+1/2},z_k)+\hat{\mathcal E}_k,
\]
with
$\hat{\mathcal E}_k=\pair{\hat\cF(z_{k+1/2})-\hat\cF(z_k)}{z_{k+1/2}-z_{k+1}}
=\mathcal E_k+\pair{e_{k+1/2}-e_k}{u_{k+1/2}-u_{k+1}}$, where
$\mathcal E_k$ is the error term from the proof of \Cref{lem:PD-one-step},
bounded there by $\theta[\Phi(z_{k+1/2},z_k)+\Phi(z_{k+1},z_{k+1/2})]$. For
the new term, with $e:=\tn{u_{k+1/2}-u_{k+1}}$, Young's inequality
with weight $1-\theta$ gives
\begin{align*}
    \pair{e_{k+1/2}-e_k}{u_{k+1/2}-u_{k+1}}
    &\le\tns{e_{k+1/2}-e_k}\,e\\
    &\le\frac{\tau(\epsilon_k+\epsilon_{k+1/2})^2}{2(1-\theta)}
      +(1-\theta)\frac{e^2}{2\tau}\\
    &\le\frac{\tau(\epsilon_k+\epsilon_{k+1/2})^2}{2(1-\theta)}
      +(1-\theta)\Phi(z_{k+1},z_{k+1/2}).
\end{align*}
These bounds cancel the term $-\Phi(z_{k+1},z_{k+1/2})$ and leave
$-(1-\theta)\Phi(z_{k+1/2},z_k)$. Finally
\[
\begin{aligned}
\pair{\hat\cF(z_{k+1/2})}{z_{k+1/2}-z'}
&=\pair{\cF(z_{k+1/2})}{z_{k+1/2}-z'}
  +\pair{e_{k+1/2}}{u_{k+1/2}-u'},\\
\bigl|\pair{e_{k+1/2}}{u_{k+1/2}-u'}\bigr|
&\le\epsilon_{k+1/2}\tn{u_{k+1/2}-u'}.
\end{aligned}
\]

(ii) By \Cref{prop:monotone},
$\pair{\cF(z')}{z_{k+1/2}-z'}\le\pair{\cF(z_{k+1/2})}{z_{k+1/2}-z'}$. Average
\eqref{eq:one-step-inexact} over $k<K$, telescope, drop the nonpositive term
and bound $\tn{u_{k+1/2}-u'}\le B_K(u')$.

(iii) Put $z'=z^*$ in \eqref{eq:one-step-inexact}. By
\Cref{ass:LL,prop:monotone} and the PD-VI,
\[
    \pair{\cF(z_{k+1/2})}{z_{k+1/2}-z^*}
    \ge\pair{\cF(z^*)}{z_{k+1/2}-z^*}\ge0.
\]
This uses both monotonicity and the PD-VI. With
$\Phi_k:=\Phi(z^*,z_k)$ and $\delta_k:=\tau(\epsilon_k+\epsilon_{k+1/2})^2/(2(1-\theta))$,
\[
    \Phi_{k+1}\le\Phi_k+\delta_k+\epsilon_{k+1/2}\tn{u_{k+1/2}-u^*}.
\]
Now $u_{k+1/2}-u^*=(u_k-u^*)-\tau\tilde w_k$, and
$\tn{\tilde w_k}\le\tns{\cA \pi_k-\rho}+\epsilon_k=\tns{\cA(\pi_k-\pi^*)}+\epsilon_k\le 2L_A\sqrt{1+T}+\epsilon_k$
by \Cref{lem:elementary} and $\cA \pi^*=\rho$, while
$\tn{u_k-u^*}\le\sqrt{2\tau\Phi_k}$. Hence
$\tn{u_{k+1/2}-u^*}\le\sqrt{2\tau\Phi_k}+2\tau L_A\sqrt{1+T}+\tau\epsilon_k$,
which is the second claim once the first is proved. With $p_k:=\Phi_k+M_\tau$ we
have $2\tau L_A\sqrt{1+T}=\sqrt{2\tau M_\tau}\le\sqrt{2\tau p_k}$, so
\[
    p_{k+1}\le p_k+2\sqrt{2\tau p_k}\,\epsilon_{k+1/2}+\tau\epsilon_k\epsilon_{k+1/2}+\delta_k
    \le\bigl(\sqrt{p_k}+\sqrt{2\tau}\,\epsilon_{k+1/2}\bigr)^2+\tau\epsilon_k\epsilon_{k+1/2}+\delta_k .
\]
Taking square roots, using $\sqrt{a+b}\le\sqrt a+\sqrt b$ and
$\sqrt{\epsilon_k\epsilon_{k+1/2}}\le\epsilon_k+\epsilon_{k+1/2}$,
\[
    \sqrt{p_{k+1}}\le\sqrt{p_k}+\Bigl(\sqrt{2\tau}+\sqrt\tau+\sqrt{\tfrac{\tau}{2(1-\theta)}}\Bigr)(\epsilon_k+\epsilon_{k+1/2}),
\]
and summing over $k$ proves the first claim.
\end{proof}

\subsection{Proof of Lemma~\ref{lem:ritz-error}}
\begin{proof}
Testing the Riesz equation $\pp{Rr}{v}=\pair rv$ on $v\in V_h$ shows that $w_h$ is the orthogonal projection of $Rr$ onto $V_h$ for $\pp{\cdot}{\cdot}$, which is also the minimizer of $\cQ_r$ over the subspace. Pythagoras gives $\tn{Rr}^2=\tn{w_h}^2+\tn{w_h-Rr}^2$, and $\tn{Rr}=\tns r$ by the Riesz isometry; the best-approximation identity is the characterization of the orthogonal projection. In a basis, $w_h=\sum_ic_iv_i$ with $\mathbf G\mathbf c=\mathbf b$, so $\tn{w_h}^2=\mathbf c^\top\mathbf G\mathbf c=\mathbf b^\top\mathbf G^{-1}\mathbf b$.
\end{proof}

\subsection{Proof of Theorem~\ref{thm:inexact}}
\begin{proof}
By \eqref{eq:energy-identity}, each left-hand side in
\eqref{eq:inexact-energy-tolerances} is the squared norm of the
corresponding correction error. Thus \eqref{eq:inexact} holds with
the stated summable tolerances, and \Cref{lem:inexact} applies.
By \Cref{lem:inexact}(iii), both dual sequences are bounded in $\HY$.
The corrector equation gives
\[
 u_K-u_0
 =-\tau R\left(
       K(\cA\bar\pi_K-\rho)+\sum_{k<K}e_{k+1/2}\right).
\]
The Riesz isometry and the summability assumption therefore imply
\[
 \tns{\cA\bar\pi_K-\rho}
 \le\frac1K\left[
       \frac{\tn{u_K-u_0}}{\tau}
       +\sum_{k<K}\epsilon_{k+1/2}\right]
 \le\frac CK,
\]
which proves \eqref{eq:inexact-forward-rate}.

The summable error sequence is also square summable. Moreover,
\[
 B_K(u')
 \le \sup_{k\ge0}\tn{u_{k+1/2}-u_0}+\tn{u'-u_0}.
\]
Inserting this bound into \Cref{lem:inexact}(ii), and bounding the
partial sums by the corresponding infinite sums, proves
\eqref{eq:inexact-main-comparison} with a constant independent of
$K$ and $z'$. 

Finally, take $z'=(\pi',0)$ with $\pi'\in\Pizero$.
Feasibility gives $\cF(\pi',0)=(\cC_{\cM\pi'},0)$, so
\eqref{eq:inexact-main-comparison} becomes the asserted primal
comparison estimate.
\end{proof}

\section{Proofs for Sections~\ref{sec:LL} and~\ref{sec:finite}}\label{sec:proofs-LL}

\subsection{Proof of Proposition~\ref{prop:monotone}}
\begin{proof}
Expanding the definition,
\[
    \pair{\cF(z)-\cF(z')}{z-z'}
    =\pair{\cC_{\cM \pi}-\cC_{\cM \pi'}}{\pi-\pi'}-\pair{\cA^*(u-u')}{\pi-\pi'}+\pair{\cA(\pi-\pi')}{u-u'}.
\]
The last two terms cancel, since $\pair{\cA^*v}{\pi'}=\pair{\cA \pi'}{v}$
by definition. For the Lipschitz bounds, \Cref{ass:Lip} gives directly
\[
    \nG{\cC_{\bm\mu}-\cC_{\bm\mu'}}^2
    \le L_g^2\lVert\mu_T-\mu'_T\rVert_1^2+L_f^2\int_0^T\lVert\mu_t-\mu'_t\rVert_1^2dt
    \le L_C^2\,\nsharp{\bm\mu-\bm\mu'}^2 .
\]
For $\bm\mu=\cM \pi$ and $\bm\mu'=\cM \pi'$, use
$\lVert\mu_t-\mu'_t\rVert_1\le\lVert \pi_t-\pi'_t\rVert_1$, since taking
marginals does not increase total variation, to obtain $\nsharp{\cM \pi-\cM \pi'}\le\nsharp{\pi-\pi'}$.
\end{proof}

\subsection{Proof of Lemma~\ref{lem:PD-three-point}}

We use the static identity of \cite[Lemma~3.1]{KLMP}. Its proof uses only
that the tilting function is bounded and Borel, so we record it in that
form. Let $S$ be a compact metric space, $\mu\in\cP(S)$ a base measure,
$\phi$ bounded Borel on $S$, and
$d\mu^+/d\mu=e^{-\lambda \phi}/\int e^{-\lambda \phi}d\mu$. Then
$|\log\frac{d\mu^+}{d\mu}|\le2\lambda\lVert \phi\rVert_\infty$, so
$\KL(\mu^+\|\mu)\le2\lambda\lVert \phi\rVert_\infty$. For every $\mu'$ with
$\KL(\mu'\|\mu)<\infty$, we have
$\KL(\mu'\|\mu^+)\le\KL(\mu'\|\mu)+2\lambda\lVert \phi\rVert_\infty$ and
\begin{equation}\label{eq:static-3pt}
    \lambda\pair{\phi}{\mu^+-\mu'}=\KL(\mu'\|\mu)-\KL(\mu'\|\mu^+)-\KL(\mu^+\|\mu).
\end{equation}

\begin{proof}[Proof of Lemma~\ref{lem:PD-three-point}]
Let $\gamma=(G,\Gamma)$. For a.e.\ $t$, $\Gamma_t$ is bounded, and on the
remaining null set we set $\pi^+_t:=\pi_t$. Then $t\mapsto \pi^+_t$ is a
measurable kernel and $\pi^+\in\Pi$. By the static facts,
\[
    \bD(\pi^+\|\pi)\le2\lambda\Bigl(\lVert G\rVert_\infty+\int_0^T\lVert\Gamma_t\rVert_\infty dt\Bigr)<\infty .
\]
Let $\pi'$ satisfy $\bD(\pi'\|\pi)<\infty$. Then
$\bD(\pi'\|\pi^+)\le\bD(\pi'\|\pi)+2\lambda(\lVert G\rVert_\infty+\int\lVert\Gamma_t\rVert_\infty)<\infty$.
Apply \eqref{eq:static-3pt} to $\nu$ with $\phi=G$ and to $\pi_t$ with
$\phi=\Gamma_t$ for a.e.\ $t$. Every term is integrable in $t$, so we may
integrate and divide by $\lambda$:
\[
    \pair{\gamma}{\pi^+-\pi'}=\tfrac1\lambda\bigl[\bD(\pi'\|\pi)-\bD(\pi'\|\pi^+)-\bD(\pi^+\|\pi)\bigr].
\]
For the dual part, $u^+-u=-\tau R\ell$. Polarization gives
\[
    \pair{\ell}{u^+-u'}=\tfrac1\tau\pp{u-u^+}{u^+-u'}
    =\tfrac1{2\tau}\bigl[\tn{u'-u}^2-\tn{u'-u^+}^2-\tn{u^+-u}^2\bigr].
\]
Adding the two identities proves the three-point identity. Rearranged, it
reads
\[
    \pair{\gamma}{\pi'}+\pair{\ell}{u'}+\Phi(z',z)
    =\pair{\gamma}{\pi^+}+\pair{\ell}{u^+}+\Phi(z^+,z)+\Phi(z',z^+).
\]
If $\bD(\pi'\|\pi)=\infty$, the left side is $+\infty$. Moreover
$\Phi(z',z^+)=0$ if and only if $z'=z^+$. So $z^+$ is the unique minimizer.
\end{proof}

\subsection{Proof of Lemma~\ref{lem:PD-one-step}}
\begin{proof}
Both $\pi_{k+1/2}$ and $\pi_{k+1}$ are tilts of $\pi_k$. So
\Cref{lem:PD-three-point} and induction give
$\bD(\pi'\|\pi_k)<\infty$ and $\bD(\pi'\|\pi_{k+1/2})<\infty$ for every $k$, as well
as $\bD(\pi_{k+1}\|\pi_k)<\infty$.

Apply \Cref{lem:PD-three-point} to the corrector with comparison $z'$:
\[
    \pair{\cF(z_{k+1/2})}{z_{k+1}-z'}=\Phi(z',z_k)-\Phi(z',z_{k+1})-\Phi(z_{k+1},z_k).
\]
Apply it to the predictor with comparison $z_{k+1}$:
\[
    \pair{\cF(z_k)}{z_{k+1/2}-z_{k+1}}=\Phi(z_{k+1},z_k)-\Phi(z_{k+1},z_{k+1/2})-\Phi(z_{k+1/2},z_k).
\]
Adding, and using
$\pair{\cF(z_{k+1/2})}{z_{k+1/2}-z'}=\pair{\cF(z_{k+1/2})}{z_{k+1}-z'}+\pair{\cF(z_k)}{z_{k+1/2}-z_{k+1}}+\mathcal E_k$,
we get
\[
    \pair{\cF(z_{k+1/2})}{z_{k+1/2}-z'}
    =\Phi(z',z_k)-\Phi(z',z_{k+1})-\Phi(z_{k+1},z_{k+1/2})-\Phi(z_{k+1/2},z_k)+\mathcal E_k,
\]
with the error term
\[
\begin{aligned}
    \mathcal E_k&:=\pair{\cF(z_{k+1/2})-\cF(z_k)}{z_{k+1/2}-z_{k+1}}\\
    &=\pair{\cC_{\cM \pi_{k+1/2}}-\cC_{\cM \pi_k}-\cA^*(u_{k+1/2}-u_k)}{\pi_{k+1/2}-\pi_{k+1}}
    +\pair{\cA(\pi_{k+1/2}-\pi_k)}{u_{k+1/2}-u_{k+1}} .
\end{aligned}
\]
Set $a:=\nsharp{\pi_{k+1/2}-\pi_k}$, $b:=\tn{u_{k+1/2}-u_k}$,
$c:=\nsharp{\pi_{k+1/2}-\pi_{k+1}}$ and $e:=\tn{u_{k+1/2}-u_{k+1}}$. By
\Cref{prop:monotone,lem:elementary},
$\mathcal E_k\le L_Cac+L_Abc+L_Aae$. Young's inequality gives
\[
    L_Cac\le\lambda L_C\,\tfrac{a^2+c^2}{2\lambda},\qquad
    L_Abc\le\sqrt{\lambda\tau}L_A\Bigl(\tfrac{b^2}{2\tau}+\tfrac{c^2}{2\lambda}\Bigr),\qquad
    L_Aae\le\sqrt{\lambda\tau}L_A\Bigl(\tfrac{a^2}{2\lambda}+\tfrac{e^2}{2\tau}\Bigr).
\]
Hence
\[
    \mathcal E_k\le\theta\Bigl[\tfrac{a^2}{2\lambda}+\tfrac{b^2}{2\tau}+\tfrac{c^2}{2\lambda}+\tfrac{e^2}{2\tau}\Bigr]
    \le\theta\bigl[\Phi(z_{k+1/2},z_k)+\Phi(z_{k+1},z_{k+1/2})\bigr],
\]
where the last step uses Pinsker's inequality \eqref{eq:pinsker-dyn}.
Substituting proves the claim.
\end{proof}

\subsection{Proof of Theorem~\ref{thm:PD-ergodic-minty}}
\begin{proof}
Since $\theta\le1$, \Cref{lem:PD-one-step} gives
$\pair{\cF(z_{k+1/2})}{z_{k+1/2}-z'}\le\Phi(z',z_k)-\Phi(z',z_{k+1})$. By
\Cref{prop:monotone} and \Cref{ass:LL},
$\pair{\cF(z')}{z_{k+1/2}-z'}\le\pair{\cF(z_{k+1/2})}{z_{k+1/2}-z'}$.
Summing over $k=0,\ldots,K-1$ and using $\Phi(z',z_K)\ge0$,
\[
    \sum_{k=0}^{K-1}\pair{\cF(z')}{z_{k+1/2}-z'}\le\Phi(z',z_0).
\]
The map $w\mapsto\pair{\cF(z')}{w-z'}$ is affine, so the left side equals
$K\pair{\cF(z')}{\bar z_K-z'}$.
\end{proof}

\subsection{Proof of Corollary~\ref{cor:primal-minty-residual}}
\begin{proof}
(a) Let $\pi'\in\Pizero$ and $u'\in\HY$. Since $\cA \pi'=\rho$,
\eqref{eq:gap-expansion} and \Cref{thm:PD-ergodic-minty} give
\[
    \pair{\cC_{\cM \pi'}}{\bar \pi_K-\pi'}-\pair{\cA\bar\pi_K-\rho}{u'}
    \le\frac{\bD(\pi'\|\pi_0)}{\lambda K}+\frac{\tn{u'-u_0}^2}{2\tau K}.
\]
Choose $u'=u_0-\tau KR (\cA\bar\pi_K-\rho)$. Then
\[
    -\pair{\cA\bar\pi_K-\rho}{u'}=-\pair{\cA\bar\pi_K-\rho}{u_0}+\tau K\tns{\cA\bar\pi_K-\rho}^2,
    \qquad
    \tn{u'-u_0}^2=\tau^2K^2\tns{\cA\bar\pi_K-\rho}^2.
\]
This yields the first inequality in \eqref{eq:master}. For the second inequality, apply
\[
\pair{\cA\bar\pi_K-\rho}{u_0}-\frac{\tau K}{2}\tns{\cA\bar\pi_K-\rho}^2\le\tn{u_0}\tns{\cA\bar\pi_K-\rho}-\frac{\tau K}{2}\tns{\cA\bar\pi_K-\rho}^2\le\frac{\tn{u_0}^2}{2\tau K}.
\]

(b) Since $\pair{\cC_{\cM \pi'}}{\bar \pi_K-\pi'}\ge-\nG{\cC_{\cM \pi'}}\nsharp{\bar \pi_K-\pi'}\ge-2C_0\sqrt{1+T}$, by \eqref{eq:master}, the residual $r:=\tns{\cA\bar\pi_K-\rho}$ satisfies $\frac{\tau K}2r^2\le a_0+b_0\,r$ with $a_0=\bD(\pi'\|\pi_0)/(\lambda K)+2C_0\sqrt{1+T}$ and $b_0=\tn{u_0}$. For $a>0$ and $a_0,b_0\ge0$, the inequality $a r^2\le a_0+b_0r$ implies $r\le b_0/a+\sqrt{a_0/a}$. Taking $a=\tau K/2$ gives the claim.

For the sharper bound, by \Cref{prop:PDVI-MFE}(ii), $(\pi^*,u^*)$ solves the PD-VI. Testing with $z=(\bar \pi_K,u^*)$ and using $\cA \pi^*=\rho$,
\[
    \pair{\cC_{\cM \pi^*}}{\bar \pi_K-\pi^*}\ge\pair{\cA(\bar \pi_K-\pi^*)}{u^*}=\pair{\cA\bar\pi_K-\rho}{u^*}\ge-\tn{u^*}\tns{\cA\bar\pi_K-\rho}.
\]
Put $\pi'=\pi^*$ in \eqref{eq:master} and apply the same quadratic estimate, with $a_0=\bD(\pi^*\|\pi_0)/(\lambda K)$ and $b_0=\tn{u_0}+\tn{u^*}$.
\end{proof}

\subsection{Proof of Lemma~\ref{lem:topology}}
\begin{proof}
For $\pi\in\Pi$ let $\bm\pi(dt,dy,da):=dt\,\pi_t(dy,da)$ be its running family viewed as a single finite measure on $[0,T]\times E\times A$, so that $\pi^n\Rightarrow\pi$ is weak convergence of $(\nu^n,\bm\pi^n)$ to $(\nu,\bm\pi)$. Let $\mathfrak M_T$ be the set of finite nonnegative Borel measures on $[0,T]\times E\times A$ whose $t$-marginal is Lebesgue measure on $[0,T]$. By the disintegration theorem, every $\bm\pi\in\mathfrak M_T$ has the form $\bm\pi(dt,dy,da)=dt\,\pi_t(dy,da)$ for a measurable kernel $(\pi_t)$, unique up to $dt$-null sets. So $\pi\mapsto(\nu,\bm\pi)$ is a bijection from $\Pi$ onto $\cP(E)\times\mathfrak M_T$. Measures of total mass $T$ on the compact metric space $[0,T]\times E\times A$ form a compact metrizable set for the weak topology. The marginal constraint is weakly closed, so $\mathfrak M_T$ is compact metrizable, and so is $\cP(E)$.

Let $\pi^n\Rightarrow \pi$ and let $\Gamma$ and
$\Lambda(t):=\lVert\Gamma_t\rVert_\infty\in L^1$ be as in the statement.
Changing $\Gamma$ on a $dt$-null set of times, where we set $\Gamma_t:=0$,
changes none of the integrals, because all $t$-marginals are Lebesgue. We
may therefore assume $\Gamma_t\in C(E\times A)$ for every $t$.
Fix $\epsilon>0$ and choose $M$ with $\int_{\{\Lambda>M\}}\Lambda\,dt<\epsilon$.
The truncation $\Gamma^M:=(\Gamma\wedge M)\vee(-M)$ is still Carath\'eodory.
Because all $t$-marginals equal Lebesgue measure,
\[
    \Bigl|\int_0^T\!\!\int(\Gamma_t-\Gamma^M_t)\,d\pi^n_t\,dt\Bigr|\le\epsilon
    \quad\text{for all }n,
\]
and likewise for $\pi$. By the Scorza-Dragoni theorem (see, e.g.,
\cite{FonsecaLeoni2007}), there is a compact $T_\epsilon\subset[0,T]$ with
$|[0,T]\setminus T_\epsilon|<\epsilon/M$ such that $\Gamma^M$ is continuous
on $T_\epsilon\times E\times A$. Extend it by Tietze to a continuous
$\tilde\Gamma$ on $[0,T]\times E\times A$ with $|\tilde\Gamma|\le M$. Then
$|\int\!\int(\Gamma^M-\tilde\Gamma)\,d\pi^n_tdt|\le2\epsilon$ for all $n$, and
likewise for $\pi$. Finally,
$\int\tilde\Gamma\,d\bm\pi^n\to\int\tilde\Gamma\,d\bm\pi$ by weak convergence. Since
$\epsilon$ is arbitrary, the first claim follows.

For $\gamma=(G,\Gamma)\in\cG$, we have $G\in C(E)$ and
$\int\lVert\Gamma_t\rVert_\infty dt\le\sqrt T\nG{\gamma}$, so
$\pi\mapsto\pair{\gamma}{\pi}$ is continuous. Finally,
$\pair{\cA \pi-\rho}{u}=\pair{\cA^*u}{\pi}-\pair{\rho}{u}$ with
$\cA^*u\in\cG$ by \Cref{lem:elementary}.
\end{proof}

\subsection{Proof of Theorem~\ref{thm:cluster-points}}
\begin{proof}
(i) By \Cref{lem:topology}, limit points exist. Let
$\bar \pi_{K_j}\Rightarrow\bar \pi$. For
$u\in\HY$, \Cref{lem:topology} and \Cref{cor:primal-minty-residual}(b) give
\[
    |\pair{\cA\bar \pi-\rho}{u}|=\lim_j|\pair{\cA\bar\pi_{K_j}-\rho}{u}|\le\lim_j\tns{\cA\bar\pi_{K_j}-\rho}\tn{u}=0,
\]
so $\bar \pi\in\Pirho$. For $\pi'\in\Pizero$, the cost vector
$\cC_{\cM \pi'}$ belongs to $\cG$. Passing to the limit in
\Cref{cor:primal-minty-residual}(a) with \Cref{lem:topology} gives
$\pair{\cC_{\cM \pi'}}{\bar \pi-\pi'}\le0$.

(ii) Let $\pi'\in\Pirho$ and let $\pi'_n\in\Pizero$ be as in
\Cref{ass:density}. By (i), $\pair{\cC_{\cM \pi'_n}}{\bar \pi-\pi'_n}\le0$. Split
\[
    \pair{\cC_{\cM \pi'_n}}{\bar \pi-\pi'_n}-\pair{\cC_{\cM \pi'}}{\bar \pi-\pi'}
    =\pair{\cC_{\cM \pi'_n}-\cC_{\cM \pi'}}{\bar \pi-\pi'_n}+\pair{\cC_{\cM \pi'}}{\pi'-\pi'_n}.
\]
The second term tends to $0$ by \Cref{lem:topology}, since
$\pi'_n\Rightarrow \pi'$. By \Cref{prop:monotone}, the first term is bounded by
\[
    \nG{\cC_{\cM \pi'_n}-\cC_{\cM \pi'}}\,\nsharp{\bar \pi-\pi'_n}
    \le L_C\,\nsharp{\cM \pi'_n-\cM \pi'}\,2\sqrt{1+T}\to0 .
\] Hence $\pair{\cC_{\cM \pi'}}{\bar \pi-\pi'}\le0$ for all $\pi'\in\Pirho$, that is, $\GapM(\bar\pi)=0$. By \Cref{prop:eps-relations}(ii), $\bar \pi$ is an MFE.

Suppose now the MFE $\pi^*$ is unique but $\bar \pi_K\not\Rightarrow \pi^*$. By
compactness there is a subsequence converging to some $\bar \pi\ne \pi^*$. By
what we just proved, $\bar \pi$ is an MFE, a contradiction.
\end{proof}

\subsection{Proof of Proposition~\ref{prop:eps-relations}}
\begin{proof}
(i) For $\pi'\in\Pirho$, \Cref{ass:LL} gives
\[
    \pair{\cC_{\cM \pi'}}{\pi-\pi'}
    =\pair{\cC_{\cM \pi}}{\pi-\pi'}-\pair{\cC_{\cM \pi}-\cC_{\cM \pi'}}{\pi-\pi'}
    \le\pair{\cC_{\cM \pi}}{\pi-\pi'} .
\]
Take the supremum over $\pi'\in\Pirho$.

(ii) We adapt \cite[Prop.~3.2(i)]{bigi2023approximate} and
\cite[Prop.~3.5]{KLMP}. First suppose $\varepsilon>0$. Let $\delta:=4\sqrt{(1+T)L_C\,\varepsilon}$ and
suppose there is $\tilde \pi\in\Pirho$ with
$\pair{\cC_{\cM \pi}}{\tilde \pi-\pi}<-\delta$. Let
$s_0:=(\varepsilon/(4L_C(1+T)))^{1/2}\in(0,1]$ and
$\pi_{s_0}:=(1-s_0)\pi+s_0\tilde \pi\in\Pirho$. Taking $\pi'=\pi_{s_0}$ in
$\GapM(\pi)\le\varepsilon$ and using \Cref{prop:monotone},
\[
    -\varepsilon\le s_0\pair{\cC_{\cM \pi_{s_0}}}{\tilde \pi-\pi}
    <-\delta s_0+s_0L_C\nsharp{\pi_{s_0}-\pi}\nsharp{\tilde \pi-\pi}
    \le-\delta s_0+4L_C(1+T)s_0^2=-2\varepsilon+\varepsilon,
\]
because $\delta s_0=2\varepsilon$ and $4L_C(1+T)s_0^2=\varepsilon$. This is
a contradiction, so $\Expl(\pi)\le\delta$. If $\GapM(\pi)=0$, this holds for
every $\varepsilon\in(0,4L_C(1+T)]$, and letting $\varepsilon\downarrow0$
gives $\Expl(\pi)=0$.
\end{proof}

\subsection{Proof of Theorem~\ref{thm:quantitative-eps-MFE}}
\begin{proof}
By \Cref{ass:entropy-bound}, $\Pizero=\Pirho$ and
$\bD(\pi'\|\pi_0)\le\bar D$ for all $\pi'\in\Pirho$.

(i) Put $\pi'=\hat \pi_K\in\Pirho$ in \eqref{eq:master}. By
\Cref{ass:projection},
\[
    \pair{\cC_{\cM\hat \pi_K}}{\bar \pi_K-\hat \pi_K}\ge-C_0\nsharp{\bar \pi_K-\hat \pi_K}\ge-C_0C_{\mathrm{proj}}\tns{\cA\bar\pi_K-\rho}.
\]
Hence $\frac{\tau K}{2}\tns{\cA\bar\pi_K-\rho}^2\le\frac{\bar D}{\lambda K}+(\tn{u_0}+C_0C_{\mathrm{proj}})\tns{\cA\bar\pi_K-\rho}$.
The quadratic argument from the proof of
\Cref{cor:primal-minty-residual}(b) gives (i).

(ii) This follows from (i) and \Cref{ass:projection}.

(iii) For $\pi'\in\Pirho$, \Cref{cor:primal-minty-residual}(a) and (ii) give
\[
    \pair{\cC_{\cM \pi'}}{\hat \pi_K-\pi'}
    =\pair{\cC_{\cM \pi'}}{\bar \pi_K-\pi'}+\pair{\cC_{\cM \pi'}}{\hat \pi_K-\bar \pi_K}
    \le\frac{\bar D}{\lambda K}+\frac{\tn{u_0}^2}{2\tau K}+C_0C_{\mathrm{proj}}\omega_K=\varepsilon_K .
\]
Take the supremum over $\pi'\in\Pirho$.

(iv) This follows from (iii) and \Cref{prop:eps-relations}(ii).

For $u_0=0$, $\varepsilon_K=\frac1K\bigl[\frac{\bar D}{\lambda}+\frac{2C_0^2C_{\mathrm{proj}}^2}{\tau}+C_0C_{\mathrm{proj}}\sqrt{2\bar D/(\lambda\tau)}\bigr]$. The inequality $\sqrt{2ab}\le a+b/2$, applied with $a=\bar D/\lambda$ and $b=C_0^2C_{\mathrm{proj}}^2/\tau$, gives the stated bound. For the stated number of iterations,
$\varepsilon_K\le\min\{4L_C(1+T),\,\delta^2/(16(1+T)L_C)\}$, so (iv) gives $\Expl(\hat\pi_K)\le4\sqrt{(1+T)L_C\,\delta^2/(16(1+T)L_C)}=\delta$.
\end{proof}

\subsection{Proof of Proposition~\ref{prop:cost-correction}}
\label{subsec:proof-cost-correction}
\begin{proof}
The entropy assumption gives $\Pizero=\Pirho$. Use the feasible comparison $\hat\pi_K=\projr(\bar\pi_K)$ in
\eqref{eq:master}. By \eqref{eq:cost-correction}, with $q=\hat\pi_K$,
\[
 \pair{\cC_{\cM\hat\pi_K}}{\bar\pi_K-\hat\pi_K}
 \ge-C_{\mathrm c}\tns{\cA\bar\pi_K-\rho}^{\alpha}.
\]
Therefore
\[
 \frac{\tau K}{2}\tns{\cA\bar\pi_K-\rho}^2
 \le\frac{\bar D}{\lambda K}
       +\tn{u_0}\tns{\cA\bar\pi_K-\rho}+C_{\mathrm c}\tns{\cA\bar\pi_K-\rho}^{\alpha}.
\]
If $\tns{\cA\bar\pi_K-\rho}>0$, at least one of the three terms on the right is at
least $\tau K\tns{\cA\bar\pi_K-\rho}^2/6$. Solving the corresponding inequality
gives one of the three bounds defining $\omega_K$, so
$\tns{\cA\bar\pi_K-\rho}\le\omega_K$. The conclusion also holds when $\tns{\cA\bar\pi_K-\rho}=0$.

For every $q\in\Pirho$, \Cref{cor:primal-minty-residual}(a) gives
\[
 \pair{\cC_{\cM q}}{\bar\pi_K-q}
 \le\frac{\bar D}{\lambda K}+\frac{\tn{u_0}^2}{2\tau K}.
\]
Add the correction and apply \eqref{eq:cost-correction} to obtain
\[
 \pair{\cC_{\cM q}}{\hat\pi_K-q}
 \le\frac{\bar D}{\lambda K}+\frac{\tn{u_0}^2}{2\tau K}
       +C_{\mathrm c}\tns{\cA\bar\pi_K-\rho}^{\alpha}
 \le\varepsilon_K.
\]
Taking the supremum over $q$ proves the Minty-gap estimate. The
approximate-equilibrium bound follows from
\Cref{prop:eps-relations}(ii). Since $0<\alpha\le1$, the slowest
decay in $\omega_K$ is $K^{-1/(2-\alpha)}$, and the slowest decay
in $\varepsilon_K$ is $K^{-\alpha/(2-\alpha)}$, giving the stated
orders.
\end{proof}

\subsection{Proof of Theorem~\ref{thm:strong-LL}}
\begin{proof}
(i) Since $\bD(\pi^*\|\pi_0)<\infty$ and $z^*\in\cZ$, \Cref{lem:PD-one-step}
applies with $z'=z^*$:
\[
    \pair{\cF(z_{k+1/2})}{z_{k+1/2}-z^*}
    \le\Phi(z^*,z_k)-\Phi(z^*,z_{k+1})
    -(1-\theta)\bigl[\Phi(z_{k+1/2},z_k)+\Phi(z_{k+1},z_{k+1/2})\bigr].
\]
By the PD-VI and \Cref{prop:monotone,ass:strong-LL},
\[
    \pair{\cF(z_{k+1/2})}{z_{k+1/2}-z^*}
    \ge\kappa\lVert\cM\pi_{k+1/2}-\bm\mu^*\rVert_{W_1,2}^2.
\]
Inserting this lower
bound into the one-step estimate and rearranging gives (i). Since
$\theta\le1$, every subtracted term is nonnegative, so
$(\Phi(z^*,z_k))_k$ is nonincreasing, and
$\tn{u_k-u^*}^2/(2\tau)\le\Phi(z^*,z_k)\le\Phi(z^*,z_0)$.

(ii) Summing (i) over $k=0,\ldots,K-1$, dropping the $(1-\theta)$ terms
and using $\Phi(z^*,z_K)\ge0$,
\[
    \kappa\sum_{k=0}^{K-1}\lVert\cM \pi_{k+1/2}-\bm\mu^*\rVert_{W_1,2}^2
    \le\Phi(z^*,z_0)-\Phi(z^*,z_K)\le\Phi(z^*,z_0).
\]
Letting $K\to\infty$ gives the series bound, and the bound on the minimum
follows from $K\min_{k<K}\le\sum_{k<K}$. For the average, $\cM$ is affine,
so $\cM\bar \pi_K=\frac1K\sum_{k=0}^{K-1}\cM \pi_{k+1/2}$. For each $t$ the
map $\mu\mapsto W_1(\mu,\mu^*_t)$ is convex on $\cP(E)$, being the norm of
$\mu-\mu^*_t$ in the dual of the $1$-Lipschitz functions, and the square of a
nonnegative convex function is convex. Hence
$\pi\mapsto\lVert\cM \pi-\bm\mu^*\rVert_{W_1,2}^2$ is convex on $\Pi$, and
Jensen's inequality gives
\[
    \lVert\cM\bar \pi_K-\bm\mu^*\rVert_{W_1,2}^2
    \le\frac1K\sum_{k=0}^{K-1}\lVert\cM \pi_{k+1/2}-\bm\mu^*\rVert_{W_1,2}^2
    \le\frac{\Phi(z^*,z_0)}{\kappa K}.
\]

(iii) The series in (ii) converges, so its terms tend to zero. If
$\theta<1$, summing (i) also gives
$(1-\theta)\sum_{k\ge0}\Phi(z_{k+1/2},z_k)\le\Phi(z^*,z_0)$, hence
$\bD(\pi_{k+1/2}\|\pi_k)\le\lambda\Phi(z_{k+1/2},z_k)\to0$, and Pinsker's
inequality \eqref{eq:pinsker-dyn} gives $\nsharp{\pi_{k+1/2}-\pi_k}\to0$. By the
triangle inequality for $\lVert\cdot\rVert_{W_1,2}$ and \eqref{eq:flat-vs-d1},
\[
\begin{aligned}
    \lVert\cM \pi_k-\bm\mu^*\rVert_{W_1,2}
    &\le\lVert\cM \pi_k-\cM \pi_{k+1/2}\rVert_{W_1,2}
       +\lVert\cM \pi_{k+1/2}-\bm\mu^*\rVert_{W_1,2}\\
    &\le\tfrac{D_E}{2}\nsharp{\pi_k-\pi_{k+1/2}}
       +\lVert\cM \pi_{k+1/2}-\bm\mu^*\rVert_{W_1,2}\to0 .
\end{aligned}
\]

(iv) Let $\pi^*,\pi^{**}\in\Pirho$ be MFE. \Cref{def:MFE} with the competitor
$\pi^{**}$ for $\pi^*$ and $\pi^*$ for $\pi^{**}$ gives
$\pair{\cC_{\cM \pi^*}}{\pi^*-\pi^{**}}\le0$ and
$\pair{\cC_{\cM \pi^{**}}}{\pi^{**}-\pi^*}\le0$. Adding,
$\pair{\cC_{\cM \pi^*}-\cC_{\cM \pi^{**}}}{\pi^*-\pi^{**}}\le0$, so
\Cref{ass:strong-LL} yields $\kappa\lVert\cM \pi^*-\cM \pi^{**}\rVert_{W_1,2}^2\le0$.
Thus $W_1(\mu^*_t,\mu^{**}_t)=0$ for a.e.\ $t$ and $W_1(\mu^*_T,\mu^{**}_T)=0$,
that is, $\cM \pi^*=\cM \pi^{**}$. Finally, the primal part of every PD-VI
solution is an MFE by \Cref{prop:PDVI-MFE}(i), so the flow $\bm\mu^*$ in
(ii)--(iii) is this common equilibrium flow.
\end{proof}

\subsection{Proof of Theorem~\ref{thm:convex-diffusion}}\label{subsec:proof-convex-diffusion}
\begin{proof}
All constants below are independent of $K$ and of the smoothing
parameters. The comparison of control costs in Steps~1--2 follows the
method of \cite[Section~3.2]{BertucciCecchin2024}, using
\eqref{eq:master} and comparison measures with finite entropy. Steps~3--4 estimate the
additional forward-equation residual of the averaged iterate.
We write $\Gamma^*$ for the running HJB residual,
\[
 \Gamma^*(t,y,a)
 :=f(t,y,a,\mu_t^*)+\partial_tu^*(t,y)
      +a\,\partial_yu^*(t,y)+\tfrac12\sigma^2\partial_{yy}u^*(t,y).
\]
The terminal HJB residual is zero by \eqref{eq:convex-HJB}.

\emph{Step 1: feasible comparisons with finite entropy.}
The interior first-order condition is
$\partial_af_1(t,y,a^*)+\partial_yu^*(t,y)=0$.
Taylor's formula and (C2) therefore give
\begin{equation}\label{eq:convex-slack}
 \frac{c_{f_1}}{2}|a-a^*(t,y)|^2
 \le\Gamma^*(t,y,a)
 \le\frac{C_{f_1}}{2}|a-a^*(t,y)|^2.
\end{equation}
For $0<\delta\le\delta_0$, define the probability kernel
\[
 \alpha^\delta_t(da\mid y)
   :=\frac{\mathbf1_{\{|a-a^*(t,y)|\le\delta\}}}{2\delta}\,da
\]
and the pair
$q^\delta=(\mu_T^*,(\mu_t^*(dy)\alpha^\delta_t(da\mid y))_t)$.
Its conditional mean action is $a^*(t,y)$.
Since the generator is affine in $a$ and its diffusion coefficient is
independent of $a$, testing the forward equation shows exactly that
\begin{equation}\label{eq:convex-comparison-feasible}
 \cA q^\delta=\cA\pi^*=\rho,\qquad
 \cM q^\delta=\bm\mu^*.
\end{equation}
The variance of $\alpha^\delta_t(\cdot\mid y)$ is $\delta^2/3$.
Consequently \eqref{eq:convex-slack} gives
\begin{equation}\label{eq:convex-comparison-cost}
 0\le\pair{\cC_{\bm\mu^*}-\cA^*u^*}{q^\delta}
 \le\frac{C_{f_1} T}{6}\delta^2.
\end{equation}
Using the density upper bound $M$ and the reference lower bound $c_0$,
with $0\log0=0$, we also have
\begin{align}
 \bD(q^\delta\|\pi_0)
 &\le \log(M/c_0)+T\log\!\left(\frac{M}{2c_0\delta}\right)
 \le C+T\log(1/\delta),\label{eq:convex-comparison-entropy}
\end{align}
for a constant $C$. Indeed, the terminal density ratio is at most
$M/c_0$, and the running density ratio, on its support, is at most
$M/(2c_0\delta)$; integrating their logarithms proves the bound.
In particular $q^\delta\in\Pizero$.

\emph{Step 2: control and forward-equation errors.}
By \eqref{eq:convex-slack} and \eqref{eq:convex-comparison-cost},
\[
 \pair{\cC_{\bm\mu^*}-\cA^*u^*}{\bar\pi_K-q^\delta}
 \ge\frac{c_{f_1}}{2}E_K-\frac{C_{f_1} T}{6}\delta^2.
\]
Moreover
$\pair{\cA^*u^*}{\bar\pi_K-q^\delta}=\pair{\cA\bar\pi_K-\rho}{u^*}$.
Insert $q^\delta$ into \eqref{eq:master}, using
$\cM q^\delta=\bm\mu^*$, to obtain
\[
 \frac{c_{f_1}}{2}E_K+\frac{\tau K}{2}\tns{\cA\bar\pi_K-\rho}^2
 \le \frac{\bD(q^\delta\|\pi_0)}{\lambda K}
       +\frac{C_{f_1} T}{6}\delta^2+\pair{\cA\bar\pi_K-\rho}{u_0-u^*}.
\]
Cauchy--Schwarz and Young's inequality imply
\begin{equation}\label{eq:convex-master}
 \frac{c_{f_1}}{2}E_K+\frac{\tau K}{4}\tns{\cA\bar\pi_K-\rho}^2
 \le \frac{C+T\log(1/\delta)}{\lambda K}
       +\frac{C_{f_1} T}{6}\delta^2
       +\frac{\tn{u_0-u^*}^2}{\tau K}.
\end{equation}
Take $\delta=\delta_0/\sqrt K$. This proves
\eqref{eq:convex-energy-rate}, including $K=1$.

\emph{Step 3: estimates for a backward linear equation.}
We establish the estimates needed to turn a forward-equation residual into a
Wasserstein error. For terminal data $\phi\in H^{m+1}(\mathbb T)$
and forcing $\psi\in L^2(0,T;H^m(\mathbb T))$, let $v$ solve
\begin{equation}\label{eq:convex-adjoint}
 -\partial_tv-\tfrac12\sigma^2\partial_{yy}v-a^*(t,y)\partial_yv=\psi,
 \qquad v(T,\cdot)=\phi.
\end{equation}
Then
\begin{equation}\label{eq:convex-adjoint-energy}
 \tn v\le C\bigl(\|\phi\|_{H^{m+1}}
                    +\|\psi\|_{L^2(0,T;H^m)}\bigr).
\end{equation}
Here is an energy proof, also giving existence in the required space.
Reverse time, writing $w(s)=v(T-s)$,
$b(s,y)=a^*(T-s,y)$ and $F(s)=\psi(T-s)$.
The equation is
$\partial_sw-\tfrac12\sigma^2\partial_{yy}w-b\partial_yw=F$.
For smooth data, differentiate $j$ times in space and pair with
$\partial_y^jw$, for $0\le j\le m+1$. The leading drift term integrates
to $-\frac12\int(\partial_yb)|\partial_y^jw|^2$.
The remaining drift terms are bounded by
$C\|b\|_{W^{m+1,\infty}}\|w\|_{H^{m+1}}^2$.
At the highest order the forcing is interpreted by
\[
 \pair{\partial_y^{m+1}F}{\partial_y^{m+1}w}
 =-\pair{\partial_y^mF}{\partial_y^{m+2}w}.
\]
Young's inequality absorbs this term into diffusion; the lower orders
are estimated directly. Summing the resulting inequalities gives
\[
 \frac{d}{ds}\|w(s)\|_{H^{m+1}}^2
 +\tfrac12\sigma^2\|\partial_yw(s)\|_{H^{m+1}}^2
 \le C\|w(s)\|_{H^{m+1}}^2+C\|F(s)\|_{H^m}^2.
\]
Gronwall's inequality controls $w$ in
$L^\infty(0,T;H^{m+1})\cap L^2(0,T;H^{m+2})$.
The equation then controls $\partial_sw$ in $L^2(0,T;H^m)$,
proving \eqref{eq:convex-adjoint-energy}.
Fourier--Galerkin approximation and weak compactness give existence;
the energy estimate for differences gives uniqueness for the stated
data. Spatial and temporal approximation of the coefficients
and data justifies the calculation under the stated regularity of $a^*$.

For spatially Lipschitz data there is also the bound
\begin{equation}\label{eq:convex-adjoint-gradient}
 \|\partial_yv\|_{L^2(0,T;L^\infty)}
 \le C\bigl(\operatorname{Lip}(\phi)
             +\|\operatorname{Lip}(\psi_t)\|_{L^2(0,T)}\bigr).
\end{equation}
To see this first for smooth data, differentiate the reversed equation.
The derivative $z=\partial_yw$ solves
$\partial_sz-\tfrac12\sigma^2\partial_{yy}z-b\partial_yz
  =(\partial_yb)z+\partial_yF$.
The maximum principle, with
$B=\|\partial_ya^*\|_{L^\infty([0,T]\times\mathbb T)}$, yields
\[
 \|\partial_yv(t)\|_\infty
 \le e^{B(T-t)}\operatorname{Lip}(\phi)
    +\int_t^T e^{B(s-t)}\operatorname{Lip}(\psi_s)\,ds.
\]
Cauchy--Schwarz in time proves \eqref{eq:convex-adjoint-gradient}.
Approximation extends it to Lipschitz terminal data and forcing
measurable in time with square-integrable Lipschitz seminorm.

\emph{Step 4: stability for an arbitrary probability pair.}
Let $\pi=(\nu,(\pi_t)_t)\in\Pi$, put $\mu_t=\pi_t(\cdot\times A)$,
and define
\[
 E(\pi)=\int_0^T\!\!\int|a-a^*(t,y)|^2\,d\pi_t\,dt,\qquad
 \eta(\pi)=\tns{\cA\pi-\rho},\qquad
 d(\pi)=\lVert\cM\pi-\bm\mu^*\rVert_{W_1,2} .
\]
We will prove
\begin{equation}\label{eq:convex-residual-stability}
 d(\pi)\le C\bigl(E(\pi)^{1/2}+\eta(\pi)^{1/(m+1)}\bigr).
\end{equation}
Fix a point $y_0\in\mathbb T$ and Lipschitz functions $\phi,\psi_t$
vanishing at $y_0$, with $t\mapsto\psi_t$ measurable and
\[
 \operatorname{Lip}(\phi)^2+
       \int_0^T\operatorname{Lip}(\psi_t)^2\,dt\le1.
\]
Their pairing with the population difference is
\[
 I(\phi,\psi)
 :=\int\phi\,d(\nu-\mu_T^*)
      +\int_0^T\!\!\int\psi_t\,d(\mu_t-\mu_t^*)\,dt.
\]
Kantorovich--Rubinstein duality and $L^2$ duality give
\begin{equation}\label{eq:convex-flow-duality}
 d(\pi)=\sup_{\phi,\psi} I(\phi,\psi)
\end{equation}
over these tests. For completeness, choose a countable uniformly dense
family of $1$-Lipschitz functions vanishing at $y_0$. Their pairings are
measurable in $t$, so one can choose measurably, within any prescribed
error, a maximizing function at each time. If
$w(t)=W_1(\mu_t,\mu_t^*)$ and $w_T=W_1(\nu,\mu_T^*)$, multiply these
functions by $w(t)/d(\pi)$ and $w_T/d(\pi)$, respectively, when
$d(\pi)>0$. Their squared seminorms sum to at most one, and their
pairings tend to $d(\pi)$. The reverse bound follows from
Cauchy--Schwarz, proving \eqref{eq:convex-flow-duality}.

Let $\phi_\epsilon,\psi_\epsilon$ be spatial convolutions with a smooth
periodic mollifier of radius $0<\epsilon\le\epsilon_0$, for a fixed
$\epsilon_0>0$ smaller than the injectivity radius of the torus.
Normalization at $y_0$ bounds the original sup norms by their Lipschitz
seminorms. The mollifier estimates are
\begin{align*}
 \|\phi_\epsilon-\phi\|_\infty
     &\le C\epsilon\operatorname{Lip}(\phi),&
 \|\phi_\epsilon\|_{H^{m+1}}
     &\le C\epsilon^{-m}\operatorname{Lip}(\phi),\\
 \|\psi_{\epsilon,t}-\psi_t\|_\infty
     &\le C\epsilon\operatorname{Lip}(\psi_t),&
 \|\psi_\epsilon\|_{L^2(0,T;H^m)}
     &\le C\epsilon^{-(m-1)}
                  \|\operatorname{Lip}(\psi_t)\|_{L^2(0,T)}.
\end{align*}
These follow by writing each spatial derivative of order $j\ge1$
as the convolution of the bounded first weak derivative with the
$(j-1)$st derivative of the mollifier.
Convolution does not increase the Lipschitz seminorms.
Let $v_\epsilon$ solve \eqref{eq:convex-adjoint} with these mollified
data. Step~3 gives
\begin{equation}\label{eq:convex-smoothed-tests}
 \tn{v_\epsilon}\le C\epsilon^{-m},\qquad
 \|\partial_yv_\epsilon\|_{L^2(0,T;L^\infty)}\le C.
\end{equation}
Subtracting the forward identities for $\pi$ and the feasible
reference $\pi^*$, and using the backward equation, gives the exact
identity
\[
 I(\phi_\epsilon,\psi_\epsilon)
 =\pair{\cA\pi-\rho}{v_\epsilon}
   +\int_0^T\!\!\int
       (a-a^*(t,y))\,\partial_yv_\epsilon(t,y)\,\pi_t(dy,da)\,dt.
\]
Since each $\pi_t$ is a probability measure, Cauchy--Schwarz and
\eqref{eq:convex-smoothed-tests} bound the right side in absolute
value by $C(\eta(\pi)\epsilon^{-m}+E(\pi)^{1/2})$.
The mollification error in $I$ is at most $C\epsilon$, because all
state measures have mass one and the Lipschitz seminorms have the
stated $L^2$ bound. Taking the supremum in
\eqref{eq:convex-flow-duality} proves
\[
 d(\pi)\le C\bigl(E(\pi)^{1/2}+\epsilon
                         +\eta(\pi)\epsilon^{-m}\bigr).
\]
For $0<\eta(\pi)\le\epsilon_0^{m+1}$ choose
$\epsilon=\eta(\pi)^{1/(m+1)}$; for $\eta(\pi)=0$ let
$\epsilon\downarrow0$. For larger residuals use
$d(\pi)\le\operatorname{diam}(\mathbb T)\sqrt{1+T}$ and increase $C$.
This proves \eqref{eq:convex-residual-stability}.

\emph{Step 5: rates for the state and action distributions.}
Apply \eqref{eq:convex-residual-stability} to $\bar\pi_K$ and insert
\eqref{eq:convex-energy-rate} to obtain \eqref{eq:convex-flow-rate}.
For the last assertion, map each $(y,a)$ to $(y,a^*(t,y))$.
This couples $\bar\pi_{K,t}$ to
$(\operatorname{Id},a^*(t,\cdot))_\#\bar\mu_{K,t}$ at cost at most
$\int|a-a^*(t,y)|\,d\bar\pi_{K,t}$.
The map $y\mapsto(y,a^*(t,y))$ is Lipschitz with constant
$1+\|\partial_ya^*\|_\infty$ for the product distance. Hence
\[
 W_{1,\mathbb T\times A}(\bar\pi_{K,t},\pi_t^*)
 \le \int|a-a^*(t,y)|\,d\bar\pi_{K,t}
       +(1+\|\partial_ya^*\|_\infty)
                        W_1(\bar\mu_{K,t},\mu_t^*).
\]
The $L^2$ norm of the first term is at most $\sqrt{E_K}$.
Minkowski's inequality, followed by
\eqref{eq:convex-energy-rate} and \eqref{eq:convex-flow-rate}, proves
the bound for \eqref{eq:convex-occupation-distance}.
\end{proof}

\subsection{Proof of Theorem~\ref{thm:finite-state}}\label{subsec:proof-finite}
\begin{proof}
(a) For $\nu\in\cP(E)$,
\[
    \KL(\nu\|\nu_0)
    =\sum_y\nu(y)\log\frac{\nu(y)}{\nu_0(y)}
    \le\log\frac1{\underline\nu},
\]
and likewise $\KL(\pi_t\|\pi_{0,t})\le\log(1/\underline\pi)$ for a.e.\ $t$. Hence
$\bD(\pi\|\pi_0)\le\log(1/\underline\nu)+T\log(1/\underline\pi)$ for every $\pi\in\Pi$, so
$\Pizero=\Pirho$. The set $\Pirho$ is nonempty. Let $\alpha$ be a feedback that is
Borel in $t$. The rates $Q^\alpha_t$ are Borel in $t$, with off-diagonal entries
in $[0,\bar q]$ and zero row sums. The forward equation $\dot\mu_t=\mu_tQ^\alpha_t$,
$\mu_0=\rho$, is linear with bounded measurable coefficients, so by the
Carath\'eodory existence and uniqueness theorem it has a unique absolutely
continuous solution on $[0,T]$. It lies in $\cP(E)$: the zero row sums keep
$\sum_y\mu_t(y)=1$, and $B_t:=Q^\alpha_t+\bar qI$ has nonnegative entries, so
$e^{\bar qt}\mu_t=\rho+\int_0^te^{\bar qs}\mu_sB_s\,ds$ is the limit of Picard
iterates started at $\rho$, all of which are nonnegative. The product rule for
absolutely continuous functions gives, for $u\in\HY$ and a.e. $t$,
\[
    \frac{d}{dt}\pair{u(t)}{\mu_t}
    =\pair{\partial_t u(t)}{\mu_t}
     +\pair{u(t)}{\mu_tQ^\alpha_t}.
\]
After integration this yields $\cA \pi=\rho$ for
$\pi=(\mu_T,(\mu_t\otimes\alpha_t)_t)$.

(b) The right side of \eqref{eq:hjb-finite} is Borel in $t$, continuous
in $V$, bounded by $\lVert f\rVert_\infty+2\bar q|V|_\infty$ and Lipschitz
in $V$ with constant $2\bar q$, uniformly in $t$, as a minimum of affine
maps with uniformly bounded coefficients. By the Carath\'eodory existence
and uniqueness theorem the terminal value problem has a unique absolutely
continuous solution on $[0,T]$, with $\dot V\in L^\infty$. Hence
$V\in W^{1,\infty}(0,T;\R^E)\subset H^1(0,T;H)=\HY$. Its HJB residuals at
$z^*=(\pi^*,V)$ are $G^{z^*}=g(\cdot,\mu^*_T)-V(T)=0$ and
\[
\begin{aligned}
    \Gamma^{z^*}(t,y,a)&=f(t,y,a,\mu^*_t)+\dot V(t,y)+\LL^a_tV(t,\cdot)(y)\\
    &=f(t,y,a,\mu^*_t)+\LL^a_tV(t,\cdot)(y)-\min_{a'\in A}\bigl[f(t,y,a',\mu^*_t)+\LL^{a'}_tV(t,\cdot)(y)\bigr]\ge0 .
\end{aligned}
\]
For the duality gap, let $\pi\in\Pirho$. Since $\pair{\cA^*V}{\pi}=\pair{\cA \pi}{V}=\pair{\rho}{V}$,
\[
    \pair{\cC_{\bm\mu^*}}{\pi}-\pair{\rho}{V}=\pair{\cC_{\bm\mu^*}-\cA^*V}{\pi}=\pair{G^{z^*}}{\nu}+\int_0^T\!\!\int\Gamma^{z^*}_t\,d\pi_t\,dt\ge0 .
\]
Let $\alpha^*_t(y)$ be a Borel selection of the argmin in
\eqref{eq:hjb-finite}, which exists since $A$ is compact and the bracket
is continuous in $a$ and Borel in $t$, and let $\pi^{\alpha^*}\in\Pirho$ be the
occupation measure of the chain with rates $Q^{\alpha^*}$, constructed as
in (a). Then $\Gamma^{z^*}(t,y,\alpha^*_t(y))=0$, so the display is an
equality at $\pi=\pi^{\alpha^*}$: $\inf_{\pi\in\Pirho}\pair{\cC_{\bm\mu^*}}{\pi}=\pair{\rho}{V}$.
Since $\pi^*$ is an MFE it attains this infimum, so
$\pair{\cC_{\bm\mu^*}}{\pi^*}=\pair{\rho}{V}$, and $V$ is an optimal dual variable.
\Cref{prop:PDVI-MFE}(ii) gives the rest.

(c) Write $e_t:=\mu_t-\tilde\mu_t$. For $\phi\in\R^E$,
\[
    \pair{\phi}{\LL^*_t\pi_t}=\sum_{y,a}\pi_t(y,a)(Q^a_t\phi)(y)=\sum_y\mu_t(y)\sum_a\alpha_t(a\mid y)(Q^a_t\phi)(y)=\pair{\phi}{\mu_tQ^\alpha_t},
\]
so $\LL^*_t\pi_t=\mu_tQ^\alpha_t$ and $d_s=\mu_s-\rho-\int_0^s\mu_tQ^\alpha_tdt$.
Subtracting $\tilde\mu_s=\rho+\int_0^s\tilde\mu_tQ^\alpha_tdt$ gives
$e_s=d_s+\int_0^se_tQ^\alpha_tdt$. Set $b_s:=\int_0^se_tQ^\alpha_tdt$. Then
$b$ is absolutely continuous, $b_0=0$ and $\dot b_s=b_sQ^\alpha_s+d_sQ^\alpha_s$
for a.e.\ $s$. Let $P_{t,s}$ be the transition matrix of the chain from time $t$ to $s$,
$\partial_sP_{t,s}=P_{t,s}Q^\alpha_s$, $P_{t,t}=I$, which exists by
Carath\'eodory's theorem since $Q^\alpha$ is bounded and Borel in $t$; it
is a stochastic matrix, so $\lVert \xi P_{t,s}\rVert_1\le\lVert \xi\rVert_1$ for
every row vector $\xi$. Variation of constants gives
$b_s=\int_0^sd_tQ^\alpha_tP_{t,s}\,dt$, and since
$\lVert \xi Q^\alpha_t\rVert_1\le2\bar q\lVert \xi\rVert_1$, the bounded
measurable function $a_s:=\lVert e_s\rVert_1$ satisfies
$a_s\le\lVert d_s\rVert_1+2\bar q\int_0^s\lVert d_t\rVert_1dt$.
By Young's convolution inequality with the kernel $2\bar q\mathbf 1_{[0,T]}$,
whose $L^1$ norm is $2\bar qT$, and with
$\delta:=(\int_0^T\lVert d_t\rVert_1^2dt)^{1/2}$,
\[
    \lVert a\rVert_{L^2(0,T)}\le(1+2\bar qT)\,\delta .
\]
Now $\lVert \pi_t-\tilde\mu_t\otimes\alpha_t\rVert_1=\lVert e_t\otimes\alpha_t\rVert_1=a_t$,
because $\pi_t=\mu_t\otimes\alpha_t$ by construction, and
$\nu-\tilde\mu_T=d_T+\int_0^Te_tQ^\alpha_tdt$ gives
$\lVert\nu-\tilde\mu_T\rVert_1\le\delta_T+2\bar q\sqrt T\lVert a\rVert_{L^2}$
with $\delta_T:=\lVert d_T\rVert_1$. Hence, with $c:=2\bar q\sqrt T(1+2\bar qT)$ and
Cauchy--Schwarz in the form $(s+ct)^2\le(1+c^2)(s^2+t^2)$,
\[
    \nsharp{\projr(\pi)-\pi}^2\le(\delta_T+c\delta)^2+(1+2\bar qT)^2\delta^2
    \le\bigl(1+c^2+(1+2\bar qT)^2\bigr)(\delta_T^2+\delta^2) .
\]
Finally $\lVert \xi\rVert_1\le\sqrt N|\xi|_2$ on $\R^E$ and
$\tns{\cA \pi-\rho}^2=|d_T|_2^2+\int_0^T|d_t|_2^2dt$ by \Cref{lem:riesz}, so
$\delta_T^2+\delta^2\le N\tns{\cA \pi-\rho}^2$, which is the claim.
Feasibility of $\projr(\pi)$ follows as in (a), and measurability of
$t\mapsto\alpha_t$ from that of $t\mapsto \pi_t$.
\end{proof}

\subsection{Proof of Theorem~\ref{thm:mesh}}\label{subsec:proof-mesh}
\begin{proof}
\emph{Discrete embedding.} For $i,j\in\mathbb Z/N\mathbb Z$ and the
increasing path from $j$ to $i$,
\[
    \phi_i^2-\phi_j^2=\sum_{k=j}^{i-1}(\phi_{k+1}+\phi_k)(\phi_{k+1}-\phi_k)
    =h\sum_{k=j}^{i-1}(\phi_{k+1}+\phi_k)\,D\phi_k
    \le\lVert \phi_{\cdot+1}+\phi\rVert_h\lVert D\phi\rVert_h\le2\lVert \phi\rVert_h\lVert D\phi\rVert_h
\]
by Cauchy--Schwarz over all $N$ indices. Multiplying by $h$ and summing over
$j$ gives $2L\phi_i^2\le\lVert \phi\rVert_h^2+4L\lVert \phi\rVert_h\lVert D\phi\rVert_h$, so
\begin{equation}\label{eq:embed}
    \lVert \phi\rVert_\infty^2\le\frac1{2L}\lVert \phi\rVert_h^2+2\lVert \phi\rVert_h\lVert D\phi\rVert_h
    \le c_H^2\bigl(\lVert \phi\rVert_h^2+\lVert D\phi\rVert_h^2\bigr).
\end{equation}
Applied to the periodic grid function $D^j\phi$, and since every binomial
coefficient in \eqref{eq:Hm} is at least $1$,
\begin{equation}\label{eq:embed-j}
    \lVert D^j\phi\rVert_\infty\le c_H\bigl(\lVert D^j\phi\rVert_h^2+\lVert D^{j+1}\phi\rVert_h^2\bigr)^{1/2}\le c_H\nH \phi
    \qquad(0\le j\le m-1).
\end{equation}

(i) (S2) with this $c_H$ is \eqref{eq:embed-j} with $j=0$. In \eqref{eq:gen-h}, $D^-\phi$ and $\Delta_h\phi$ are index shifts of $D\phi$ and $D^2\phi$ and have the same sup norms, so
\[
    \lVert\LL^a_{h,t}\phi\rVert_\infty
    \le\lVert b\rVert_\infty\lVert D\phi\rVert_\infty+\tfrac12\lVert\sigma^2\rVert_\infty\lVert D^2\phi\rVert_\infty
    \le c_H\bigl(\lVert b\rVert_\infty+\tfrac12\lVert\sigma^2\rVert_\infty\bigr)\nH \phi
\]
by \eqref{eq:embed-j} with $j=1,2$, which requires $m\ge3$. Measurability in
$t$ and continuity in $(y,a)$ are inherited from $b,\sigma$ on the finite
grid, and the formula for $L_A$ is \Cref{lem:elementary}.

(ii) For finite $A$, the entropy estimate of \Cref{thm:finite-state}(a)
uses the lower bounds on the reference masses and is independent of the
inner product on $\R^{E_h}$. Under the bounds with a common $c$,
$\log(1/\underline\nu)\le\log(N/c)$ and $\log(1/\underline\pi)\le\log(N|A|/c)$.

(iii) For compact $A$, the HJB and verification argument of
\Cref{thm:finite-state}(b) applies: the minimum is attained by a measurable
selector since the minimized expression is Borel in time and continuous
in the action. The resulting $V\in W^{1,\infty}(0,T;\R^{E_h})$ belongs
to $H^1(0,T;H_h^m)$ by equivalence of norms on $\R^{E_h}$.

(iv) From $|\phi_{i+1}-\phi_i|^2\le2(\phi_{i+1}^2+\phi_i^2)$ we get
$\lVert D\phi\rVert_h\le\frac2h\lVert \phi\rVert_h$, hence
$\nH \phi^2\le\sum_j\binom mj(4/h^2)^j\lVert \phi\rVert_h^2=(1+4/h^2)^m\lVert \phi\rVert_h^2$.
For a measure $\xi$ on $E_h$,
\begin{align*}
    \lVert \xi\rVert_{H^*}=\sup_{\nH \phi\le1}\pair{\phi}{\xi}
    &\ge\sup_{\lVert \phi\rVert_h\le(1+4/h^2)^{-m/2}}(\phi,\xi/h)_h\\
    &=\Bigl(1+\frac4{h^2}\Bigr)^{-m/2}\lVert \xi/h\rVert_h\\
    &=\Bigl(1+\frac4{h^2}\Bigr)^{-m/2}h^{-1/2}|\xi|_2.
\end{align*}
and $\lVert \xi\rVert_1\le\sqrt N|\xi|_2=\sqrt{2L/h}\,|\xi|_2$, so
\begin{equation}\label{eq:inverse}
\begin{aligned}
    \lVert \xi\rVert_1
    &\le\sqrt{2L}\,\Bigl(1+\frac4{h^2}\Bigr)^{m/2}
      \lVert \xi\rVert_{H^*}.
\end{aligned}
\end{equation}
The proof of \Cref{thm:finite-state}(c) also applies to compact $A$:
integrating the rates against the conditional action law gives the same
matrix $Q^\alpha_t$ and the same estimate on the state marginals.
Thus, for the total-variation
quantities $\delta_T$ and $\delta$ defined there,
$\nsharp{\projr(\pi)-\pi}^2\le\bigl(1+(1+4\bar q^2T)(1+2\bar qT)^2\bigr)(\delta_T^2+\delta^2)$, and by
\eqref{eq:inverse} and \Cref{lem:riesz},
$\delta_T^2+\delta^2\le2L(1+4/h^2)^m\tns{\cA \pi-\rho}^2$, which is
\eqref{eq:CPi-h}. The bound $O(h^{-m-4})$ follows from $\bar q=O(h^{-2})$.
\end{proof}

\subsection{Proof of Proposition~\ref{prop:mesh-limit}}
\label{subsec:proof-mesh-limit}

Throughout this subsection, $\Pi_h(\rho_h)$ denotes the feasible set
for the chain \eqref{eq:rates}. We identify its occupation measures
with measures on the torus. All diffusion flows below are taken in
their narrowly continuous versions, including their initial and
terminal values.

\begin{lemma}
\label{lem:mesh-recovery}
Under the coefficient assumptions of \Cref{prop:mesh-limit}, for every
$\pi=(\nu,(\pi_t)_t)\in\Pi(\rho)$ there are
$\pi_N\in\Pi_h(\rho_h)$ such that
\begin{equation}\label{eq:mesh-recovery}
    \pi_N\Rightarrow\pi,
    \qquad
    \sup_{0\le t\le T}
       W_1\bigl((\cM\pi_N)_t,(\cM\pi)_t\bigr)\longrightarrow0.
\end{equation}
If $K_N$ satisfies $\log N/K_N\to0$ and the references satisfy
\Cref{prop:mesh-limit}, the approximating sequence can moreover be chosen
so that
\begin{equation}\label{eq:mesh-recovery-entropy}
    \frac{\bD(\pi_N\|\pi_0^N)}{K_N}\longrightarrow0.
\end{equation}
\end{lemma}

\begin{proof}
Write $\mu_t=(\cM\pi)_t$ and disintegrate
$\pi_t(dy,da)=\mu_t(dy)\alpha_t(da\mid y)$. Choose a Borel version of this probability kernel and set
\[
    B(t,y):=\int_A b(t,y,a)\,\alpha_t(da\mid y).
\]
Testing the forward equation with $v(t,y)=\chi(t)\varphi(y)$ shows
that $t\mapsto\int\varphi\,d\mu_t$ has an absolutely continuous
version for each smooth periodic $\varphi$, with derivative
$\int\LL_t^a\varphi\,d\pi_t$. A countable uniformly dense family
of smooth functions therefore determines a narrowly continuous
marginal representative; the boundary terms give $\mu_0=\rho$ and
$\mu_T=\nu$. The superposition principle realizes
this flow as the marginals of a weak solution of
\begin{equation}\label{eq:mesh-recovery-sde}
    dX_t=B(t,X_t)\,dt+\sigma(t,X_t)\,dW_t.
\end{equation}
For completeness, the integrability hypothesis of
\cite[Theorem~2.5]{Trevisan2016} holds because both coefficients are
bounded. Its Euclidean statement applies on the torus by a smooth
embedding of the circle in $\R^2$: the pushed-forward coefficients
are bounded, and the superposition paths remain on the circle since
all their one-time marginals are supported there.

We will also use weak uniqueness for \eqref{eq:mesh-recovery-sde},
with any bounded Borel drift in place of $B$. The zero-drift equation
has a unique strong solution because $\sigma$ is uniformly Lipschitz
in space. The ratio of any bounded drift to $\sigma$ is bounded.
Girsanov's theorem, applied in both directions, therefore gives
existence and uniqueness in law with that drift. One can perform the
argument on the real line with periodic coefficients and then project
onto the torus.

\emph{Step 1: approximation of the relaxed feedback.}
Let $\eta(dt,dy)=dt\,\mu_t(dy)$. We construct kernels
$\alpha^j_t(da\mid y)=w_j(t,y,a)\zeta(da)$, weakly continuous in
$(t,y)$, such that $0<\varepsilon_j\le w_j\le M_j<\infty$ and
\begin{equation}\label{eq:mesh-feedback-approximation}
 \int|B^j-B|^2\,d\eta\longrightarrow0,
 \qquad B^j(t,y):=\int_A b(t,y,a)\alpha^j_t(da\mid y).
\end{equation}
They also satisfy, for every bounded Borel function $F(t,y,a)$
continuous in $(y,a)$ for each $t$,
\begin{equation}\label{eq:mesh-feedback-test-approximation}
 \int\left|\int_A F(t,y,a)
           [\alpha^j_t-\alpha_t](da\mid y)\right|\eta(dt,dy)
       \longrightarrow0.
\end{equation}

For an interval $A$, partition it into $n_j$ intervals
$I_{j,l}$ of equal length, where $n_j\to\infty$, assigning endpoints
consistently. Put
$p_{j,l}(t,y)=\alpha_t(I_{j,l}\mid y)$ and
$\kappa_{j,l}:=n_j\mathbf1_{I_{j,l}}\zeta$.
Approximate the simplex-valued vector $p_j$ by a jointly Lipschitz
simplex-valued vector $v_j$ so that
\[
 \int\left(\sum_{l=1}^{n_j}|v_{j,l}-p_{j,l}|\right)^2d\eta
       \longrightarrow0.
\]
This approximation follows from density of continuous functions in
$L^2(\eta)$, uniform smooth approximation on the compact
time--torus cylinder, and Lipschitz Euclidean projection onto the
simplex. The approximation accuracy can be chosen separately for
each finite dimension $n_j$.
For $\varepsilon_j\in(0,1)$ with $\varepsilon_j\downarrow0$, define
\[
 \alpha^j_t(da\mid y)
   =(1-\varepsilon_j)\sum_l v_{j,l}(t,y)\kappa_{j,l}(da)
        +\varepsilon_j\zeta(da).
\]
Then $\varepsilon_j\le w_j\le n_j+1=:M_j$, and the kernel is
Lipschitz in $(t,y)$ in total variation. Replacing a point in a bin
by its uniform law changes the integral of $F(t,y,\cdot)$ by at
most its action modulus of continuity at the bin diameter.
For each $t$, this modulus tends to zero uniformly in $y$, by
compactness and joint continuity. The modulus is measurable in $t$
and bounded, so dominated convergence, the approximation of the
weights, and $\varepsilon_j\to0$ prove
\eqref{eq:mesh-feedback-test-approximation}. Applying the same
argument to $F=b$, with squares and boundedness, proves
\eqref{eq:mesh-feedback-approximation}.
If $A$ is finite, use its singletons as the fixed partition,
$\kappa_{j,l}=\delta_{a_l}$, and approximate the corresponding
probability vector in the same way. Then one can take
$M_j=|A|+1$, and both convergence statements follow directly.

Let $P$ and $P^j$ be the path laws with drifts $B$ and $B^j$,
respectively, the same diffusion coefficient and initial law $\rho$.
Girsanov's formula and contraction of relative entropy under the
state-path map give
\[
 D_{\rm KL}(P\|P^j)
 \le\frac12\mathbb E_P\int_0^T
       \left|\frac{B^j-B}{\sigma}(t,X_t)\right|^2dt
 \le C\int|B^j-B|^2\,d\eta\longrightarrow0.
\]
The stochastic exponentials are true martingales because their
integrands are bounded. Pinsker's inequality gives
$\|P^j-P\|_1\to0$. Write
$\mu^j_t=P^j\circ X_t^{-1}$ and
$\pi^j_t(dy,da)=\mu^j_t(dy)\alpha^j_t(da\mid y)$. Then
\[
 \sup_t\|\mu^j_t-\mu_t\|_1\le\|P^j-P\|_1\longrightarrow0.
\]
For a bounded continuous test $F$ on $[0,T]\times E\times A$,
\[
\begin{split}
 \left|\int_0^T\!\int F\,d(\pi^j_t-\pi_t)\,dt\right|
 &\le T\|F\|_\infty\|P^j-P\|_1\\
 &\quad+\int\left|\int_A F(t,y,a)
         [\alpha^j_t-\alpha_t](da\mid y)\right|d\eta
 \longrightarrow0.
\end{split}
\]
Thus $\pi^j=(\mu^j_T,(\pi^j_t)_t)$ converges weakly to $\pi$, and
its marginal flow converges uniformly in time in $W_1$.

\emph{Step 2: approximation for a fixed feedback.}
Fix $j$. Let $X^{N,j}$ be the continuous-time chain starting from
$\rho_h$ and using the rates
\[
    Q_t^{N,j}(y_i,y_{i\pm1})
       =\int_A Q_t(y_i,y_{i\pm1},a)\,\alpha^j_t(da\mid y_i).
\]
These are bounded Borel functions of time for each fixed $N$, so the
chain and its forward equation are well defined. Let $\mu^{N,j}_t$
be its law and put
\[
    \pi^{N,j}_t(dy,da)
       =\mu^{N,j}_t(dy)\alpha^j_t(da\mid y),
    \qquad
    \pi^{N,j}=(\mu^{N,j}_T,(\pi^{N,j}_t)_t).
\]
Then $\pi^{N,j}\in\Pi_h(\rho_h)$ exactly.

Here is the weak-convergence argument for these chains, including
the role of their measurable time dependence. A real-valued lift
of $X^{N,j}$ has drift $B^j(t,X^{N,j}_t)$, jumps of size $h$, and
martingale predictable quadratic variation with density
\[
    \sigma^2(t,X^{N,j}_t)
       +h\int_A |b(t,X^{N,j}_t,a)|\,\alpha^j_t(da\mid X^{N,j}_t).
\]
The drift and this density are bounded uniformly in $N$. The
martingale second-moment bound, also between bounded stopping times,
gives tightness on $D([0,T];E)$; the vanishing jump sizes imply that
all limit laws have continuous paths. For a smooth periodic
$\varphi$, Taylor expansion gives
\begin{equation}\label{eq:mesh-fixed-feedback-consistency}
    \sup_{t,i}\left|
       \int_A\LL^a_{h,t}\varphi(y_i)\,\alpha^j_t(da\mid y_i)
       -B^j(t,y_i)\varphi'(y_i)
       -\tfrac12\sigma^2(t,y_i)\varphi''(y_i)
       \right|\le C_\varphi h.
\end{equation}
The integrand in the limiting generator is bounded, Borel in time
and continuous in the state for each time. If paths converge in the
Skorokhod topology to a continuous path, they converge uniformly;
pointwise-in-time continuity in the state and dominated convergence
therefore give convergence of the time-integrated generators.
Moreover $W_1(\rho_h,\rho)\le h/2$, so every limit has initial
law $\rho$. Passing to the chain martingale identities, multiplied by bounded
continuous functions of finitely many earlier coordinates, identifies
every limit as a solution with drift $B^j$ and diffusion $\sigma$.
The martingales are uniformly integrable because the test functions
and their generators are uniformly bounded. Weak uniqueness proved
above identifies this limit as $P^j$, and hence the entire sequence
of chain path laws converges to $P^j$. This is the continuous-time
martingale-problem argument underlying Markov-chain approximation
\cite{KushnerDupuis2001}, with Borel time dependence.

A Skorokhod representation, the continuity of the limit paths and
boundedness of the torus distance now imply
\begin{equation}\label{eq:mesh-fixed-feedback-flow}
    \sup_{0\le t\le T}
       W_1(\mu^{N,j}_t,\mu^j_t)\longrightarrow0.
\end{equation}
Indeed, this supremum is bounded by the expectation of the uniform
distance of the coupled paths. Integrating continuous functions of
$(t,y,a)$ against $\alpha^j_t(da\mid y)$ in the same coupling gives
$\pi^{N,j}\Rightarrow\pi^j$.

\emph{Step 3: the diagonal sequence and entropy.}
For a fixed $j$, the running comparison has density
$\mu^{N,j}_t(y_i)w_j(t,y_i,a)$ relative to counting measure on the
grid times $\zeta$. The reference dominates $(c/N)\zeta$ at each
state, so
\[
\begin{split}
 D_{\rm KL}(\pi^{N,j}_t\|\pi^N_{0,t})
 &\le \log(N/c)+\sum_i\mu^{N,j}_t(y_i)\log\mu^{N,j}_t(y_i)\\
 &\quad+\sum_i\mu^{N,j}_t(y_i)
             \int_Aw_j(t,y_i,a)\log w_j(t,y_i,a)\,\zeta(da)\\
 &\le\log(N/c)+\log M_j.
\end{split}
\]
The terminal entropy is at most $\log(N/c)$. Therefore
\begin{equation}\label{eq:mesh-recovery-fixed-entropy}
 \bD(\pi^{N,j}\|\pi_0^N)
 \le (1+T)\log(N/c)+T\log M_j.
\end{equation}
This estimate also holds if the reference has a singular part,
since only its stated lower bound as a measure is used.

Fix a metric for $\Rightarrow$ as in \Cref{lem:topology}.
Choose integers $N_j$ increasing so that, for all $N\ge N_j$,
both the distance between $\pi^{N,j}$ and $\pi^j$ in this metric and
the supremum in \eqref{eq:mesh-fixed-feedback-flow} are at most
$1/j$. If an iteration sequence $K_N$ is prescribed as in the
lemma, increase $N_j$ further so that
$K_N\ge j(1+T\log M_j)$ for all $N\ge N_j$; this is possible
because $K_N\to\infty$. Let $j(N)\to\infty$ sufficiently slowly
that $N\ge N_{j(N)}$, and put $\pi_N=\pi^{N,j(N)}$.
Steps~1--2 give \eqref{eq:mesh-recovery}. If $K_N$ is prescribed,
\eqref{eq:mesh-recovery-fixed-entropy} gives
\[
 \frac{\bD(\pi_N\|\pi_0^N)}{K_N}
 \le (1+T)\frac{\log(N/c)}{K_N}+\frac1{j(N)}
 \longrightarrow0,
\]
which proves \eqref{eq:mesh-recovery-entropy}.
\end{proof}

\begin{proof}[Proof of \Cref{prop:mesh-limit}]
Write $\bar\pi_N=\bar\pi_{N,K_N}$ and
$r_N=\cA_h\bar\pi_N-\rho_h$, with its dual norm taken in
$Y_h^*=(H^1(0,T;H_h^m))^*$. The comparison obtained by using the fixed relaxed policy $\zeta$
in the chain has entropy at most
$D_N^0:=(1+T)\log(N/c)$. Apply \eqref{eq:master} to this comparison,
use $u_0^N=0$ and
$|\pair{\cC_{\cM\pi'}}{\bar\pi_N-\pi'}|\le 2C_0\sqrt{1+T}$, and obtain
\begin{equation}\label{eq:mesh-limit-residual}
    \|r_N\|_{Y_h^*}^2
       \le \frac{2D_N^0}{\lambda\tau K_N^2}
             +\frac{4C_0\sqrt{1+T}}{\tau K_N}
       \longrightarrow0.
\end{equation}
For every feasible comparison $\pi'\in\Pi_h(\rho_h)$ with finite
entropy, \eqref{eq:master} also gives
\begin{equation}\label{eq:mesh-limit-minty}
    \pair{\cC_{\cM\pi'}}{\bar\pi_N-\pi'}
       \le\frac{\bD(\pi'\|\pi_0^N)}{\lambda K_N}.
\end{equation}
All constants in these estimates are independent of the grid.

\emph{Step 1: compactness and feasibility.}
Compactness follows from \Cref{lem:topology}, since all grid pairs
belong to the same compact space $\Pi$ of probability pairs. Pass to a
subsequence with $\bar\pi_N\Rightarrow\bar\pi$. Let
$v\in C^1([0,T];C^\infty(E))$, periodic in space, and let $v_h$ be
its restriction to the grid. Finite differences of a sampled smooth
function are bounded by the corresponding derivatives, so
\[
    \sup_N\|v_h\|_{Y_h}<\infty.
\]
Taylor expansion in space, boundedness of $b,\sigma$ and
\eqref{eq:mesh-limit-residual} show that
\[
\begin{split}
    &\int_E v(T,y)\,\bar\nu_N(dy)
      -\int_E v(0,y)\,\rho_h(dy)\\
    &\qquad
      -\int_0^T\!\int_{E\times A}
         (\partial_t v+\LL_t^a v)(t,y)
           \,\bar\pi_{N,t}(dy,da)\,dt
       \longrightarrow0.
\end{split}
\]
The running integrand is bounded, Borel in time and continuous in
$(y,a)$ for each time. Thus \Cref{lem:topology} passes this identity
to the limit for the stated Borel time dependence. Also $\rho_h\Rightarrow\rho$. We conclude that
\[
    \int_Ev(T,y)\,\bar\nu(dy)
       -\int_Ev(0,y)\,\rho(dy)
    =\int_0^T\!\int_{E\times A}
       (\partial_t v+\LL_t^a v)(t,y)
         \,\bar\pi_t(dy,da)\,dt.
\]
Smooth space--time functions are dense in
$H^1(0,T;H^m(E))$, and both sides are continuous in this norm.
Consequently $\bar\pi\in\Pi(\rho)$. This forward equation gives
the narrowly continuous marginal representative with initial law
$\rho$ and terminal law $\bar\nu$, as in the proof of
\Cref{lem:mesh-recovery}.

\emph{Step 2: compactness of the population flows.}
We prove the stronger convergence in \eqref{eq:mesh-population-convergence}.
Put $s=m+1$ and use the usual negative Sobolev spaces on the torus.
For $\varphi\in H^s(E)$, restriction to the grid satisfies
$\|\varphi|_{E_h}\|_{H_h^m}\le C\|\varphi\|_{H^s}$:
each finite difference of order $j\le m$ is bounded by
$\|\partial_y^j\varphi\|_\infty$, which is controlled by $H^s$.
Also $H^s(E)\subset C^3(E)$, and Taylor's formula gives
\[
 \sup_{t,i,a}
 |\LL^a_{h,t}(\varphi|_{E_h})(y_i)-\LL^a_t\varphi(y_i)|
 \le Ch\|\varphi\|_{H^s}.
\]
Let $Z=H^1(0,T;H^s(E))$, with the norm defined for the dual
space in Section~\ref{sec:setting}. Define $r_N^c\in Z^*$ as the forward-equation residual of $\bar\pi_N$,
viewed as measures on the torus, with the diffusion generator and
initial law $\rho_h$. The preceding estimates imply
\begin{equation}\label{eq:mesh-continuum-residual}
 \|r_N^c\|_{Z^*}\le C\bigl(\|r_N\|_{Y_h^*}+h\bigr)\longrightarrow0.
\end{equation}
Indeed, the restriction estimate controls the discrete residual
tested on $v|_{E_h}$, and generator consistency controls the
remaining time integral by $Ch\|v\|_Z$.

Define the distribution-valued curve
\[
 B_N(t):=\rho_h+\int_0^t\LL_{t'}^*\bar\pi_{N,t'}\,dt',
 \qquad
 \pair{\LL_{t'}^*\bar\pi_{N,t'}}{\varphi}
 :=\int_{E\times A}\LL_{t'}^a\varphi(y)\,
                       \bar\pi_{N,t'}(dy,da).
\]
Bounded coefficients and $H^m(E)\subset C^2(E)$ imply that $B_N$
is uniformly bounded and Lipschitz as an $H^{-m}(E)$-valued curve.
The compact embedding $H^{-m}(E)\subset H^{-s}(E)$ and
Arzel\`a--Ascoli therefore give relative compactness in
$C([0,T];H^{-s}(E))$.
Writing $\bar\mu_{N,t}=\bar\pi_{N,t}(\cdot\times A)$, integration
in time gives
\[
 \pair{r_N^c}{v}
 =\pair{\bar\nu_N-B_N(T)}{v(T)}
   -\int_0^T\pair{\bar\mu_{N,t}-B_N(t)}{\partial_tv(t)}\,dt.
\]
The map $v\mapsto(v(T),\partial_tv)$ is an onto isometry from
$Z$ to $H^s\times L^2(0,T;H^s)$. Consequently
\begin{equation}\label{eq:mesh-population-sobolev}
 \|\bar\nu_N-B_N(T)\|_{H^{-s}}^2
 +\int_0^T\|\bar\mu_{N,t}-B_N(t)\|_{H^{-s}}^2\,dt
 =\|r_N^c\|_{Z^*}^2\longrightarrow0.
\end{equation}
Along any further subsequence for which $B_N$ converges uniformly
in $H^{-s}$, the weak occupation limit identifies its limit with
$\bar\mu_t=(\cM\bar\pi)_t$ for a.e.\ $t$. The extracted limit is continuous
in $H^{-s}$ by Arzel\`a--Ascoli, and $\bar\mu$ is continuous by the
forward equation, so they agree at every time,
including $T$. Every such subsequence has the same limit. Thus
$\bar\mu_N\to\bar\mu$ in $L^2(0,T;H^{-s})$ and
$\bar\nu_N\to\bar\nu$ in $H^{-s}$ along the original
subsequence for which the occupation measures converge.

Let $\mu,\mu'$ be probability measures on the torus. Use the
Kantorovich--Rubinstein formula with $1$-Lipschitz tests normalized
to vanish at a fixed point. Spatial mollification gives
\[
 W_1(\mu,\mu')\le
 C\left(\epsilon+\epsilon^{-(s-1)}
                   \|\mu-\mu'\|_{H^{-s}}\right),
 \qquad 0<\epsilon\le\epsilon_0.
\]
Optimizing in $\epsilon$, and using boundedness of $W_1$ for the
remaining range, yields
$W_1(\mu,\mu')\le C\|\mu-\mu'\|_{H^{-s}}^{1/s}$.
H\"older's inequality in time now proves
\eqref{eq:mesh-population-convergence} with $\bar\pi$ in place of
$\pi^*$. In particular the population flows are relatively compact
in the metric \eqref{eq:flat-norm}.

\emph{Step 3: passage to the limit in the Minty inequality.}
Fix an arbitrary $\pi'\in\Pi(\rho)$ and use \Cref{lem:mesh-recovery}
to choose feasible chain occupation measures $\pi'_N$ approximating
$\pi'$, with
$\bD(\pi'_N\|\pi_0^N)/K_N\to0$. Write
$\bm\mu'_N=\cM\pi'_N$ and $\bm\mu'=\cM\pi'$. The uniform
convergence of these flows in $W_1$ and joint continuity of the
costs imply, for each $t$,
\[
    \sup_{y,a}
       |f(t,y,a,\mu'_{N,t})-f(t,y,a,\mu'_t)|\longrightarrow0,
    \qquad
    \sup_y|g(y,\mu'_{N,T})-g(y,\mu'_T)|\longrightarrow0.
\]
The costs are bounded, so dominated convergence yields
\begin{equation}\label{eq:mesh-recovery-costs}
    \nG{\cC_{\bm\mu'_N}-\cC_{\bm\mu'}}\longrightarrow0.
\end{equation}
Regard these cost vectors as functions on the full torus before
restricting them to grid states. Since
$\nsharp{\bar\pi_N-\pi'_N}\le2\sqrt{1+T}$, replacing the cost vector
in \eqref{eq:mesh-limit-minty} by $\cC_{\bm\mu'}$ changes its
left side by a quantity tending to zero. For this fixed cost vector,
\Cref{lem:topology} and the weak convergence of both occupation
pairs give
\[
    \pair{\cC_{\cM\pi'}}{\bar\pi-\pi'}\le0
    \qquad\text{for every }\pi'\in\Pi(\rho).
\]
Thus $\GapM(\bar\pi)=0$, and $\bar\pi$ is an MFE by \Cref{prop:eps-relations}(ii). The superposition and uniqueness arguments above realize it as a relaxed-control equilibrium of the diffusion game. Every convergent subsequence has this property, and uniqueness of the equilibrium occupation pair therefore implies convergence of the entire sequence.
If only the equilibrium population flow is unique, Step~2 and compactness show that every subsequence has a further subsequence converging to that flow in the metric \eqref{eq:flat-norm}. This proves convergence of the entire population sequence in that metric.

\end{proof}

\subsection{Proof of Theorem~\ref{thm:convex-grid} and Corollary~\ref{cor:convex-grid-balanced}}
\label{subsec:proof-convex-grid}

\begin{proof}
Write $R_h\phi=(\phi(y_i))_i$ for restriction to the grid.
All constants below are independent of $N,K$ and the smoothing
parameters. The action interval and the horizon are fixed.
We combine the comparison of control costs used by Bertucci and Cecchin
\cite[Section~3.2]{BertucciCecchin2024} with the Mirror-Prox estimate
\eqref{eq:master}.

\emph{Step 1: restriction to the grid and estimates for the backward equation.}
Restriction satisfies the uniform bound
\begin{equation}\label{eq:grid-sampling-bound}
 \|R_h\phi\|_{H_h^m}\le C\|\phi\|_{H^m},\qquad m\ge1.
\end{equation}
To prove it, for $j\ge1$ use the repeated integral formula
\[
 D_h^j\phi(y_i)
 =h^{-j}\int_{[0,h]^j}
       \phi^{(j)}(y_i+s_1+\cdots+s_j)\,ds_1\cdots ds_j.
\]
Jensen's inequality, followed by integration in $s_1$ and summation
over the grid cells, gives
\[
 h\sum_i|D_h^j\phi(y_i)|^2
 \le h^{1-j}\int_{[0,h]^j}
          \sum_i|\phi^{(j)}(y_i+s_1+\cdots+s_j)|^2\,ds
 =\|\phi^{(j)}\|_{L^2}^2.
\]
For $j=0$, comparison with the values inside each adjacent cell gives
\[
 h\sum_i|\phi(y_i)|^2
 \le2\|\phi\|_{L^2}^2+2h^2\|\phi'\|_{L^2}^2.
\]
These estimates first hold for smooth functions and then extend by
density to $H^m$; point evaluation is continuous in $H^1$ on the
one-dimensional torus. Summing with the weights in \eqref{eq:Hm}
proves \eqref{eq:grid-sampling-bound}. Since restriction commutes
with the weak time derivative, it also proves
\begin{equation}\label{eq:grid-time-sampling}
 \|R_hv\|_{Y_h}\le C\|v\|_Y,
 \qquad Y=H^1(0,T;H^m).
\end{equation}

Taylor's formula gives, uniformly over $a\in A$,
\begin{equation}\label{eq:grid-smooth-consistency}
 \|\LL_h^aR_hv-R_h\LL^av\|_\infty
 \le Ch\|v\|_{C^3}.
\end{equation}
Indeed, the upwind first derivatives have error at most
$h\|v''\|_\infty/2$, and the centered second derivative has error
at most $Ch\|v'''\|_\infty$.

For later use let $v$ solve the backward equation
\[
 -\partial_tv-\tfrac12\sigma^2\partial_{yy}v-a^*\partial_yv=\psi,
 \qquad v(T)=\phi.
\]
The energy and gradient estimates
\eqref{eq:convex-adjoint-energy}--\eqref{eq:convex-adjoint-gradient}
apply. There is also the bound
\begin{equation}\label{eq:grid-adjoint-C3}
 \int_0^T\|v(t)\|_{C^3}\,dt
 \le C\left(\|\phi\|_{C^3}
               +\int_0^T\|\psi_t\|_{C^3}\,dt\right).
\end{equation}
To see this, reverse time and differentiate the equation $j$ times
in space for $j=0,1,2,3$. The highest-order drift term is
$a^*\partial_y(\partial_y^jv)$; the term with one derivative of
$a^*$ is a bounded multiple of $\partial_y^jv$, and all other
commutators involve lower derivatives of $v$ and derivatives of
$a^*$ of order at most three. The maximum principle, inductively
in $j$, and Gronwall's inequality give
\[
 \|v(t)\|_{C^3}
 \le C\left(\|\phi\|_{C^3}
               +\int_t^T\|\psi_s\|_{C^3}\,ds\right).
\]
Integration proves \eqref{eq:grid-adjoint-C3}. Approximation justifies
these estimates for the time-measurable forcing used below. The
constants are uniform when $T$ is replaced by any terminal time
in $[0,T]$.

Fix $0<\epsilon\le\epsilon_0$, with $\epsilon_0$ smaller than the
injectivity radius of the torus. Mollify normalized Lipschitz data
$\phi,\psi_t$ in space, where
\[
 \operatorname{Lip}(\phi)^2
       +\int_0^T\operatorname{Lip}(\psi_t)^2\,dt\le1,
 \qquad \phi(y_0)=\psi_t(y_0)=0.
\]
For the solution $v_\epsilon$ with mollified data,
\eqref{eq:grid-time-sampling}, the estimates proved in
\Cref{subsec:proof-convex-diffusion}, and
\eqref{eq:grid-adjoint-C3} yield
\begin{equation}\label{eq:grid-smoothed-adjoint}
\begin{gathered}
 \|R_hv_\epsilon\|_{Y_h}\le C\epsilon^{-m},\qquad
 \int_0^T\|v_\epsilon(t)\|_{C^3}\,dt\le C\epsilon^{-2},\\
 \|\partial_yv_\epsilon\|_{L^2_tL^\infty_y}
      +\|\partial_yv_\epsilon(0)\|_\infty\le C.
\end{gathered}
\end{equation}
The last exponent follows because a spatial derivative of order
three of a mollified Lipschitz function is bounded by
$C\epsilon^{-2}$ times its Lipschitz seminorm.

\emph{Step 2: a feasible chain comparison.}
For $0<\delta\le\delta_0$, use the kernel
\[
 \alpha^\delta_t(da\mid y)
    =\frac{\mathbf1_{\{|a-a^*(t,y)|\le\delta\}}}{2\delta}\,da
\]
at the grid points. Let $\mu^{h,\delta}$ be the law of the
continuous-time chain with the averaged rates and initial law
$\rho_h$, and set
\[
 q^{h,\delta}
   :=\bigl(\mu_T^{h,\delta},
      (\mu_t^{h,\delta}(dy)\alpha^\delta_t(da\mid y))_t\bigr).
\]
Thus $\cA_hq^{h,\delta}=\rho_h$ exactly. The conditional mean
action is $a^*(t,y)$ and its conditional variance is $\delta^2/3$.

We first prove, uniformly in $\delta$,
\begin{equation}\label{eq:grid-reference-recovery}
 \sup_{0\le t\le T}W_1(\mu_t^{h,\delta},\mu_t^*)
 \le Ch^{1/3}.
\end{equation}
Fix a terminal time $t\le T$ and a $1$-Lipschitz test $\phi$
vanishing at $y_0$. Solve the homogeneous backward equation on
$[0,t]$ with terminal value $\phi_\epsilon$, its spatial mollification.
The chain and diffusion forward identities give
\[
\begin{split}
 \int\phi_\epsilon\,d(\mu_t^{h,\delta}-\mu_t^*)
 &=\int v(0)\,d(\rho_h-\rho)\\
 &\quad+\int_0^t\!\int
    \left[\int_A\LL_h^aR_hv\,\alpha^\delta_s(da\mid y)
          -\tfrac12\sigma^2 v_{yy}-a^* v_y\right](s,y)
             \,\mu_s^{h,\delta}(dy)\,ds.
\end{split}
\]
The identity $\int a\,\alpha^\delta_t(da\mid y)=a^*(t,y)$
and consistency \eqref{eq:grid-smooth-consistency} bound the
second term by $Ch\epsilon^{-2}$. The first term is bounded by $Ch$,
since $W_1(\rho_h,\rho)\le h/2$ and the backward gradient bound
is uniform. The mollification error is at most $C\epsilon$.
Taking the supremum over $\phi$ gives
\[
 W_1(\mu_t^{h,\delta},\mu_t^*)
 \le C(\epsilon+h\epsilon^{-2}+h),
\]
with constants independent of $t$ and $\delta$. Choosing
$\epsilon=\epsilon_0h^{1/3}$ proves \eqref{eq:grid-reference-recovery}.
Consistency is uniform over actions, so the constant in
\eqref{eq:grid-reference-recovery} is independent of $\delta$.

The grid references in (G3) give
\begin{equation}\label{eq:grid-reference-entropy}
 \bD(q^{h,\delta}\|\pi_0^N)
 \le (1+T)\log(N/c)
       +T\log\!\left(\frac{a_+-a_-}{2\delta}\right)
 \le C+(1+T)\log N+T\log(1/\delta).
\end{equation}
Indeed, terminal atom masses are at most one, and the density of
the running comparison at each state relative to $\zeta$ is at
most $(a_+-a_-)/(2\delta)$. Divide by the reference lower bounds
and integrate the logarithm.

\emph{Step 3: control and forward-equation errors.}
Let $u_h^*=R_hu^*$, and define the running HJB residual on the grid
\[
 \Gamma_h^*(t,y_i,a)
 :=f(t,y_i,a,\mu_t^*)+\partial_tu^*(t,y_i)
                          +\LL_h^au_h^*(t,y_i).
\]
The terminal HJB residual is zero. Assumption (G1) gives
\[
 u^*\in C([0,T];H^{m+1})\subset C([0,T];C^3),
 \qquad \sup_N\|u_h^*\|_{Y_h}<\infty,
\]
where the second bound follows from \eqref{eq:grid-time-sampling}. Combining
\eqref{eq:grid-smooth-consistency} with the HJB residual bound for the diffusion
\eqref{eq:convex-slack}, we obtain
\begin{equation}\label{eq:grid-slack-coercivity}
 \frac{c_{f_1}}{2}|a-a^*(t,y_i)|^2-Ch
 \le\Gamma_h^*(t,y_i,a)
 \le\frac{C_{f_1}}{2}|a-a^*(t,y_i)|^2+Ch.
\end{equation}

Abbreviate $\bar\pi=\bar\pi_{N,K}$,
$r=\cA_h\bar\pi-\rho_h$, $\eta=\|r\|_{Y_h^*}$ and
$E=E_{N,K}$. From \eqref{eq:grid-slack-coercivity} and the
conditional variance of the comparison kernel,
\[
 \pair{\cC_{\bm\mu^*}-\cA_h^*u_h^*}
           {\bar\pi-q^{h,\delta}}
 \ge\frac{c_{f_1}}{2}E-C(\delta^2+h).
\]
All cost vectors here are evaluated at grid states. By (G2) and
\eqref{eq:grid-reference-recovery},
\[
 \nG{\cC_{\cM q^{h,\delta}}-\cC_{\bm\mu^*}}
 \le Ch^{1/3}.
\]
Both pairs have probability measures as their components, so the norm
of their difference is at most $2\sqrt{1+T}$. Consequently
\[
 \pair{\cC_{\cM q^{h,\delta}}}
           {\bar\pi-q^{h,\delta}}
 \ge\frac{c_{f_1}}{2}E
       +\pair{r}{u_h^*}-C(\delta^2+h^{1/3}).
\]
Insert the feasible, finite-entropy comparison $q^{h,\delta}$ into
the estimate \eqref{eq:master}. Cauchy--Schwarz, Young's
inequality and the uniform bound on $\|u_h^*\|_{Y_h}$ give
\begin{equation}\label{eq:grid-master-coercive}
 \frac{c_{f_1}}{2}E+\frac{\tau K}{4}\eta^2
 \le C\left[
       \frac{1+\log N+\log(1/\delta)}K
       +\delta^2+h^{1/3}\right].
\end{equation}
Take $\delta=\delta_0/\sqrt K$. This proves
\eqref{eq:convex-grid-energy}.

\emph{Step 4: a grid pair that need not satisfy the forward equation.}
We prove, for any grid probability pair $\pi$, that
\begin{equation}\label{eq:grid-raw-stability}
 \lVert\cM\pi-\bm\mu^*\rVert_{W_1,2}
 \le C\left[
       E(\pi)^{1/2}+\eta_h(\pi)^{1/(m+1)}+h^{1/3}
       \right],
\end{equation}
where
$E(\pi)=\int_0^T\!\int|a-a^*(t,y)|^2\,d\pi_t\,dt$ and
$\eta_h(\pi)=\|\cA_h\pi-\rho_h\|_{Y_h^*}$.
Use the normalized terminal and running Lipschitz tests from
\eqref{eq:convex-flow-duality}, and let $v_\epsilon$ solve the
backward equation with their spatial mollifications. With
$I(\phi,\psi)$ denoting their pairing against
$\cM\pi-\bm\mu^*$, subtraction of the chain and diffusion forward
identities yields
\[
\begin{split}
 I(\phi_\epsilon,\psi_\epsilon)
 &=\pair{\cA_h\pi-\rho_h}{R_hv_\epsilon}
       +\int v_\epsilon(0)\,d(\rho_h-\rho)\\
 &\quad+\int_0^T\!\!\int
       (a-a^*(t,y))\partial_yv_\epsilon(t,y)
          \,\pi_t(dy,da)\,dt\\
 &\quad+\int_0^T\!\!\int
       (\LL_h^aR_hv_\epsilon-\LL^av_\epsilon)(t,y)
          \,\pi_t(dy,da)\,dt.
\end{split}
\]
The four terms on the right are bounded, respectively, by
\[
 C\eta_h(\pi)\epsilon^{-m},\qquad Ch,\qquad
 CE(\pi)^{1/2},\qquad Ch\epsilon^{-2},
\]
using \eqref{eq:grid-smoothed-adjoint}, Cauchy--Schwarz in the
control term and \eqref{eq:grid-smooth-consistency}. The error from
removing the mollifiers is at most $C\epsilon$. Taking the supremum
over normalized Lipschitz tests therefore gives
\[
 \lVert\cM\pi-\bm\mu^*\rVert_{W_1,2}
 \le C\left[E(\pi)^{1/2}+\epsilon
       +\eta_h(\pi)\epsilon^{-m}+h\epsilon^{-2}+h\right].
\]
If $\eta_h(\pi)\le1$, choose
\[
 \epsilon=\epsilon_0\max\{\eta_h(\pi)^{1/(m+1)},h^{1/3}\}.
\]
Each smoothing term is bounded by
$C(\eta_h(\pi)^{1/(m+1)}+h^{1/3})$. If
$\eta_h(\pi)>1$, the same conclusion follows from the bounded
diameter of the population metric. This proves
\eqref{eq:grid-raw-stability}.

Apply \eqref{eq:grid-raw-stability} to $\bar\pi_{N,K}$ and use
\eqref{eq:convex-grid-energy}. Since $h\le1$, the remaining
$h^{1/3}$ is bounded by the first term on the right of
\eqref{eq:convex-grid-flow-rate}. Increasing $C$ gives that estimate.
For the state--action estimate, first couple $(y,a)$ to
$(y,a^*(t,y))$, then couple the two state laws and use the
Lipschitz constant of $y\mapsto(y,a^*(t,y))$. This gives
\[
 W_{1,\mathbb T\times A}(\pi_t,\pi_t^*)
 \le \int|a-a^*(t,y)|\,d\pi_t
       +(1+\operatorname{Lip}(a^*_t))W_1(\mu_t,\mu_t^*).
\]
Minkowski's and Cauchy--Schwarz inequalities in time bound the
distance \eqref{eq:convex-occupation-distance} by
$CE(\pi)^{1/2}+C\lVert\cM\pi-\bm\mu^*\rVert_{W_1,2}$.
Finally, for $K\ge h^{-m/3}=N^{m/3}$, monotonicity in $K$ of
$\log(eNK)/K$ and $\log(eNK)/K^2$ gives
\[
 \frac{\log(eNK)}K\le C_m h^{1/3},\qquad
 \frac{\log(eNK)}{K^2}+\frac{h^{1/3}}K
       \le C_m h^{(m+1)/3}.
\]
Substitution into \eqref{eq:convex-grid-flow-rate} proves
\Cref{cor:convex-grid-balanced}.
\end{proof}

\clearpage
\section{Logical structure of the results}\label{app:logical-structure}

\begin{figure}[!htbp]
\centering
\begin{tikzpicture}[
  x=1cm,y=1cm,
  every node/.style={font=\scriptsize},
  box/.style={draw,rounded corners=2pt,align=center,inner sep=3pt,minimum height=7mm,text width=2.45cm},
  ass/.style={box,fill=gray!30},
  res/.style={box,fill=white},
  main/.style={box,fill=gray!10,double,double distance=1pt,thick},
  ver/.style={box,fill=white,dashed},
  arr/.style={-{Latex[length=1.6mm]},semithick},
  verif/.style={arr,dashed}
]
\node[ass] (S)   at (3,0)     {\textbf{Assumption~\ref{ass:standing}}\\ standing (S1)--(S5)};
\node[res] (PDVI) at (0,-1.5) {\textbf{Proposition~\ref{prop:PDVI-MFE}}\\ PD-VI $\Rightarrow$ MFE;\\ converse with an optimal dual function};
\node[res] (Rz)   at (3,-1.5) {\textbf{Lemmas~\ref{lem:riesz},~\ref{lem:riesz-explicit}}\\ dual step is explicit};
\node[res] (El)   at (6,-1.5) {\textbf{Lemma~\ref{lem:elementary}}\\ $\nG{\cA^*u}\le L_A\tn{u}$};
\node[ass] (Lip)  at (9,-1.5) {\textbf{Assumption~\ref{ass:Lip}}\\ Lipschitz costs};
\node[ass] (LL)   at (12,-1.5){\textbf{Assumption~\ref{ass:LL}}\\ Lasry--Lions monotone};
\node[main] (Ix)  at (0,-3)    {\textbf{Theorem~\ref{thm:inexact}}\\ Lemma~\ref{lem:inexact};\\ summable errors give $O(K^{-1})$ estimates};
\node[res] (TP)  at (3,-3)    {\textbf{Lemma~\ref{lem:PD-three-point}}\\ three-point identity};
\node[res,text width=2.7cm] (OS)  at (6,-3)    {\textbf{Lemma~\ref{lem:PD-one-step}}\\ one-step estimate, $\theta\le1$};
\node[res] (Mon) at (12,-3)   {\textbf{Proposition~\ref{prop:monotone}}\\ $\cF$ monotone};
\node[res] (Top) at (0,-4.5)  {\textbf{Lemma~\ref{lem:topology}}\\ compactness, continuity};
\node[res] (EM)  at (6,-4.5)  {\textbf{Theorem~\ref{thm:PD-ergodic-minty}}\\ ergodic Minty estimate};
\node[res] (Rel) at (9,-4.5)  {\textbf{Proposition~\ref{prop:eps-relations}}\\ $\GapM$ controls $\Expl$};
\node[res] (Cor) at (6,-6)    {\textbf{Corollary~\ref{cor:primal-minty-residual}}\\ primal Minty + forward-equation residual};
\node[ass,text width=2.9cm] (PE)  at (12,-6)   {\textbf{Assumptions~\ref{ass:projection},~\ref{ass:entropy-bound}}\\ total-variation correction, entropy bound};
\node[main] (CP) at (3,-7.5)  {\textbf{Theorem~\ref{thm:cluster-points}}\\ cluster points are MFE\\ + density (\ref{ass:density})};
\node[main] (QE) at (9,-7.5)  {\textbf{Theorem~\ref{thm:quantitative-eps-MFE}}\\ corrected average is an $O(K^{-1/2})$-MFE};
\node[main] (SL) at (6,-7.5)  {\textbf{Theorem~\ref{thm:strong-LL}}\\ flow rates: strong LL (\ref{ass:strong-LL})\\ + finite-entropy PD-VI solution};
\node[ver] (FS)  at (12,-7.5) {\textbf{Theorem~\ref{thm:finite-state}}\\ finite state/action;\\ positive reference};
\node[ver] (MI)  at (12,-9)   {\textbf{Theorem~\ref{thm:mesh}}\\ periodic grid:\\ mesh-independent step sizes};
\node[ass,text width=2.8cm] (CD) at (3,-9)
  {\textbf{Assumption~\ref{ass:convex-diffusion}}\\ smooth equilibrium;\\ convex control cost};
\node[main,text width=3cm] (CR) at (3,-10.5)
  {\textbf{Theorem~\ref{thm:convex-diffusion}}\\ Wasserstein convergence rate\\ for the diffusion game};
\node[main,text width=2.8cm] (ML) at (12,-10.5)
  {\textbf{Proposition~\ref{prop:mesh-limit}}\\ chain occupations converge\\ to diffusion equilibria};
\node[main,text width=2.45cm] (CC) at (0,-9)
  {\textbf{Proposition~\ref{prop:cost-correction}}\\ cost bound \eqref{eq:cost-correction}\\ + entropy (\ref{ass:entropy-bound});\\ using Prop.~\ref{prop:eps-relations}:\\ approximate MFE\\ at rate $O(K^{-\alpha/[2(2-\alpha)]})$};
\node[ass,text width=2.8cm] (CG) at (3,-12.6)
  {\textbf{Assumption~\ref{ass:convex-grid}}\\ smooth equilibrium;\\ convex control cost; $W_1$-Lipschitz coupling};
\node[main,text width=3cm] (QR) at (9,-12.6)
  {\textbf{Theorem~\ref{thm:convex-grid}, Corollary~\ref{cor:convex-grid-balanced}}\\ error in terms of $h$ and $K$;\\ $K\ge h^{-m/3}$ gives $O(h^{1/6})$};
\draw[arr] (S) -- (PDVI);
\draw[arr] (S) -- (Rz);
\draw[arr] (S) -- (El);
\draw[arr] (Rz) -- (TP);
\draw[arr] (Rz) -- (Ix);
\draw[arr] (El) -- (OS);
\draw[arr] (TP) -- (OS);
\draw[arr] (Lip) -- (OS);
\draw[arr] (LL) -- (Mon);
\draw[arr] (OS) -- (EM);
\draw[arr] (OS.west) -- (4.5,-3) -- (4.5,-7.5) -- (SL.west);
\draw[arr] (Mon) -- (EM);
\draw[arr] (EM) -- (Cor);
\draw[arr] (Top) -- (CP);
\draw[arr] (Cor) -- (CP);
\draw[arr] (Cor) -- (QE);
\draw[arr] (Rel) -- (QE);
\draw[arr] (PE) -- (QE);
\draw[verif] (FS) -- (PE);
\draw[verif] (MI) -- (FS);
\draw[arr] (CD) -- (CR);
\draw[arr] (Cor.east) -- (7.5,-6) -- (7.5,-10.5) -- (CR.east);
\draw[arr] (MI) -- (ML);
\draw[arr] (Cor.east) -- (7.5,-6) -- (7.5,-11.6) -- (12,-11.6) -- (ML.south);
\draw[arr] (CG) -- (QR);
\draw[arr] (MI.east) -- (13.7,-9) -- (13.7,-12.6) -- (QR.east);
\draw[arr] (Cor.west) -- (0,-6) -- (CC.north);
\end{tikzpicture}
\caption{Logical structure of the paper. Shaded boxes: assumptions, placed
next to the results that use them; the standing assumptions underlie
everything. White boxes: results; double-framed boxes: the main results.
Dashed boxes: verification of the hypotheses for finite-state games and for
the Markov-chain approximation of diffusions (dashed arrows: ``verifies'').
The total-variation correction bound in Assumption~\ref{ass:projection}
is used in Theorem~\ref{thm:quantitative-eps-MFE};
Proposition~\ref{prop:cost-correction} combines
Corollary~\ref{cor:primal-minty-residual} and Proposition~\ref{prop:eps-relations}
under the weaker cost bound \eqref{eq:cost-correction} and the entropy
bound in Assumption~\ref{ass:entropy-bound}. For $\alpha=1$, it gives
the same $O(K^{-1/2})$ approximate-equilibrium rate.
The diffusion results use occupation-measure approximations
and backward parabolic estimates. Theorem~\ref{thm:inexact} combines
the energy identity \eqref{eq:energy-identity} with the error estimates
in Lemma~\ref{lem:inexact}, under the existence of a finite-entropy
primal--dual solution and summable correction errors.}
\label{fig:map}
\end{figure}
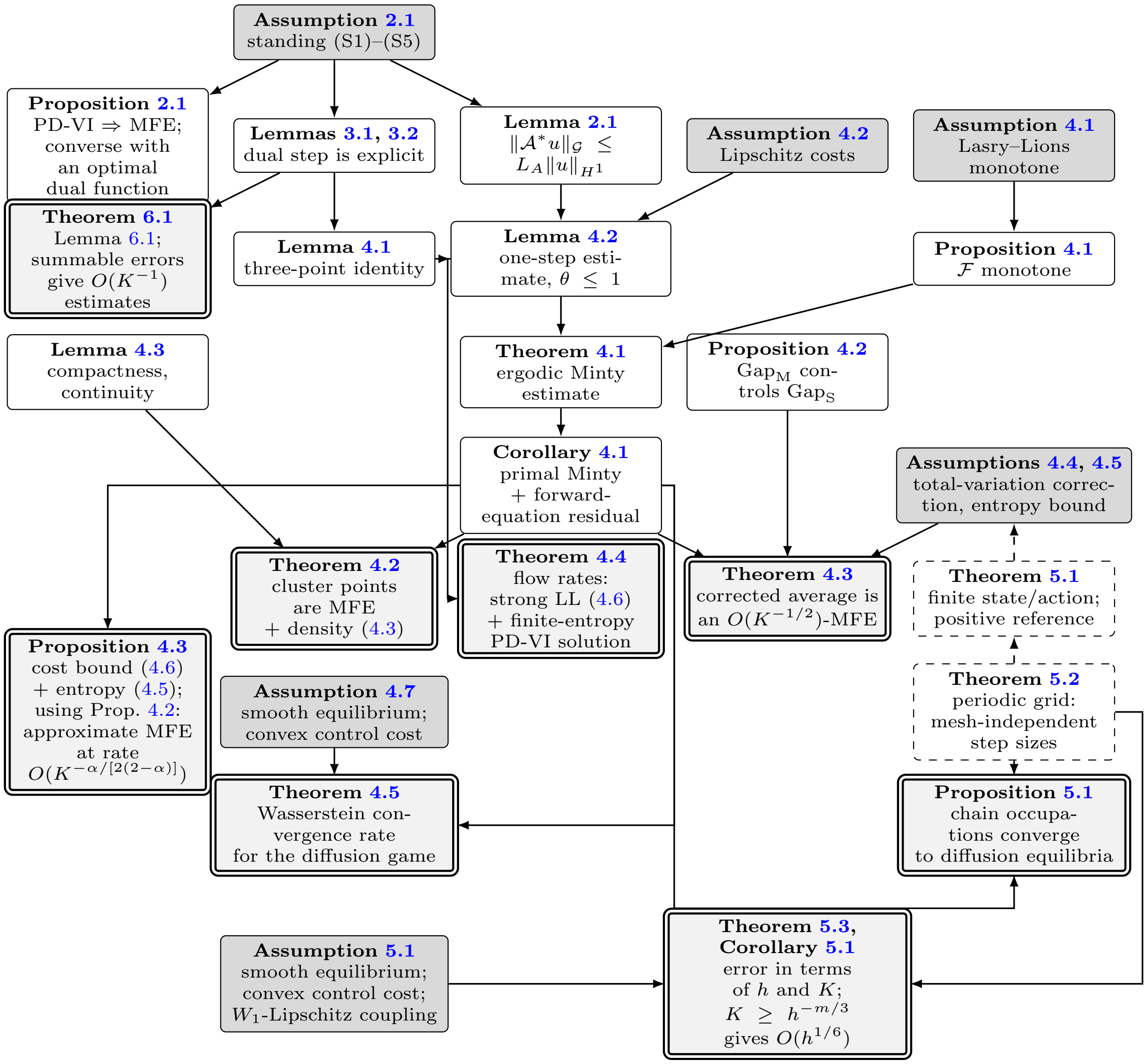

\clearpage
\bibliographystyle{plain}
\bibliography{references}

@article{korpelevich1976extragradient,
  title={The extragradient method for finding saddle points and other problems},
  author={Korpelevich, Galina M},
  journal={Matecon},
  volume={12},
  pages={747--756},
  year={1976}
}

@article{nemirovski2004prox,
  title={Prox-method with rate of convergence {$O(1/t)$} for variational inequalities with {Lipschitz} continuous monotone operators and smooth convex-concave saddle point problems},
  author={Nemirovski, Arkadi},
  journal={SIAM Journal on Optimization},
  volume={15},
  number={1},
  pages={229--251},
  year={2004},
  publisher={SIAM}
}

@article{lasry2007mean,
  title={Mean field games},
  author={Lasry, Jean-Michel and Lions, Pierre-Louis},
  journal={Japanese Journal of Mathematics},
  volume={2},
  number={1},
  pages={229--260},
  year={2007},
  publisher={Springer}
}

@article{huang2006large,
  title={Large population stochastic dynamic games: closed-loop {McKean--Vlasov} systems and the {Nash} certainty equivalence principle},
  author={Huang, Minyi and Malham{\'e}, Roland P and Caines, Peter E},
  journal={Communications in Information and Systems},
  volume={6},
  number={3},
  pages={221--252},
  year={2006},
  doi={10.4310/CIS.2006.v6.n3.a5}
}

@article{lauriere2022learning,
  title={Learning in mean field games: A survey},
  author={Lauri{\`e}re, Mathieu and Perrin, Sarah and P{\'e}rolat, Julien and Girgin, Sertan and Muller, Paul and {\'E}lie, Romuald and Geist, Matthieu and Pietquin, Olivier},
  journal={arXiv preprint arXiv:2205.12944},
  year={2022}
}

@incollection{lauriere2021numerical,
  title={Numerical methods for mean field games and mean field type control},
  author={Lauri{\`e}re, Mathieu},
  booktitle={Mean Field Games},
  editor={Delarue, Fran{\c{c}}ois},
  series={Proceedings of Symposia in Applied Mathematics},
  volume={78},
  pages={221--282},
  publisher={American Mathematical Society},
  address={Providence, RI},
  year={2021}
}

@incollection{achdou2020mean,
  title={Mean field games and applications: Numerical aspects},
  author={Achdou, Yves and Lauri{\`e}re, Mathieu},
  booktitle={Mean Field Games: Cetraro, Italy 2019},
  editor={Cardaliaguet, Pierre and Porretta, Alessio},
  series={Lecture Notes in Mathematics},
  volume={2281},
  pages={249--307},
  year={2020},
  publisher={Springer}
}

@article{achdou2010mean,
  title={Mean field games: numerical methods},
  author={Achdou, Yves and Capuzzo-Dolcetta, Italo},
  journal={SIAM Journal on Numerical Analysis},
  volume={48},
  number={3},
  pages={1136--1162},
  year={2010},
  publisher={SIAM}
}

@article{achdou2013mean,
  title={Mean field games: convergence of a finite difference method},
  author={Achdou, Yves and Camilli, Fabio and Capuzzo-Dolcetta, Italo},
  journal={SIAM Journal on Numerical Analysis},
  volume={51},
  number={5},
  pages={2585--2612},
  year={2013},
  publisher={SIAM}
}

@article{cacace2021policy,
  title={A policy iteration method for mean field games},
  author={Cacace, Simone and Camilli, Fabio and Goffi, Alessandro},
  journal={ESAIM: Control, Optimisation and Calculus of Variations},
  volume={27},
  pages={85},
  year={2021},
  doi={10.1051/cocv/2021081},
  publisher={EDP Sciences}
}

@article{briceno2018proximal,
  title={Proximal methods for stationary mean field games with local couplings},
  author={Brice{\~n}o-Arias, Luis M and Kalise, Dante and Silva, Francisco J},
  journal={SIAM Journal on Control and Optimization},
  volume={56},
  number={2},
  pages={801--836},
  year={2018},
  publisher={SIAM}
}

@article{liu2021computational,
  title={Computational methods for first-order nonlocal mean field games with applications},
  author={Liu, Siting and Jacobs, Matthew and Li, Wuchen and Nurbekyan, Levon and Osher, Stanley J},
  journal={SIAM Journal on Numerical Analysis},
  volume={59},
  number={5},
  pages={2639--2668},
  year={2021},
  doi={10.1137/20M1334668},
  publisher={SIAM}
}

@article{angiuli2019cemracs,
  title={{CEMRACS} 2017: numerical probabilistic approach to {MFG}},
  author={Angiuli, Andrea and Graves, Christy V and Li, Houzhi and Chassagneux, Jean-Fran{\c{c}}ois and Delarue, Fran{\c{c}}ois and Carmona, Ren{\'e}},
  journal={ESAIM: Proceedings and Surveys},
  volume={65},
  pages={84--113},
  year={2019},
  doi={10.1051/proc/201965084},
  publisher={EDP Sciences}
}

@article{cardaliaguet2017learning,
  title={Learning in mean field games: the fictitious play},
  author={Cardaliaguet, Pierre and Hadikhanloo, Saeed},
  journal={ESAIM: Control, Optimisation and Calculus of Variations},
  volume={23},
  number={2},
  pages={569--591},
  year={2017},
  publisher={EDP Sciences}
}

@article{lavigne2023generalized,
  title={Generalized conditional gradient and learning in potential mean field games},
  author={Lavigne, Pierre and Pfeiffer, Laurent},
  journal={Applied Mathematics \& Optimization},
  volume={88},
  pages={89},
  year={2023},
  doi={10.1007/s00245-023-10056-8},
  publisher={Springer}
}

@article{BonnansLiuPfeiffer2023Theta,
  title={Error estimates of a theta-scheme for second-order mean field games},
  author={Bonnans, J. Fr\'ed\'eric and Liu, Kang and Pfeiffer, Laurent},
  journal={ESAIM: Mathematical Modelling and Numerical Analysis},
  volume={57},
  number={4},
  pages={2493--2528},
  year={2023},
  doi={10.1051/m2an/2023059}
}

@article{BonnansLavignePfeiffer2023Discrete,
  title={Discrete potential mean field games: duality and numerical resolution},
  author={Bonnans, J. Fr\'ed\'eric and Lavigne, Pierre and Pfeiffer, Laurent},
  journal={Mathematical Programming},
  volume={202},
  number={1--2},
  pages={241--278},
  year={2023},
  doi={10.1007/s10107-023-01934-8}
}

@article{LiuPfeiffer2025Mesh,
  title={A mesh-independent method for second-order potential mean field games},
  author={Liu, Kang and Pfeiffer, Laurent},
  journal={IMA Journal of Numerical Analysis},
  volume={45},
  number={2},
  pages={1226--1266},
  year={2025},
  doi={10.1093/imanum/drae061}
}

@article{hadikhanloo2017learning,
  title={Learning in anonymous nonatomic games with applications to first-order mean field games},
  author={Hadikhanloo, Saeed},
  journal={arXiv preprint arXiv:1704.00378},
  year={2017}
}

@article{perrin2020fictitious,
  title={Fictitious play for mean field games: Continuous time analysis and applications},
  author={Perrin, Sarah and P{\'e}rolat, Julien and Lauri{\`e}re, Mathieu and Geist, Matthieu and Elie, Romuald and Pietquin, Olivier},
  journal={Advances in Neural Information Processing Systems},
  volume={33},
  pages={13199--13213},
  year={2020}
}

@inproceedings{xie2021learning,
  title={Learning While Playing in Mean-Field Games: Convergence and Optimality},
  author={Xie, Qiaomin and Yang, Zhuoran and Wang, Zhaoran and Minca, Andreea},
  booktitle={Proceedings of the 38th International Conference on Machine Learning},
  series={Proceedings of Machine Learning Research},
  volume={139},
  pages={11436--11447},
  year={2021},
  publisher={PMLR}
}

@article{yu2024convergence,
  title={Convergence Analysis and Acceleration of Fictitious Play for General Mean-Field Games via the Best Response},
  author={Yu, Jiajia and Cheng, Xiuyuan and Liu, Jian-Guo and Zhao, Hongkai},
  journal={arXiv preprint arXiv:2411.07989},
  year={2024}
}

@article{hadikhanloo2022learning,
author  = {Hadikhanloo, Saeed and Laraki, Rida and Mertikopoulos, Panayotis and Sorin, Sylvain},
title   = {Learning in Nonatomic Games, {Part I}: Finite Action Spaces and Population Games},
journal = {Journal of Dynamics and Games},
volume  = {9},
number  = {4},
pages   = {433--460},
year    = {2022},
doi     = {10.3934/jdg.2022018}
}

@inproceedings{perolat2022scaling,
author    = {P{\'e}rolat, Julien and Perrin, Sarah and Elie, Romuald and Lauri{\`e}re, Mathieu and Piliouras, Georgios and Geist, Matthieu and Tuyls, Karl and Pietquin, Olivier},
title     = {Scaling Mean Field Games by Online Mirror Descent},
booktitle = {Proceedings of the 21st International Conference on Autonomous Agents and Multiagent Systems},
year      = {2022},
pages     = {1028--1037}
}

@article{guo2025continuous,
  title={Continuous-time mean field games: a primal-dual characterization},
  author={Guo, Xin and Hu, Anran and Zhang, Jiacheng and Zhang, Yufei},
  journal={arXiv preprint arXiv:2503.01042v3},
  year={2025}
}

@unpublished{BayraktarHuZhang2026StochasticDuality,
  title={Characterization of Mean Field Equilibria via Stochastic Duality},
  author={Bayraktar, Erhan and Hu, Anran and Zhang, Yufei},
  note={Preprint},
  year={2026}
}

@article{BayraktarBudhirajaCohen2019RateControl,
  title={Rate Control under Heavy Traffic with Strategic Servers},
  author={Bayraktar, Erhan and Budhiraja, Amarjit and Cohen, Asaf},
  journal={The Annals of Applied Probability},
  volume={29},
  number={1},
  pages={1--35},
  year={2019},
  doi={10.1214/17-AAP1349},
  eprint={1605.09010},
  archivePrefix={arXiv},
  primaryClass={math.PR}
}

@article{BayraktarBudhirajaCohen2018Numerical,
  title={A Numerical Scheme for a Mean Field Game in Some Queueing Systems Based on {Markov} Chain Approximation Method},
  author={Bayraktar, Erhan and Budhiraja, Amarjit and Cohen, Asaf},
  journal={SIAM Journal on Control and Optimization},
  volume={56},
  number={6},
  pages={4017--4044},
  year={2018},
  doi={10.1137/17M1154357},
  eprint={1708.08343},
  archivePrefix={arXiv},
  primaryClass={math.OC}
}

@article{ChassagneuxCrisanDelarue2019Numerical,
  title={Numerical method for {FBSDEs} of {McKean--Vlasov} type},
  author={Chassagneux, Jean-Fran{\c{c}}ois and Crisan, Dan and Delarue, Fran{\c{c}}ois},
  journal={The Annals of Applied Probability},
  volume={29},
  number={3},
  pages={1640--1684},
  year={2019},
  doi={10.1214/18-AAP1429},
  eprint={1703.02007},
  archivePrefix={arXiv},
  primaryClass={math.PR}
}

@article{BayraktarCohen2018FiniteState,
  title={Analysis of a Finite State Many Player Game Using Its Master Equation},
  author={Bayraktar, Erhan and Cohen, Asaf},
  journal={SIAM Journal on Control and Optimization},
  volume={56},
  number={5},
  pages={3538--3568},
  year={2018},
  doi={10.1137/17M113887X},
  eprint={1707.02648},
  archivePrefix={arXiv},
  primaryClass={math.OC}
}

@article{guo2019learning,
  title={Learning mean-field games},
  author={Guo, Xin and Hu, Anran and Xu, Renyuan and Zhang, Junzi},
  journal={Advances in Neural Information Processing Systems},
  volume={32},
  pages={4966--4976},
  year={2019}
}

@inproceedings{lauriere2022scalable,
  title={Scalable deep reinforcement learning algorithms for mean field games},
  author={Lauri{\`e}re, Mathieu and Perrin, Sarah and Girgin, Sertan and Muller, Paul and Jain, Ayush and Cabannes, Theophile and Piliouras, Georgios and P{\'e}rolat, Julien and Elie, Romuald and Pietquin, Olivier and Geist, Matthieu},
  booktitle={Proceedings of the 39th International Conference on Machine Learning},
  series={Proceedings of Machine Learning Research},
  volume={162},
  pages={12078--12095},
  year={2022},
  publisher={PMLR}
}

@article{magnino2026solving,
  title={Solving continuous mean field games: Deep reinforcement learning for non-stationary dynamics},
  author={Magnino, Lorenzo and Shao, Kai and Wu, Zida and Shen, Jiacheng and Lauri{\`e}re, Mathieu},
  journal={Advances in Neural Information Processing Systems},
  volume={38},
  pages={104325--104354},
  year={2025}
}

@article{hynd2023evolution,
  title={Evolution of mixed strategies in monotone games},
  author={Hynd, Ryan},
  journal={SIAM Journal on Optimization},
  volume={33},
  number={4},
  pages={2750--2771},
  year={2023},
  doi={10.1137/22M1486066},
  publisher={SIAM}
}

@article{yardim2025variational,
  title={A Variational Inequality Approach to Independent Learning in Static Mean-Field Games},
  author={Yardim, Batuhan and Cayci, Semih and He, Niao},
  journal={ACM/IMS Journal of Data Science},
  volume={2},
  number={2},
  pages={1--43},
  year={2025},
  publisher={ACM, New York, NY}
}

@article{tangpi2025particle,
  title={Particle system approximation of {Nash} equilibria in large games},
  author={Tangpi, Ludovic and Touzi, Nizar},
  journal={arXiv preprint arXiv:2510.19211},
  year={2025}
}

@article{isobe2026last,
  title={Last iterate convergence in monotone mean field games},
  author={Isobe, Noboru and Abe, Kenshi and Ariu, Kaito},
  journal={Advances in Neural Information Processing Systems},
  volume={38},
  pages={18086--18119},
  year={2025}
}

@article{hadikhanloo2019finite,
  title={Finite mean field games: fictitious play and convergence to a first order continuous mean field game},
  author={Hadikhanloo, Saeed and Silva, Francisco J},
  journal={Journal de Math{\'e}matiques Pures et Appliqu{\'e}es},
  volume={132},
  pages={369--397},
  year={2019},
  publisher={Elsevier}
}

@article{meynard2026extragradient,
  title={Extragradient methods for mean field games of controls and mean field type {FBSDEs}},
  author={Meynard, Charles},
  journal={arXiv preprint arXiv:2602.14621},
  year={2026}
}

@article{abdulaziz2026bregman,
  title={Bregman-projected mirror methods for regularized stationary mean-field games},
  author={Abdulaziz, Hussain Al and Ashrafyan, Yuri and Gevorgyan, Yeva and Gomes, Diogo},
  journal={arXiv preprint arXiv:2606.19611},
  year={2026}
}

@article{guo2022entropy,
  title={Entropy regularization for mean field games with learning},
  author={Guo, Xin and Xu, Renyuan and Zariphopoulou, Thaleia},
  journal={Mathematics of Operations Research},
  volume={47},
  number={4},
  pages={3239--3260},
  year={2022},
  doi={10.1287/moor.2021.1238},
  publisher={INFORMS}
}

@inproceedings{cui2021approximately,
  title={Approximately solving mean field games via entropy-regularized deep reinforcement learning},
  author={Cui, Kai and Koeppl, Heinz},
  booktitle={Proceedings of the 24th International Conference on Artificial Intelligence and Statistics},
  series={Proceedings of Machine Learning Research},
  volume={130},
  pages={1909--1917},
  year={2021},
  publisher={PMLR}
}

@article{bigi2023approximate,
  title={Approximate variational inequalities and equilibria},
  author={Bigi, Giancarlo and Lampariello, Lorenzo and Sagratella, Simone and Sasso, Valerio Giuseppe},
  journal={Computational Management Science},
  volume={20},
  number={1},
  pages={43},
  year={2023},
  doi={10.1007/s10287-023-00476-w},
  publisher={Springer}
}

@article{GomesMohrSouza2013,
  title={Continuous time finite state mean field games},
  author={Gomes, Diogo A. and Mohr, Joana and Souza, Rafael Rigao},
  journal={Applied Mathematics \& Optimization},
  volume={68},
  number={1},
  pages={99--143},
  year={2013},
  publisher={Springer}
}

@article{Nesterov2007,
  title={Dual extrapolation and its applications to solving variational inequalities and related problems},
  author={Nesterov, Yurii},
  journal={Mathematical Programming},
  volume={109},
  number={2--3},
  pages={319--344},
  year={2007},
  publisher={Springer}
}

@book{FonsecaLeoni2007,
  title={Modern Methods in the Calculus of Variations: $L^p$ Spaces},
  author={Fonseca, Irene and Leoni, Giovanni},
  publisher={Springer},
  address={New York},
  year={2007}
}

@article{KurtzStockbridge1998,
  title={Existence of {M}arkov controls and characterization of optimal {M}arkov controls},
  author={Kurtz, Thomas G. and Stockbridge, Richard H.},
  journal={SIAM Journal on Control and Optimization},
  volume={36},
  number={2},
  pages={609--653},
  year={1998}
}

@article{BhattBorkar1996,
  title={Occupation measures for controlled {M}arkov processes: characterization and optimality},
  author={Bhatt, Abhay G. and Borkar, Vivek S.},
  journal={The Annals of Probability},
  volume={24},
  number={3},
  pages={1531--1562},
  year={1996}
}

@book{CarmonaDelarue2018,
  title={Probabilistic Theory of Mean Field Games with Applications {I}--{II}},
  author={Carmona, Ren{\'e} and Delarue, Fran{\c{c}}ois},
  series={Probability Theory and Stochastic Modelling},
  volume={83--84},
  publisher={Springer},
  year={2018}
}

@article{ChambollePock2011,
  title={A first-order primal-dual algorithm for convex problems with applications to imaging},
  author={Chambolle, Antonin and Pock, Thomas},
  journal={Journal of Mathematical Imaging and Vision},
  volume={40},
  number={1},
  pages={120--145},
  year={2011},
  publisher={Springer}
}

@article{Trevisan2016,
  title={Well-posedness of multidimensional diffusion processes with weakly differentiable coefficients},
  author={Trevisan, Dario},
  journal={Electronic Journal of Probability},
  volume={21},
  pages={1--41},
  year={2016},
  doi={10.1214/16-EJP4453}
}

@unpublished{KLMP,
  title={{Kullback--Leibler Mirror-Prox} for measure-valued variational inequalities and mean-field equilibria},
  author={Bayraktar, Erhan and Ekren, Ibrahim and Vy, Lu and Zhang, Ziqing},
  note={\href{https://arxiv.org/abs/2608.10293}{arXiv:2608.10293}},
  year={2026}
}

@book{KushnerDupuis2001,
  title={Numerical Methods for Stochastic Control Problems in Continuous Time},
  author={Kushner, Harold J. and Dupuis, Paul},
  edition={2},
  year={2001},
  publisher={Springer},
  address={New York}
}

@misc{GevorgyanGomes2026,
  title={A {R}iesz-map preconditioned extragradient method for price formation in mean field games},
  author={Gevorgyan, Yeva and Gomes, Diogo},
  year={2026},
  note={arXiv:2610.13934}
}

@misc{MengLiuLiOsher2024,
  title={A primal-dual hybrid gradient method for solving optimal control problems and the corresponding {H}amilton--{J}acobi {PDE}s},
  author={Meng, Tingwei and Liu, Siting and Li, Wuchen and Osher, Stanley},
  year={2024},
  note={arXiv:2403.02468}
}

@article{BertucciCecchin2024,
  author={Bertucci, Charles and Cecchin, Alekos},
  title={Mean Field Games Master Equations: From Discrete to Continuous State Space},
  journal={SIAM Journal on Mathematical Analysis},
  volume={56},
  number={2},
  pages={2569--2610},
  year={2024},
  doi={10.1137/23M1552528}
}

@misc{bayraktar2026policygradient,
  author={Bayraktar, Erhan and Hernandez, Martin and Yan, Qinxin and Zhu, Yuhua},
  title={Policy Gradient for Continuous-Time Mean-Field Control},
  year={2026},
  note={arXiv:2605.20718}
}

@article{DumitrescuLeutscherTankov2023,
  author  = {Dumitrescu, Roxana and Leutscher, Marcos and Tankov, Peter},
  title   = {Linear programming fictitious play algorithm for mean field games with optimal stopping and absorption},
  journal = {{ESAIM}: Mathematical Modelling and Numerical Analysis},
  year    = {2023},
  volume  = {57},
  number  = {2},
  pages   = {953--990},
  doi     = {10.1051/m2an/2023019}
}

\end{document}